\documentclass[10pt,journal,compsoc]{IEEEtran}

\usepackage{amsmath,bm, amsfonts,color,colordvi,caption,subcaption,url}
\usepackage{graphicx}
\usepackage{mathrsfs}
\usepackage{CJK}
\usepackage{amsmath}
\usepackage{makecell}
\usepackage{picinpar}
\usepackage{multirow}
\usepackage{picinpar}
\usepackage{algorithm}
\usepackage{algorithmic}
\usepackage{ragged2e}
\usepackage{array}
\usepackage{bm}

\usepackage[colorlinks,
           urlcolor=black,
             linkcolor=black,
            anchorcolor=black,
            citecolor=black
            ]{hyperref}
\usepackage[noadjust]{cite}
\usepackage{url}
\usepackage{lipsum}
\usepackage{float}
\newcounter{RomanNumber}

\makeatletter

\newcommand{\Rmnum}[1]{\expandafter\@slowromancap\romannumeral #1@}
\makeatother
\graphicspath{ {./images/} }

\newtheorem{theorem}{Theorem}[section]

\newtheorem{definition}{Definition}[section]
\newtheorem{remark}{Remark}[section]

\newtheorem{example}{Example}[section]
\newcommand{\be}{\begin{equation}}
\newcommand{\ee}{\end{equation}}
\newcommand{\ba}{\begin{eqnarray}}
\newcommand{\ea}{\end{eqnarray}}
\newcommand{\bas}{\begin{eqnarray*}}
	\newcommand{\eas}{\end{eqnarray*}}

\def\I{{\cal I}}

\def\N{{\mathbb N}}

\def\R{{\mathbb R}}

\def\bfe{{\bf 1}}
\def\bfe{{\bf 1}}

\def\bfe{{\bf 1}}

\def\bfv{{\bf v}}
\def\bfx{{\bf x}}
\def\bfw{{\bf w}}
\def\bfy{{\bf y}}
\def\bfz{{\bf z}}

\def\bfu{{\bf u}}
\def\bfr{{\bf r}}

\def\hsvm{\texttt{HSVM}}
\def\ssvm{\texttt{LSVM}}
\def\sosvm{\texttt{SSVM}}
\def\psvm{\texttt{PSVM}}
\def\rsvm{\texttt{RSVM}}
\def\logi{\texttt{LOGI}}

 \def\TNI{\texttt{TNI}}
\def\ACC{\texttt{ACC}}
\def\SWS{\texttt{SWS/ITER}}

\def\CPU{\texttt{CPU}}
\def\NSV{\texttt{NSV}}
\def\pega{\texttt{PEGA}}
\def\svrg{\texttt{SVRG}}
\def\katy{\texttt{KATY}}

\def\ADMML{{\texttt{L$_{0/1}$ADMM}}}

\hypersetup{colorlinks,breaklinks,
            citecolor=black,
            linkcolor=black}
\usepackage[nameinlink,noabbrev,capitalize]{cleveref}

\ifCLASSINFOpdf
\else
\fi

\ifCLASSOPTIONcompsoc
\else
  \usepackage{cite}
\fi

\begin{document}
\title{Support Vector Machine Classifier via ~~~~ \\ $L_{0/1}$ Soft-Margin Loss}

\author{
{Huajun Wang, Yuanhai Shao$^*$, Shenglong Zhou, Ce Zhang and  Naihua Xiu$^*$}
\IEEEcompsocitemizethanks{\IEEEcompsocthanksitem
H.J. Wang, C. Zhang, N.H. Xiu  are with the Department of Applied Mathematics, Beijing Jiaotong
University, Beijing, P.R. China.
Email: {huajunwang@bjtu.edu.cn}, czhang@bjtu.edu.cn, nhxiu@bjtu.edu.cn.\protect
\IEEEcompsocthanksitem
Y.H. Shao is with the School of Management, Hainan University, Haikou, P.R. China. Email:  shaoyuanhai@hainanu.edu.cn.\protect
\IEEEcompsocthanksitem S.L. Zhou is with  the School of Mathematical Sciences, University of Southampton, Southampton, UK. Email: shenglong.zhou@soton.ac.uk.\protect
\IEEEcompsocthanksitem * Corresponding author
}

\thanks{Manuscript received xx, xx; revised xx, xx.}}

\IEEEtitleabstractindextext{\justify
\begin{abstract}
Support vector machines (SVM) have drawn wide attention for the last two decades due to its extensive applications, so a vast body of work  has developed optimization algorithms to solve SVM with various soft-margin losses. To distinguish all, in this paper, we aim at solving an ideal soft-margin loss SVM: $L_{0/1}$ soft-margin loss SVM (dubbed as $L_{0/1}$-SVM). Many of the existing (non)convex soft-margin losses can be viewed as one of the surrogates of the $L_{0/1}$ soft-margin loss. Despite its discrete nature, we manage to establish the optimality theory for the $L_{0/1}$-SVM including the existence of the optimal solutions, the relationship between them and P-stationary points. These not only enable us to deliver a  rigorous definition of  $L_{0/1}$ support vectors but also allow  us to define a working set. Integrating such a working set, a fast alternating direction method of multipliers is then proposed with its limit point being a locally optimal solution to the  $L_{0/1}$-SVM.  Finally, numerical experiments demonstrate that our proposed method outperforms some leading classification solvers from SVM communities, in terms of faster computational speed and a fewer number of support vectors. The bigger the data size is, the more evident its advantage appears.
\end{abstract}

}

\maketitle


\IEEEdisplaynontitleabstractindextext

%
\IEEEpeerreviewmaketitle

\ifCLASSOPTIONcompsoc
\IEEEraisesectionheading{\section{Introduction}\label{sec:introduction}}
\else
\section{Introduction}
\label{sec:introduction}
\fi

%
%
%
%
\IEEEPARstart{S}{upport} vector machines (SVM) were first introduced by Vapnik and Cortes \cite{CV95} and then have been
extensively applied into machine learning, statistic, pattern recognition and so forth. The basic idea of SVM is to find a maximum margin-type hyperplane in the input space that separates the training dataset. In the paper, we focus on the binary classification problem described as follows. Suppose we are given a training set $\{(\bfx_i,y_i):i=1,2,\cdots,m\},$ where $\bfx_i\in \R^n$ are the input vectors and $y_i\in \{-1,1\}$ are the output labels. The purpose of SVM is to train a hyperplane $\langle  \bfw, \bfx\rangle+b= w_1x_1+\cdots+w_nx_n+b=0$ with $ \bfw \in \R^n$ and $b\in \R$ to be estimated by the  training set.  For any new input vector $\bfx'$, one can predict its  label $y'$ by $y'=1$ if $ \langle \bfw,\bfx'\rangle+b>0$ and $y'=-1$ otherwise. In order to find an optimal hyperplane, there are two possible scenarios: linearly separable and inseparable training data. If the training data is linearly separated in the input space, then the unique optimal hyperplane can be obtained by solving a convex quadratic programming:
\begin{eqnarray}\label{HM-SVM}
&\underset{ \bfw \in \R^n,b\in\R}{\min} & \frac{1}{2} \Vert \bfw \Vert^{2} \nonumber\\
&\mbox{s.t.}  & y_{i}(\langle \bfw , \bfx_i\rangle+b)\geq1,{ i\in\N_m,}
\end{eqnarray}
where $\N_m:=\{1,2,\cdots,m\}$. 
The above model is known as the hard-margin SVM because it requires correct classifications of all training samples. When it comes to the training data being linearly inseparable in the input space, the popular approach is to allow violations in the satisfaction of the constraints in (\ref{HM-SVM}) and penalize such violations in the objective function, namely,
\begin{eqnarray}\label{SM-SVM1}
\min_{ \bfw \in \R^n,b\in \R}~~\frac{1}{2}\Vert \bfw \Vert^{2}+C\sum_{i=1}^{m}\ell(1-y_{i}f(\bfx_i)),
\end{eqnarray}
where $C>0$ is a penalty parameter and $f(\bfx_i):=\langle \bfw , \bfx_i\rangle+b$. Here, $\ell(\cdot)$ is one of loss functions that aims at penalizing  some sufficiently incorrectly classified samples and leaving the others. The above model is known as soft-margin SVM, allowing misclassified training samples. Authors in \cite{CV95,JP2011,YY2016} have pointed out that the ideal soft-margin SVM is
 \begin{eqnarray}\label{SM-SVM}
\min_{ \bfw \in \R^n,b\in \R}~~\frac{1}{2}\Vert \bfw \Vert^{2}+C\sum_{i=1}^{m}\ell_{0/1}(1-y_{i}f(\bfx_i)),
\end{eqnarray}
where the soft-margin loss function $\ell_{0/1}(\cdot)$ is given by
\begin{eqnarray}\label{l01-loss}
\ell_{0/1}(t_i)=
\begin{cases}
1,& 1-t_i>0,\\
0,& 1-t_i\leq0,
\end{cases}
\end{eqnarray}
and $t_i = y_if(\bfx_i), i\in \N_m$.
We name \eqref{SM-SVM} as $L_{0/1}$-SVM, which minimizes the number of soft-margin misclassified samples. It is worth mentioning that the  $\ell_{0/1}(\cdot)$ loss function  arises in binary-valued regression, and is  useful in many machine learning problems: candidates include those from perceptron learning \cite{LL2007}, deep learning \cite{IB2016} {and distributionally robust supervised learning \cite{WG2018}.}
However, the $L_{0/1}$-SVM is NP-hard \cite{BK1995,EV1998}  since the $\ell_{0/1}(\cdot)$ loss is  nonconvex and  discontinuous, and up to now, it has not been fundamentally well investigated.

As far as we know, this is the first paper that establishes the optimality theory for the $L_{0/1}$-SVM and develops an effective algorithm aiming at pursuing an optimal solution to (\ref{SM-SVM}). The main contributions are summarized as follows.

(C1) We prove that the globally optimal solutions to the $L_{0/1}$-SVM exist and also establish its optimality condition aiming at finding such solutions. The condition has a close relationship to the P-stationary point which is very practical to solve the $L_{0/1}$-SVM, even though the problem is NP-hard.

(C2) Recall that the vector $\bfw^*$
that maximizes the margin can be shown to have the form:
\begin{eqnarray}\label{supp-vects}
\bfw^*= \alpha_1^* y_1 \bfx_1+\cdots+\alpha_m^* y_m \bfx_m=\sum_{i:~\alpha_i^*\neq0}\alpha_i^* y_i \bfx_i,
\end{eqnarray}
where $ \bm{\alpha}^*=(\bm{\alpha}^*_1,\bm{\alpha}^*_2,...,\bm{\alpha}^*_m)^\top$ is a solution to the dual problem of \eqref{HM-SVM}. The training vectors $\bfx_i$ corresponding to non-zero $\alpha_i^*$ are called support vectors \cite{CV95}, \cite{CSS13}. In this paper, the P-stationary point allows us to define the $L_{0/1}$ support vectors which coincide with the non-zero elements of the Lagrangian multiplier of (\ref{SM-SVM}). From the point of the optimization, the Lagrangian multiplier can be treated as a solution to the dual problem of (\ref{SM-SVM}), even though the dual problem is difficult to be derived due to the discreteness of $\ell_{0/1}(\cdot)$. Therefore, $L_{0/1}$ support vectors are standard support vectors. Furthermore,   we show that all $L_{0/1}$ support vectors fall into the support hyperplanes $\langle\bfw^*,\bfx\rangle+b^*=\pm1$, where $(\bfw^*,b^*)$ is a P-stationary point of (\ref{SM-SVM}).   Hence, the number of $L_{0/1}$ support vectors are naturally expected to be no greater than the number of the standard  support vectors. This is also testified by our numerical experiments.

(C3)
When it comes to solving the problem (\ref{SM-SVM}), we adopt the famous alternating direction method of multipliers (ADMM), where one of its sub-problems is addressed by the  $L_{0/1}$ proximal operator involved in the P-stationary point, which together with the idea of $L_{0/1}$ support vectors allows us to define a working set in each step. Indices $i$ of vectors $\bfx_i$ out of this working set will be discarded, so the proposed method has a considerably low computational complexity and thus runs super fast. We prove that the limit point of the generated sequence is a  P-stationary point and also a locally optimal solution to the problem (\ref{SM-SVM}). This means the final classifier only uses a small number of support vectors based on the statements in C2.

(C4) Comparing with some leading classification solvers for addressing the SVM problems on synthetic and real datasets, extensive numerical experiments demonstrate that our proposed method achieves better performance including higher prediction accuracy,  a fewer number of support vectors and  faster computational speed. In addition, the numerical comparison also certifies the robustness to the outliers of  the $L_{0/1}$-SVM.

The remainder of this paper is organized as follows.  In the next section, a brief overview of various soft-margin loss functions used in \eqref{SM-SVM1} will be given. \cref{sec:opt} establishes the optimality theory including the existence of a globally optimal solution to the problem (\ref{SM-SVM}) and the relationships between a P-stationary point  and  an optimal solution.  In \cref{sec:alg}, we will introduce the $L_{0/1}$ support vectors and cast a fast ADMM whose each step is integrated by a  working set strategy inspired by the $L_{0/1}$ support vectors. Numerical experiments and concluding remarks are given in the last two sections.

\section{Related work}\label{sec:sub-prox}
The discrete nature of $\ell_{0/1}(\cdot)$ in  $L_{0/1}$-SVM (\ref{SM-SVM}) limits its wide applications. Therefore, most previous work \cite{HN2009, YW2019} focus on the continuous surrogates of (\ref{SM-SVM}), namely, $\ell(\cdot)$ in \eqref{SM-SVM1} is a continuous approximation of $\ell_{0/1}(\cdot)$. We mention two typical classes of such surrogate soft-margin  loss functions \cite{YY2016}. The first one consists of the convex  soft-margin  loss functions. An impressive body of work has  designed such kinds of functions since they make the corresponding SVM problems easier to deal with. Here, we only review some popular ones.
\begin{itemize}
\item \emph{Hinge soft-margin loss function:} $\ell_{\rm hinge}(t)=\max\{0, 1-t\}.$  It is  non-differentiable at $t=1$ and unbounded. SVM with hinge soft-margin loss function was first proposed by Vapnik and Cortes  \cite{CV95}, aiming at only penalizing the samples with $t<1$. Hinge soft-margin loss SVM is the first SVM model and { is widely studied by researchers \cite{BA2002}.}

\item  \emph{Pinball soft-margin loss function:} $\ell^{\tau}_{\rm pinball}(t)=\max\{1-t,-\tau (1-t)\}$, with $0\leq \tau \leq 1,$ which  is  still non-differentiable at $t=1$ and unbounded. SVM with this soft-margin loss function was proposed in \cite{JHS2013}, \cite{HSS2014} to pay penalty for  all training samples. There is a quadratic programming solver embedded in Matlab to solve the SVM with  pinball soft-margin loss function \cite{HSS2014}.
\item\emph{Huberized hinge soft-margin loss function:} $\ell^{\tau}_{\rm HH}(t)=\max\{1-t-\tau/2,\min\{\max\{1-t,0\}^2/2\tau,\tau/2\}\}$ with $\tau>0$. It is smooth but still unbounded function. SVM with such soft-margin loss function was first proposed in \cite{WZZ2008} which can be solved by proximal gradient method \cite{XAC2016}.

\item \emph{Square soft-margin loss function \cite{SV1999,YTH2014}:} $\ell_{\rm square}(t)=(1-t)^2,$ a smooth but unbounded function.

\item \emph{Other convex and smooth soft-margin loss functions} include the squared hinge soft-margin loss function \cite{TF2008} and log soft-margin loss function \cite{FHT2000}.

    \item \emph{Other convex and nonsmooth soft-margin loss functions} include the $\varepsilon$-insensitive zone pinball soft-margin loss function \cite{HSS2014} and $\phi$-risk hinge soft-margin loss function \cite{PL2008}.
\end{itemize}

As the above loss functions are convex, their corresponding SVM models are not difficult to be dealt with \cite{BA2002,JHS2013,HSS2014,WZZ2008,XAC2016,SV1999,YTH2014,TF2008,FHT2000,PL2008,PM2004,PM2006,JH1997}.
However, the convexity often induces the unboundedness \cite{MBB2000,PNF2003}, which weakens the robustness of those loss functions to outliers from the training data. To overcome such a drawback, one can set an upper bound and enforce the loss to stop increasing after a certain point. This gives rise to the second group: the nonconvex  soft-margin loss functions.
\begin{itemize}
\item\emph{Ramp soft-margin loss function \cite{YY2006, XL2014}:} $\ell_{\rm ramp}^{\mu}(t)=\max\{0,1-t\}-\max\{0,1-(t+\mu)\}$ with $ \mu > 0,$   which is non-differentiable at $t=1-\mu$ and $t=1$ but  bounded between 0 and $\mu$. It does not penalize the case when $t>1$, while pays linear penalty when $1-\mu\leq t\leq 1$ and a fixed penalty $\mu$ when $t<1- \mu$. This makes such a function robust to outliers.

\item\emph{Truncated pinball soft-margin loss function \cite{SNZT2017}} (truncated right side of pinball loss function): $\ell^{\tau,\kappa}_{\rm Tpin}(t)=\max\{0,(1+\tau)(1-t)\}-(\max\{0,\tau(1-t+\kappa)\}-\tau \kappa),$ with $  0\leq \tau \leq 1$ and $ \kappa \geq 0.$ It is non-differentiable at $t=1$ and $t= 1+\kappa$ and unbounded. The penalty is fixed at $\kappa$ for $t>1+\kappa$ and is linear otherwise. 

\item\emph{Asymmetrical truncated pinball soft-margin loss function \cite{YD2018}} (truncated two side of pinball loss function):  $\ell_{\rm ATpin}^{\tau,\kappa,\mu}(t)=\max\{0,(1+\tau)(1-t)\}-(\max\{0,\tau(1-t+\kappa)\}+\max\{0,1-t-\mu\}-\tau \kappa)$ with $0\leq \tau \leq 1$ and $ \mu,\kappa \geq 0.$  This function is non-differentiable at $t=1-\mu, t=1+\kappa$ and $t=1$ but bounded. The penalty is fixed at $\tau\kappa$ for $t>1+\kappa$ and at $\mu$ for $t<1-\mu$ but is linear otherwise. 

\item\emph{Sigmoid soft-margin loss function \cite{PNP2000}}:  $\ell_{\rm sigmoid}(t)$ $=1/(1+\text{exp}(-\tau(1-t))$ with $\tau>0$. It is a smooth and bounded function. It penalizes all training samples.

\item \emph{Other nonconvex and smooth soft-margin loss functions} include  the smooth ramp soft-margin loss function \cite{LH2008}, {savage loss \cite{HN2009}} and One-sided cauchy soft-margin loss function \cite{YY2016}, \cite{IA2008}.

\item \emph{Other nonconvex and nonsmooth soft-margin loss functions} include the truncated logistic soft-margin loss function \cite{SY2011}, { curriculum loss \cite{YW2019} }and $\varepsilon$-insensitive truncated least square soft-margin loss function \cite{DY2016}.
\end{itemize}

Compared to convex soft-margin loss functions, most nonconvex ones are less sensitive to feature noise or outliers due to their boundedness. Apparently, nonconvexity would lead to difficulties of computations in terms of solving the corresponding SVM models \cite{MBB2000,PNF2003,XL2014,YY2006,SNZT2017,YD2018,PNP2000,LH2008,IA2008,SY2011,DY2016,IN2004,SL2010,HN2009,YW2019}.

\section{ Optimality Theory of $L_{0/1}$-SVM}\label{sec:opt}
For convenience of our subsequent analysis, denote
\begin{eqnarray}
A&:=&[y_1\bfx_{1}~y_2\bfx_{2}~\cdots~y_m\bfx_{m}]^\top\in {\R}^{m\times n},\nonumber\\
 \bfy&:=&(y_{1},y_2,\cdots,y_{m})^{\top}\in {\R}^{m},\nonumber\\
\label{notation-Ay1u} {\bfe} &:=&(1,1,\cdots,1)^{\top}\in {\R}^{m},\\
 {\bfu}&:=&{\bfe}-A  \bfw -b\bfy\in {\R}^{m}, \nonumber\\
 {\bfu}_+&:=&((u_1)_+, \cdots,(u_m)_+)^\top \in \R^m,\nonumber
\end{eqnarray}
where $t_+:=\text{max}\{t, 0\}$. Moreover, the zero-norm of the vector  $\bfu$ is denoted by $\|\bfu\|_0$ which counts the number of its non-zero elements. It is easy to see that $u_i =1-y_{i}\langle \bfw , \bfx_i\rangle-y_{i}b=1-y_{i}f(\bfx_i)=1-t_i, i\in\N_m.$ Then the soft-margin loss function $\ell_{0/1}(\cdot)$ in (\ref{l01-loss})  can be rewritten as
\be \nonumber
~~~~~~~~~~~~~~~~~~~~~~~~~~~~~\ell_{0/1}(u_i)=
\begin{cases}
1,& u_i>0,\\
0,& u_i\leq0,
\end{cases}~~~i\in\N_m.~~~~~~~~~~~~~~~(4')
\ee
 This indicates
\begin{eqnarray}
\sum_{i=1}^m \ell_{0/1}(1-y_{i}f(\bfx_i))&=&\sum_{i=1}^m \ell_{0/1}(u_i)\nonumber\\
\label{L-0-1} &=&\|{\bfu}_+\|_0=:L_{0/1}({\bfu}).\end{eqnarray}
Hence,  the function $L_{0/1}({\bfu})=\|{\bfu}_+\|_0$ computes the number of all positive elements in ${\bfu}$. We call it the $L_{0/1}$ soft-margin loss function.  Borrowing these notation, the $L_{0/1}$-SVM (\ref{SM-SVM}) is equivalent to the following optimization problem,
\begin{eqnarray} \label{ReM-SVM}
\min_{ \bfw \in {\R}^n, b\in {\R}} f( \bfw ;b):=\frac{1}{2}\Vert \bfw  \Vert^2+C\|({\bfe}-A \bfw -b{\bf{y}})_{+}\|_{0},
\end{eqnarray}
or the following problem with an extra variable $\bfu$,
\begin{eqnarray} \label{reformulate}
&\underset{ \bfw \in {\R}^n,b \in \R,{\bfu} \in {\R}^m}{\min} & \frac{1}{2} \Vert \bfw  \Vert^2+C\|{\bfu}_{+}\|_{0}\\ \nonumber
&\mbox{s.t.} &{\bfu}+A  \bfw +b\bfy={\bfe}.
\end{eqnarray}

Recall the sparse optimization problem $\min_{\bfv\in \R^m}\{g(\bfv)+C\|\bfv\|_0\}$, where $C>0$ is a given penalty parameter and $g:\R^m\rightarrow \R$ is smooth or nonsmooth function. Due to the combinatorial nature of $\|\bfv\|_0$, the above sparse optimization problem is generally NP-hard. However, this problem has wide applications in linear and nonlinear compressive sensing, robust linear regression, deep learning, etc. Hence it has been extensively studied by a lot of researchers in different communities. More recently, by utilizing continuous optimization theory, the optimality conditions and algorithms for such a problem are successfully established by some researchers in optimization community \cite{TM2008,TM2009,ZY2013,ZS2014,AN2018,HL2020}.

Observe the $L_{0/1}$-SVM model (\ref{ReM-SVM}) or (\ref{reformulate}). We found that it has same structure as the above sparse optimization model with difference between $\|(\cdot)_+\|_0$ and $\|(\cdot)\|_0$. Similarly, by utilizing continuous optimization theory, we do the optimality analysis of (\ref{ReM-SVM}) or (\ref{reformulate}) in this section.

\subsection{Existence of $L_{0/1}$-SVM  Minimizer }\label{sub:exist-L01}
Firstly, we show the existence of a global minimizer (a minimizer is often phrased as an optimal solution) to (\ref{ReM-SVM}), a premise of the optimality condition of the $L_{0/1}$-SVM.
\begin{theorem}\label{exist-1} Given  $b\in\I:=[-M,M]$ with $0<M<+\infty$. Then the  globally optimal solution to (\ref{ReM-SVM}) exists and the solution set is bounded. \end{theorem}
The proof of \cref{exist-1} is given in Supplement S.1.
For any $b\in \I$, since  $y_i\in\{-1,1\}$, we have the following observations
\begin{eqnarray*}
  f(\text{\bf0};b)= C\|({\bfe}-b\bfy)_{+}\|_0
                  = \left\{\begin{array}{lccr}
Cm_{-},& b&\geq& 1,\\
Cm_{+} ,&b&\leq& -1,\\
Cm,& |b|&<&1,
 \end{array} \right.
\end{eqnarray*}
where $m_{+}$ and $m_{-}$ are the number of positive and negative $y_i$. Therefore, let $(\bfw^*; b^*)$ be an optimal solution to (\ref{ReM-SVM}) (such a solution exists by \cref{exist-1}), then
\begin{eqnarray*}
  f(\bfw^*;b^*)\leq ~{ C\min\{m_+,m_-\}.}
\end{eqnarray*}
In numerical experiments, this gives us a clue to set some starting points $(\bfw^0; b^0)$ satisfying
\begin{eqnarray}\label{start-point}
  f(\bfw^0;b^0)\leq ~{ C\min\{m_+,m_-\}.}
\end{eqnarray}
\subsection{ First-Order Optimality Condition}\label{sec:FOOC1}
From the perspective of optimization, establishing the optimality conditions of an optimization problem is a key step in theoretical analysis, because those conditions effectively benefits for the algorithmic design. Now turn our attention on the $L_{0/1}$-SVM model \eqref{reformulate}.
\begin{definition}
[P-stationary point of (\ref{reformulate})]\label{pstationarypoint} For a given $C>0$, we say   $({ \bfw^* }; b^*;  \bfu^* )$ is a proximal stationary (P-stationary) point of (\ref{reformulate}) if  there is a Lagrangian multiplier $ {\bm{\lambda}^*}\in \R^m$ and a constant $\gamma>0$ such that
\be\label{aaalfi}
\left\{
\begin{array}{rll}   \bfw^* + A^{\top}{\bm{\lambda}^*} &=& {\bf 0},
\\ \langle \bfy, {\bm{\lambda}^*} \rangle&=& {\bf 0},
\\  \bfu^*+A   \bfw^*+b^*\bfy&=&
{\bfe},
\\ \text{prox}_{\gamma C\|({\cdot})_{+}\|_{0}}(\bfu^*-\gamma{\bm{\lambda}^*})&=&\bfu^*,
\end{array}\right.
\ee
\end{definition}
where
 \be \label{exp-proximal1LL}
[\text{Prox}_{\gamma C\|({\cdot})_{+}\|_{0}}({ \bfz^*})]_i=
\begin{cases}
0,& 0< z^*_i\leq\sqrt{2\gamma C},\\
z^*_i,& z^*_i>\sqrt{2\gamma C}~ \text{or}~z^*_i\leq0 ,
\end{cases}
\ee
and $\bfz^*:=\bfu^*-\gamma \bm{\lambda}^*$.
The above equation (\ref{exp-proximal1LL}) is termed as $L_{0/1}$ proximal operator, whose solution has been derived in Supplement S.2.

The $L_{0/1}$ proximal operator is the key in the optimality analysis (see \cref{gol-p} below) and algorithmic design (see \cref{Fast}) of $L_{0/1}$-SVM. Using the above definition, we reveal the relationship between local/global minimizer and a P-stationary  point of $L_{0/1}$-SVM.
 To proceed more, let
\begin{eqnarray}\label{BH}
B:=[A~\bfy]\in\R^{m\times(n+1)},~~~~
H:=\begin{bmatrix}
I_{n\times n} & {\bf0}  \\
{\bf0} & 0
\end{bmatrix}B^+,
\end{eqnarray}
where  $B^+\in\R^{(n+1)\times m}$ is the generalized inverse of $B$, and
$\lambda_H:= \lambda_{\max}(H^\top H)$ where $\lambda_{\max}(H^\top H)$ is the maximum eigenvalue of $H^\top H.$  Thus, we have following theorem.
\begin{theorem}
\label{gol-p}  The following relations hold for (\ref{reformulate}).
\begin{itemize}
\item[(i)] A globally optimal solution  is also a P-stationary point with $0<\gamma< 1/\lambda_H$ if $B$ is full  column rank.
\item[(ii)] A P-stationary point with $\gamma>0$ is also a locally optimal solution.
\end{itemize}
\end{theorem}
The proof of \cref{gol-p} is given in Supplement S.3.
Note that $B$ being full column rank implies  $m>n$, i.e., the number of samples is greater  than the number of features.  However, from  \cref{gol-p} (ii), if we find a P-stationary point of the problem (\ref{reformulate}),  then it must be a locally optimal solution without any assumptions. No requirement of $m> n$ is enforced. Our numerical experiments testify that our proposed algorithm based on the idea of the P-stationary point works well for both cases: $m> n$ and $m\leq n$.

\subsection{Extension}
In \cref{sec:FOOC1}, we established the first-order optimality condition for (\ref{reformulate}), i.e., (\ref{ReM-SVM}), which is an unconstrained optimization problem.  This can be regarded as a special case of the following general optimization model
\begin{eqnarray} \label{inverse1}
\min_{{\bfu}\in {\R}^{m}} ~~ g(\bfu)+C\|{\bfu}_{+}\|_{0},
\end{eqnarray}
where  $C>0$ is a given penalty parameter and $g:{\R}^m\rightarrow {\R}$ is a smooth function and gradient Lipschitz continuous with a Lipschitz constant $\gamma_g>0$.

Similarly, we introduce the  proximal stationary  point of (\ref{inverse1}) as below.
\begin{definition}
[P-stationary point of (\ref{inverse1})] For a  given $C>0$, we say $ \bfu^* $ is a proximal stationary (P-stationary) point of problem (\ref{inverse1}) if  there is a  constant $\gamma>0$ such that
\begin{eqnarray}\label{equation}
 {\bfu}^*&=\text{prox}_{\gamma C \|(\cdot)_+\|_0}({\bfu}^*-\gamma \nabla g(\bfu^*)),
\end{eqnarray}
where, $\nabla g(\cdot)$ is the gradient of $g(\cdot)$.
\end{definition}
The following theorem reveals the relationship between a local/global minimizer and a P-stationary  point of (\ref{inverse1}), whose the proof is similar to that of  the  \cref{gol-p} and thus is omitted.
\begin{theorem}\label{gol-p-station}  For problem (\ref{inverse1}), the following relations hold.
\begin{itemize}
 \item[(i)]  For a  given $C>0$, if $ {\bfu^*}$ is a global minimizer of (\ref{inverse1}) then it is a P-stationary point with $0<\gamma< 1/\gamma_g$.\label{gol-p-sta}

 \item[(ii)] For a  given $C>0$, if $g$ is convex and ${\bfu^*}$ is a P-stationary point with $\gamma>0$, then it is a local minimizer of (\ref{inverse1}).\label{gol-p-sta2}
     \end{itemize}
\end{theorem}

 The above two theorems state that under condition of convexity, the P-stationary point must be a local minimizer, which means that we could use the  P-stationary  point as a termination rule in terms of guaranteeing the local optimality of a point generated by the algorithm proposed in next section.

\section{Fast Algorithm}\label{sec:alg}
It is well known that the classifier is decided by support vectors, see \eqref{supp-vects}. If support vectors is used to design the solving algorithm, the fewer number of support vectors is, the faster the computational speed will be since fewer samples in training data are used to train the classifier. Therefore, reducing the number of support vectors tends to be important for datasets in extremely large sizes. Motivated by this, we introduce $L_{0/1}$ support vectors and working set strategy based on the theory in \cref{sec:FOOC1} and adopt the famous alternating direction method of multipliers (ADMM) to solve the
$L_{0/1}$-SVM \eqref{reformulate}.

\subsection{$L_{0/1}$ Support Vectors}\label{sec:l01-sv}
Let ($\bfw^{*};b^{*};\bfu^{*}$) be a P-stationary  point of problem \eqref{reformulate}. Then from \cref{pstationarypoint}, there is a Lagrangian multiplier  $\bm{\lambda}^{*} \in \R^{m}$ and a constant $\gamma> 0$  such that \eqref{aaalfi} holds. Let
\begin{eqnarray}\label{index}
T_{*} :=\left\{i\in \N_m:~\bfu_i^{*}-\gamma \bm{\lambda}_i^{*} \in (0,\sqrt{2\gamma C}]\right\},
\end{eqnarray}
and $\overline{T}_{*}: =\N_m \backslash T_{*}$ be its complementarity set. Let $\bfz_T\in\R^{|T|}$ be the sub-vector of $\bfz$ indexed on $T$ and $|T|$ be the cardinality of $T$.  It follows from the last equation of \eqref{aaalfi} and \eqref{exp-proximal1LL} that
\begin{eqnarray*}
{\bfu^{*}} &\overset{\eqref{aaalfi} }{=}&\text{prox}_{\gamma C \|(\cdot)_+\|_0}(\bfu^{*}-\gamma \bm{\lambda^{*}})\\
&=&
 \left[\begin{array}{c}
(\text{prox}_{\gamma C \|(\cdot)_+\|_0}(\bfu^{*}-\gamma \bm{\lambda^{*}}))_{T_{*}}\\
(\text{prox}_{\gamma C \|(\cdot)_+\|_0}(\bfu^{*}-\gamma \bm{\lambda^{*}}))_{\overline{T}_{*}}
\end{array}\right]\\
&\overset{\eqref{exp-proximal1LL}}{=}&\left[\begin{array}{c}
{\bf0}_{T_{*}}\\
(\bfu^{*}-\gamma \bm{\lambda^{*}})_{\overline{T}_{*}}
\end{array}\right].
\end{eqnarray*}
which is equivalent to
\begin{eqnarray}\label{proxm1}
\begin{bmatrix}
\bfu^{*}_{T_{*}}\\
\bm{\lambda^{*}}_{\overline{T}_{*}}
\end{bmatrix}={\bf0}.
 \end{eqnarray}
  Then $T_{*}$ in \eqref{index} turns to
  \begin{eqnarray}\label{index-T*}
T_{*} =\left\{i\in \N_m:~  \lambda_i^{*} \in \Big[-\sqrt{2 C/\gamma}, 0\Big)\right\}.
\end{eqnarray}
This and \eqref{proxm1} result in
\begin{eqnarray}\label{lambdaSV}
\bm{\lambda}_i^{*}
\begin{cases}
\in[- \sqrt{2 C/\gamma},0),& \text{for} ~~i \in T_{*},\\
=0,& \text{for}~~ i \in \overline{T}_{*}.
\end{cases}
\end{eqnarray}
Taking \eqref{lambdaSV} into the first equation of \eqref{aaalfi} derives
\begin{eqnarray}\label{wsupportvecotor}\nonumber
\bfw^* &=& -A^{\top}_{T_{*}}{\bm{\lambda}^{*}_{T_{*}}} - A^{\top}_{\overline{T}_{*}}{\bm{\lambda}^{*}_{\overline{T}_{*}}}\\
&=&- A^{\top}_{T_{*}}{\bm{\lambda}^{*}_{T_{*}}}=\underset{i \in T_{*}}{\sum} - \lambda_i^{*}y_i\bfx_{i}.
\end{eqnarray}
\begin{remark}Regarding the expression \eqref{wsupportvecotor}, we have the following comments.
\begin{itemize}
\item  Recall \eqref{supp-vects}, where $ \bm{\alpha}^{*}$ is a solution to the dual problem of \eqref{HM-SVM}. From the optimization perspective, the Lagrangian multiplier $-\bm{\lambda}^*$ actually is  a solution to the dual problem of \eqref{reformulate}. In such a sense, $\left\{\bfx_{i}: i \in T_{*}\right\}$ indeed are standard support vectors. While we call them the $L_{0/1}$ support vectors since they are selected by the $L_{0/1}$ proximal operator.

\item Furthermore, the third equation in \eqref{aaalfi} implies ${\bfe}=\bfu^*_{T_*}+(A   \bfw^*+b^*\bfy)_{T_*}=(A \bfw^*+b^*\bfy)_{T_*}
$ due to $\bfu^*_{T_*}=0$ by \eqref{proxm1}, which and the definition \eqref{notation-Ay1u}  of $A$ yield
\begin{eqnarray}\label{proxm1222}
\langle\bfw^{*},\bfx_{i}\rangle+b^{*}=\pm1, ~\text{for}~ i \in T_{*}.
 \end{eqnarray}
Interestingly, the $L_{0/1}$ support vectors must fall into the support hyperplanes $\langle\bfw^{*},\bfx\rangle+b^{*}=\pm1$. As far as we know, the hard-margin SVM has such a property for linearly separable datasets. For linearly inseparable datasets, most soft-margin SVM can not guarantee this property. However, \eqref{proxm1222} is ensured by the $L_{0/1}$-SVM regardless of the datasets being separable or inseparable. This phenomenon manifests that the $L_{0/1}$-SVM could render fewer support vectors than the other soft-margin SVM models, which is also certified by our numerical experiments.
\end{itemize}
\end{remark}
The set $T_{*}$ in \eqref{index-T*} gives us a clue to select support vectors, which is very practical in the following algorithmic design.
\subsection {\ADMML~via Selection of Working Set\label{Fast}}
In this subsection, we take advantages of ADMM and working set to solve the $L_{0/1}$-SVM \eqref{reformulate}. We firstly give the framework of ADMM as follows. The augmented Lagrangian function of the problem (\ref{reformulate}) is given by
\begin{eqnarray}
L_{\sigma}(\bfw;b;\bfu;\bm{\lambda})=\frac{1}{2}\Vert \bfw \nonumber \Vert^2+C\|\bfu_{+}\|_{0}+\langle\bm{\lambda},\bfu-\bfe+A\bfw+b\bfy\rangle \\ \nonumber +\frac{\sigma}{2}\|\bfu-\bfe+A\bfw+b\bfy\|^{2},~~~~~~~~~~~~~~~~~~~~~~~~~~~~~~~
\end{eqnarray}
where $\bm{\lambda}$ is the Lagrangian multiplier and $\sigma>0$ is the penalty parameter. Given the $k$th iteration $(\bfw^k; b^k;\bfu^k;\bm{\lambda}^k)$, the framework to update each component  is as follows:
\begin{eqnarray}\label{ADMM}
\begin{array}{l}
{\bfu}^{k+1}~=\underset{{\bfu}\in {\R}^{m}}{\rm argmin}~L_{\sigma}( \bfw ^{k},b^{k},{\bfu},{\bm{\lambda}}^{k})\\
 \bfw ^{k+1}=\underset{{\bfw}\in {\R}^{n}}{\rm argmin}~L_{\sigma}( \bfw ,b^{k},{\bfu}^{k+1},{\bm{\lambda}}^{k})+\frac{\sigma}{2}\|\bfw-\bfw^k\|^2_{D_k}\\
b^{k+1}~~=\underset{{b}\in {\R} }{\rm argmin}~L_{\sigma}( \bfw ^{k+1},b,{\bfu}^{k+1},{\bm{\lambda}}^{k})\\
{\bm{\lambda}}^{k+1}~={\bm{\lambda}}^{k}+\eta\sigma(\bfu^{k+1}-\bfe+A\bfw^{k+1}+b^{k+1}\bfy),
\end{array}
\end{eqnarray}
where $\eta>0$ is the dual step-size. 
The proximal term is
$$\|\bfw-\bfw^k\|^2_{D_k}=\langle\bfw-\bfw^k,  D_k(\bfw-\bfw^k)\rangle.$$
Note that if  $D_k$ is positive semidefinite, then the above framework is the standard semi-proximal ADMM \cite{FPST13}. However, authors in {papers} \cite{LST16,CLZS19} have also investigated ADMM with the indefinite proximal terms, namely, $D_k$ is indefinite.  The basic principle  of choosing $D_k$  is to guarantee the convexity of ${\bfw}$-subproblem  of (\ref{ADMM}).  Since $L_{\sigma}( \bfw ;b^{k};{\bfu}^{k+1};{\bm{\lambda}}^{k})$ is strongly convex with respect to $\bfw$, $D_k$ is flexible to be chosen as an indefinite matrix.

Now, let's see how $T_{*}$ in \eqref{index-T*} instructs to select the support vectors. Denote $\bfz^k:=\bfe-A\bfw^{k}-b^{k}\bfy-{\bm{\lambda}^{k}}/{\sigma}$. Define a working set $T_k$ at the $k$th step by
\begin{eqnarray}
\label{prox_T}T_k :=\Big\{i\in \N_m:~z^k_i \in \Big(0,\sqrt{2C/\sigma}~\Big]\Big\}
\end{eqnarray}
and $\overline{T}_{k}: =\N_m \backslash T_{k}$. Based on which, $D_k$ is chosen as
\begin{eqnarray}\label{dk}
D_{k}= -A_{\overline{T}_k}^\top A_{\overline{T}_k}.
\end{eqnarray}
Here, for a given set $T\subseteq\N_m$, $A_{T}\in\R^{|T|\times n}$ denotes the sub-matrix containing rows of $A$ indexed on $T$.  The working set $T_k $ and the choice of $D_{k}$ will tremendously speed up the whole computation in each step of ADMM. More precisely, we calculate each sub-problem in \eqref{ADMM} as follows.

\noindent {\bf (i) Updating $\bfu^{k+1}$}: The $\bfu$-subproblem in (\ref{ADMM}) is equivalent to the following problem
\begin{eqnarray*}\label{prox_u0}\nonumber
&&\bfu^{k+1}\\\nonumber
&=&\underset{\bfu\in \R^{m}}{\rm argmin}~C\|\bfu_{+}\|_{0}+\langle\bm{\lambda}^{k},\bfu\rangle+\frac{\sigma}{2}\|\bfu-\bfe+A\bfw^{k}+b^{k}\bfy\|^{2}\\ \nonumber
&=&\underset{\bfu\in \R^{m}}{\rm argmin}~C\|\bfu_{+}\|_{0}+\frac{\sigma}{2}\|\bfu-\bfz^k\|^{2} \\
&=&\text{Prox}_{\frac{C}{\sigma}\|\bf(\cdot)_{+}\|_{0}}(\bfz^k),
\end{eqnarray*}
where the last equation is from  \eqref{exp-proximal1LL} with $\gamma=1/\sigma$. This together with \eqref{exp-proximal1LL} and the working set \eqref{prox_T} suffices to
\begin{eqnarray}
\label{prox_u}{\bfu}^{k+1}_{T_k}= {\bf0},~~~
{\bfu}^{k+1}_{\overline{T}_k}=  \bfz^k_{\overline{T}_k}.
\end{eqnarray}
Therefore, updating  $\bfu^{k+1}$ turns to be very simple and fast.

\noindent {\bf (ii) Updating $\bfw^{k+1}$.}
The ${\bfw}$-subproblem in (\ref{ADMM}) is
\begin{eqnarray}\label{wk1}
\bfw^{k+1}=\arg\min_{\bfw \in \R^{n}}\frac{1}{2}\|\bfw\|^2 + \frac{\sigma}{2}\| \bfw -  \bfw^k\|^{2}_{-A_{\overline{T}_k}^\top A_{\overline{T}_k}} \nonumber\\
+\langle\bm{\lambda}^{k},A\bfw\rangle+\frac{\sigma}{2}\|\bfu^{k+1}-\bfe+A\bfw+b^{k}\bfy\|^{2}.~~
\end{eqnarray}
It is a convex quadratic programming problem. To solve (\ref{wk1}), we only need to find a solution to the equations
\begin{eqnarray}\label{euqation-w-0}{\bf0}&=&\bfw-\sigma A_{\overline{T}_k}^\top A_{\overline{T}_k} (\bfw -  \bfw^k)+A^\top\bm{\lambda}^{k}\nonumber\\
&+&
 \sigma A^\top(\bfu^{k+1}-\bfe+A\bfw+b^{k}\bfy),
 \end{eqnarray}
 which is equivalent to find a solution to the  equations
\begin{eqnarray}\label{euqation-w}
 (I +\sigma A_{T_k}^\top A_{T_k}) \bfw = \sigma A_{T_k}^\top \bfv_{T_k}^k ,
\end{eqnarray}
where $\bfv^k:=-({\bfu}^{k+1}+b^{k}\bfy-{\bfe}+{\bm{\lambda}}^{k}/\sigma)$. To derive \eqref{euqation-w} from \eqref{euqation-w-0}, we used two facts that
${\bfu}^{k+1}_{\overline{T}_k}= \bfz^k_{\overline{T}_k}$ by
\eqref{prox_u}  and $$A_{T_k}^\top A_{T_k}=A^\top A-A_{\overline T_k}^\top A_{\overline T_k}.$$
Therefore, the term $ A_{\overline T_k}$ vanishes in  \eqref{euqation-w}, which means the working set $T_k $ and the choice of $D_{k}$ discard the samples $\{\bfx_j, j\in \overline T_k\}$. This would fasten  the computation significantly if the selected $|T_k|$ is very small.  In practice, (\ref{euqation-w}) can be addressed efficiently by the following rules:
\begin{itemize}
\item If $n \leq |T_k|$, one could solve (\ref{euqation-w}) directly through
 \begin{eqnarray}\label{w1}
 \bfw ^{k+1}= (I +\sigma A_{T_k}^\top A_{T_k})^{-1}  \sigma A_{T_k}^\top \bfv_{T_k}^k.
\end{eqnarray}
\item If $n > |T_k|$,  the  Sherman-Morrison-Woodbury formula \cite{FPST96} enables us to calculate the inverse as
 \begin{eqnarray}\label{inverse-P}\nonumber(I +\sigma A_{T_k}^\top A_{T_k})^{-1}=I-\sigma A_{T_k}^\top ({I +\sigma A_{T_k} A_{T_k}^\top})^{-1}A_{T_k}.\end{eqnarray}
 Then we  update $ \bfw^{k+1}$ by
 \begin{eqnarray}\label{w2}
  \bfw ^{k+1}=  \sigma A_{T_k}^\top ({I +\sigma A_{T_k} A_{T_k}^\top})^{-1} \bfv_{T_k}^k .
  \end{eqnarray}
\end{itemize}

\noindent {\bf (iii) Updating $b^{k+1}$.} The $b$-subproblem in (\ref{ADMM}) is a convex quadratic programming
\begin{eqnarray}
b^{k+1}=\arg\min_{b\in \R}\ \ \langle\bm{\lambda}^{k},b\bfy\rangle +\frac{\sigma}{2}\|\bfu^{k+1}-\bfe+A\bfw^{k+1}+b\bfy\|^{2}. \nonumber
\end{eqnarray}
which is solved by
\begin{eqnarray}\label{b-k}
b^{k+1} = \langle \bfy,\bfr^k\rangle/\|\bfy\|^2=  \langle \bfy,\bfr^k\rangle/m,
\end{eqnarray}
where $\bfr^k:=-A \bfw ^{k+1}+{\bfe}-{\bfu}^{k+1}-\bm{\lambda}^{k}/\sigma $.

\noindent {\bf (iv) Updating $\bm{\lambda}^{k+1}$.} We  update $\bm{\lambda}^{k+1}$ in (\ref{ADMM}) as follows
\begin{eqnarray}\label{lambda-k}
\bm{\lambda}^{k+1}_{T_k}=\bm{\lambda}^{k}_{T_k}+\eta\sigma\bm{\varpi}^{k+1}_{T_k}, ~~~\bm{\lambda}^{k+1}_{\overline{T}_k}= {\bf 0},
\end{eqnarray}
where $\bm{\varpi}^{k+1}:=\bfu^{k+1}-\bfe+A\bfw^{k+1}+b^{k+1}\bfy$ and  setting $\bm{\lambda}^{k+1}_{\overline{T}_k}= {\bf 0}$ follows the idea in \eqref{proxm1}, namely, the part of the Lagrangian multiplier not on the working set is removed.

Overall, updating each subproblem is summarized into \cref{Alg-SQREDM}, which is called \ADMML,  an abbreviation for $L_{0/1}$-SVM solved by  ADMM.
\begin{algorithm} 
	\caption{: \ADMML\ for solving problem (\ref{reformulate}) }
	\label{Alg-SQREDM}
\begin{algorithmic}
	\STATE{Initialize ($ \bfw ^0;b^0;\bfu^0;{\bm{\lambda}}^0$). Set $C,\eta,\sigma,K>0$ and  $k=0$.}
 \WHILE{The halting condition does not hold and $k\leq K$}
 	\STATE{Update $T_k$ as in (\ref{prox_T}).}
 \STATE{Update $\bfu^{k+1}$ by (\ref{prox_u}).}
\STATE{ Update $\bfw^{k+1}$ by (\ref{w1}) if $n\leq |T_k|$ and by (\ref{w2}) otherwise.}
\STATE{Update $b^{k+1}$ by (\ref{b-k}).}
\STATE{Update ${\bm{\lambda}}^{k+1}$ by (\ref{lambda-k}).}
\STATE{Set $k=k+1$.}
 \ENDWHILE
 \RETURN the final solution ($ \bfw ^k,b^k$) to (\ref{reformulate}).
\end{algorithmic}
\end{algorithm}
\subsection{Convergence and Complexity Analysis}\label{calculation}
The following theorem shows that if the sequence  generated by \ADMML\ has a limit point, then it must be  a P-stationary  point and also a locally optimal solution to (\ref{reformulate}).
\begin{theorem}\label{convergence}Suppose $( \bfw ^*;b^*;{\bfu}^*;{\bm{\lambda}}^*)$ be the limit point of the sequence  $\{( \bfw ^k;b^k;{\bfu}^k;{\bm{\lambda}}^k)\}$ generated by \ADMML. Then $({ \bfw^* }; b^*;\bfu^*)$ is  a P-stationary  point with $\gamma=1/\sigma$ and also a locally optimal solution to the  problem (\ref{reformulate}).
\end{theorem}

The proof of  \cref{convergence} is given in Supplement S.4. Based on the authors' limited knowledge, the above convergence result is difficult to improve, because our \ADMML\ deals with $L_{0/1}$-SVM directly, whose objective function involves a discrete part $\|(\cdot)_+\|_0$. As a supplement, we mention some works on ADMM and its convergence analysis: for solving nonconvex nonsmooth optimization problems, see, e.g. \cite{GT2018,GT2015,MZ2016,RD2020}; for solving the nonconvex soft-margin loss SVMs, see, e.g. \cite{NF2014,GL2018}.

 With regard to the computational complexity in each iteration of the proposed algorithm \ADMML, we have the following observations:
\begin{itemize}
\item Updating $T_k$ by (\ref{prox_T}) needs the complexity $\mathcal{O}(m)$.
\item The main term involved in computing $\bfu^{k+1}$ by (\ref{prox_u}) is $A\bfw^{k}$, taking the complexity about $\mathcal{O}(mn)$.
\item To update $\bfw^{k+1}$,  we compute (\ref{w1}) if $n\leq |T_k|$ and  (\ref{w2}) otherwise. For the former,  the dominant computations are calculating $$A_{T_k}^\top A_{T_k}~~\text{and}~~(I +\sigma A_{T_k}^\top A_{T_k})^{-1}.$$ Their computational complexities are $\mathcal{O}(n^2|T_k| )$ and $\mathcal{O}(n^\kappa)$ with $\kappa\in(2,3)$, respectively. For the latter, the dominant computations are from $$A_{T_k} A_{T_k}^\top~~\text{and}~~({I +\sigma A_{T_k} A_{T_k}^\top})^{-1} $$ with the computational complexities $\mathcal{O}(n|T_k|^2 )$ and $\mathcal{O}(|T_k|^\kappa)$ with $\kappa\in(2,3)$, respectively. Therefore, the  complexity to update $\bfw^{k+1}$ in each step is $$\mathcal{O}(\min\{n^2,|T_k|^2\} \max\{n,|T_k|\}).$$
\item Similarly,  $A\bfw^{k+1}$ is the most expensive computation in (\ref{b-k}) to derive $b^{k+1}$. Again its complexity is $\mathcal{O}(mn)$.
\item Same as that of updating  $b^{k+1}$, achieving ${\bm{\lambda}}^{k+1}$ by (\ref{lambda-k}) takes $\mathcal{O}(mn)$ complexity.
\end{itemize}
Overall, the whole computational complexity in each step of  \ADMML\ in \cref{Alg-SQREDM} is
$$\mathcal{O}\left(mn+\min\{n^2,|T_k|^2\} \max\{n,|T_k|\}\right).$$
If the selected working sets have  low cardinalities $|T_k|$ or  $n$ is very small (i.e., $n \ll m$), \ADMML\ possesses a considerably low computational complexity.

With regard to non-asymptotic analysis for finding stationary points of nonsmooth nonconvex functions, see, e.g. \cite{JH2020}.
\section{ Numerical experiments}
In this section, we conduct numerical experiments to show the sparsity, robustness and effectiveness of the proposed \ADMML\ (available at \texttt{\url{https://github.com/Huajun-Wang/L01ADMM}})  by using MATLAB (2018b) on a laptop of 32GB of memory and Inter Core i7 2.7Ghz CPU,
against {nine} leading solvers on synthetic data and real data.

 Inspired by \cref{gol-p}, the P-stationary point is taken as a stopping criteria in the experiments. In the implementation, we terminate the proposed algorithm if the point ($ \bfw ^k;b^k;{\bfu}^k;{\bm{\lambda}}^k$) closely satisfies the conditions in (\ref{aaalfi}), i.e., $$\max\{\theta^k_1,\theta^k_2,\theta^k_3,\theta^k_4\}<\texttt{tol},$$  where \texttt{tol} is the tolerance level  and
\begin{eqnarray*}
\theta^k_1&:=&\frac{\| \bfw ^k+A^{\top}_{T_k}{\bm{\lambda}}^k_{T_k}\|}{1+\| \bfw ^k\|},~~~~\theta^k_2~:=~\frac{|\langle \bfy_{T_k},{\bm{\lambda}}^k_{T_k}\rangle|}{1+| {T_k}|},\\
\theta^k_3&:=&\frac{\|{\bfu}^k-{\bfe}+A \bfw ^k+b^k\bfy\|}{\sqrt{m}},\\
\theta^k_4&:=&\frac{\|{\bfu}^k-\text{prox}_{ C/\sigma \|(\cdot)_+\|_0}({\bfu}^k-  {\bm{\lambda}}^k/\sigma )\|}{1+\|{\bfu}^k\|}.
\end{eqnarray*}

{\bf (a) Parameters setting.} In our algorithm, the parameters $C$ and $\sigma$ control the number of support vectors (see (\ref{prox_T})), so tuning  good choices of these two parameters is crucial. Hence, the standard 10-fold cross validation is employed in training datasets to select them, where $C$ is picked from $\{2^{-7},2^{-6},\cdots,2^{7}\}$  and $\sigma$ is tuned from $\{a^{-7}, a^{-6},\cdots, a^{7}\}$ with $a=\sqrt{2}$. The parameters with the highest cross validation accuracy are picked out. In addition, we set $\eta=1.618$,  maximum iteration number $K=10^3$ and the tolerance level  \texttt{tol}$=10^{-3}$. For the starting points, set ${\bfu}^0={\bm{\lambda}}^0={\bf0}$. As  mentioned in \cref{sub:exist-L01}, we choose $ \bfw ^0={\bfe}/100$ and $b^0=0$ if it meets \eqref{start-point}, and $\bfw^0=0$ and $b^0=1$ (or $-1$) otherwise.

{\bf (b) Benchmark classifiers.} {{There is an impressive body of algorithms that have been developed to solve classification problems. However, to conduct fair comparisons, we only select nine solvers that were programmed by MATLAB.  Eight of them address the SVM problem and one deals with the $L_2$-regularized logistic regression problem.}} All their parameters are also optimized by 10-fold cross validation to maximize accuracy.

\begin{enumerate}

\item[$\hsvm$]  SVM with the hinge  soft-margin loss is implemented by LibSVM (\cite{CC2011}, \texttt{\url{https://www.csie.ntu.edu.tw/~cjlin/libsvm/}}), where the parameter $C$ is selected from
$\Omega:=\{2^{-7}, 2^{-6},\cdots, 2^{7}\}.$

\item[$\ssvm$  ] SVM with the square  soft-margin  loss \cite{SV1999} is implemented by LibLSSVM  (\cite{KJ2002}, \texttt{\url{https://www.esat.kuleuven.be/sista/lssvmlab/}}), where the parameter $C$ is selected from  $\Omega$.

\item[$\psvm$ ] SVM with the pinball soft-margin loss can be tackled by the traversal algorithm (\cite{XL2016}, \texttt{\url{https://www.esat.kuleuven.be/stadius/ADB/huang/softwarePINSVM.php}}), where $C$  is turned from a union of $\Omega$ and the one in \cite{XL2016} and $\tau$ is set as $\{-1,-0.99,\cdots,0.99\}$ from \cite{XL2016}.

\item[$\rsvm$] SVM with the ramp  soft-margin  loss can be addressed by \texttt{CCCP} (\cite{YY2006}, \texttt{\url{https://github.com/RampSVM/RSVM}}), where  the core subproblem of \texttt{CCCP} is solved by the MATLAB built-in function {\tt quadprog}, while $C$ and  $\mu$ are selected from $\Omega$ and  $\{0.1,0.2,\cdots,1\}$.

\item[$\sosvm$] SVM with  the one-sided Cauchy soft-margin  loss  is solved by  the  iteratively reweighted algorithm (\texttt{IRA}  \cite{YY2016}, \texttt{\url{https://www.esat.kuleuven.be/stadius/ADB/feng/softwareRSVC.php}}). The key subproblem of \texttt{IRA} is solved by the $\text{CVX} $, and both $C$  and  $\nu$ are  tuned from  $\Omega$. %

 \item[$\logi$]  $L_2$-regularized logistic regression is addressed by employing Newton algorithm (\cite{TP2003}, \texttt{\url{https://github.com/tminka/logreg/}}), where the parameter $C$ is selected from $\Omega$.

  {\item[$\pega$]  SVM with the hinge  soft-margin loss is solved by employing Pegasos algorithm (\cite{SYN2011}, \texttt{\url{https://github.com/bruincui/Pegasos}}), where the parameter $C$ is selected from $\Omega$. The mini-batch size is 1 and the maximum number of iterations is 2$m$.

   \item[$\svrg$] SVM with the squared hinge  soft-margin loss  is addressed by employing {\tt SVRG} algorithm (\cite{RT2013}, \texttt{\url{https://github.com/codes-kzhan/SVRG-1/blob/master/SVM/svm_SVRG.m}}) with  $C$  selected from $\Omega$.  The  mini-batch size is 1 and the number of ``passes" is $S=1$. The default epoch length is 2$m$.

   \item[$\katy$]  SVM with the squared hinge  soft-margin loss  is addressed by Katyusha algorithm (\cite{ZA2018}, \texttt{\url{https://github.com/codes-kzhan/SVRG-1/blob/master/SVM/svm_Katyusha.m}}) with all parameters selected the same as these for \svrg.}

\end{enumerate}
 {In addition, all other parameters of the above nine algorithms are set to their default values.}

{\bf (c) Evaluation criteria.} To evaluate classification performance, we report {five} evaluation criteria: the testing accuracy (\ACC), the number of support vectors (\NSV), { the size of working set per iteration (\SWS), the total number of iterations (\TNI)} and the \CPU\ time (\CPU). Let { $\{(\bfx_j^{\rm test}, y_j^{\rm test}): j=1,\cdots,m_t\}$} be the testing samples data. The testing accuracy is defined as follows
$$\ACC:=1-\frac{1}{2m_t}\sum_{j=1}^{m_t} \Big|{\rm sign}(\langle\bfw^*, \bfx_j^{\rm test} \rangle +b^*)-y_j^{\rm test}\Big|,$$
where ${\rm sign}(\overline{a})=1 $ if $\overline{a}>0$ and  ${\rm sign}(\overline{a})=-1 $ otherwise, and $(\bfw^*, b^*)$ is a solution obtained by one solver. The accuracy measures the ability of a solver to correctly predict the class labels of new input samples.  The higher \ACC\ (or the smaller  \NSV,~{ \SWS,~\TNI}~or~\CPU) is, the better performance of a solver delivers.

\subsection{Comparisons with Synthetic Data}
\begin{figure}[h]
\centering
\begin{subfigure}{.24\textwidth}
	\centering
	\includegraphics[height=2.7cm,width=4cm]{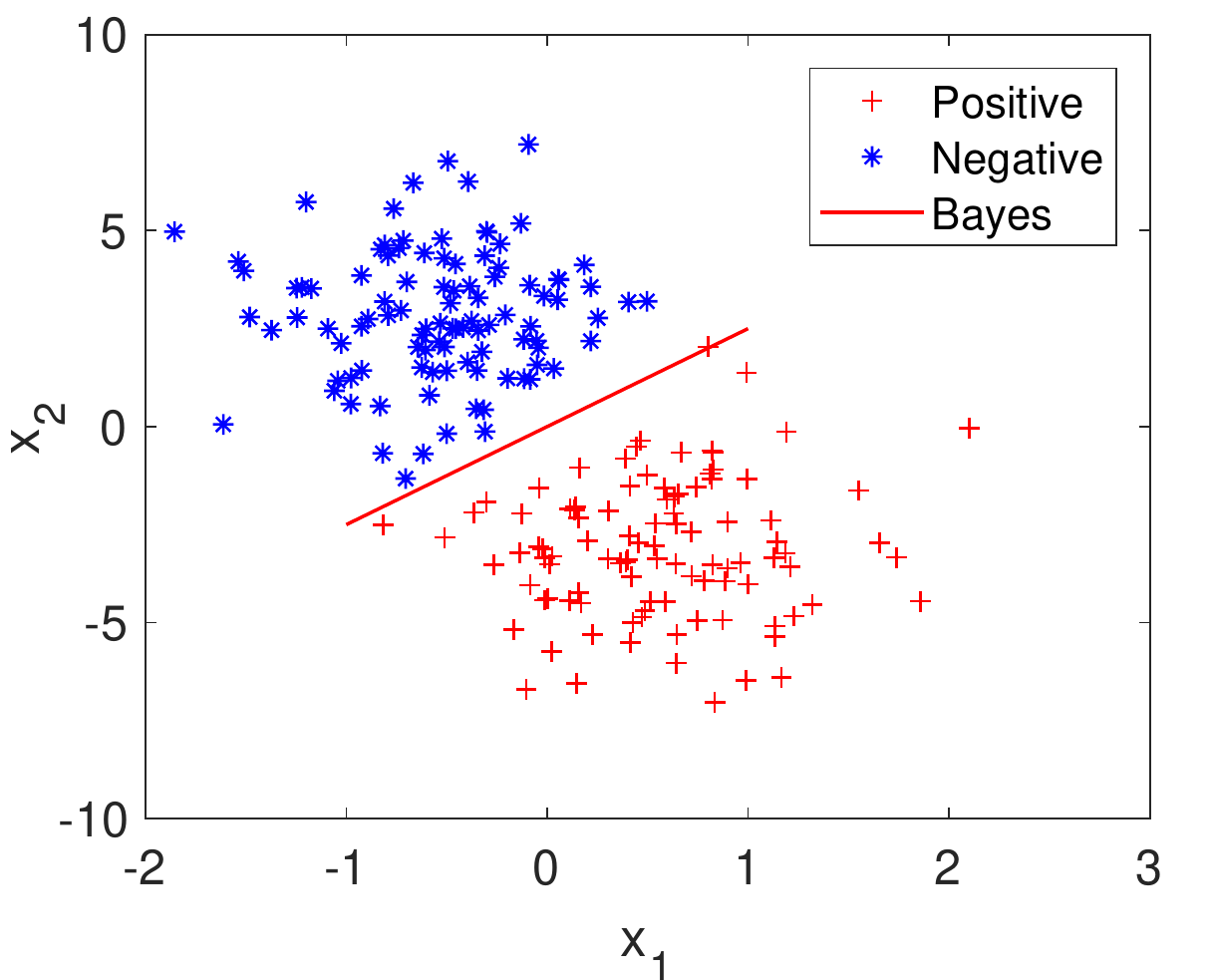}
\caption{}
\end{subfigure}%
\begin{subfigure}{.24\textwidth}
	\centering
	\includegraphics[height=2.7cm,width=4cm]{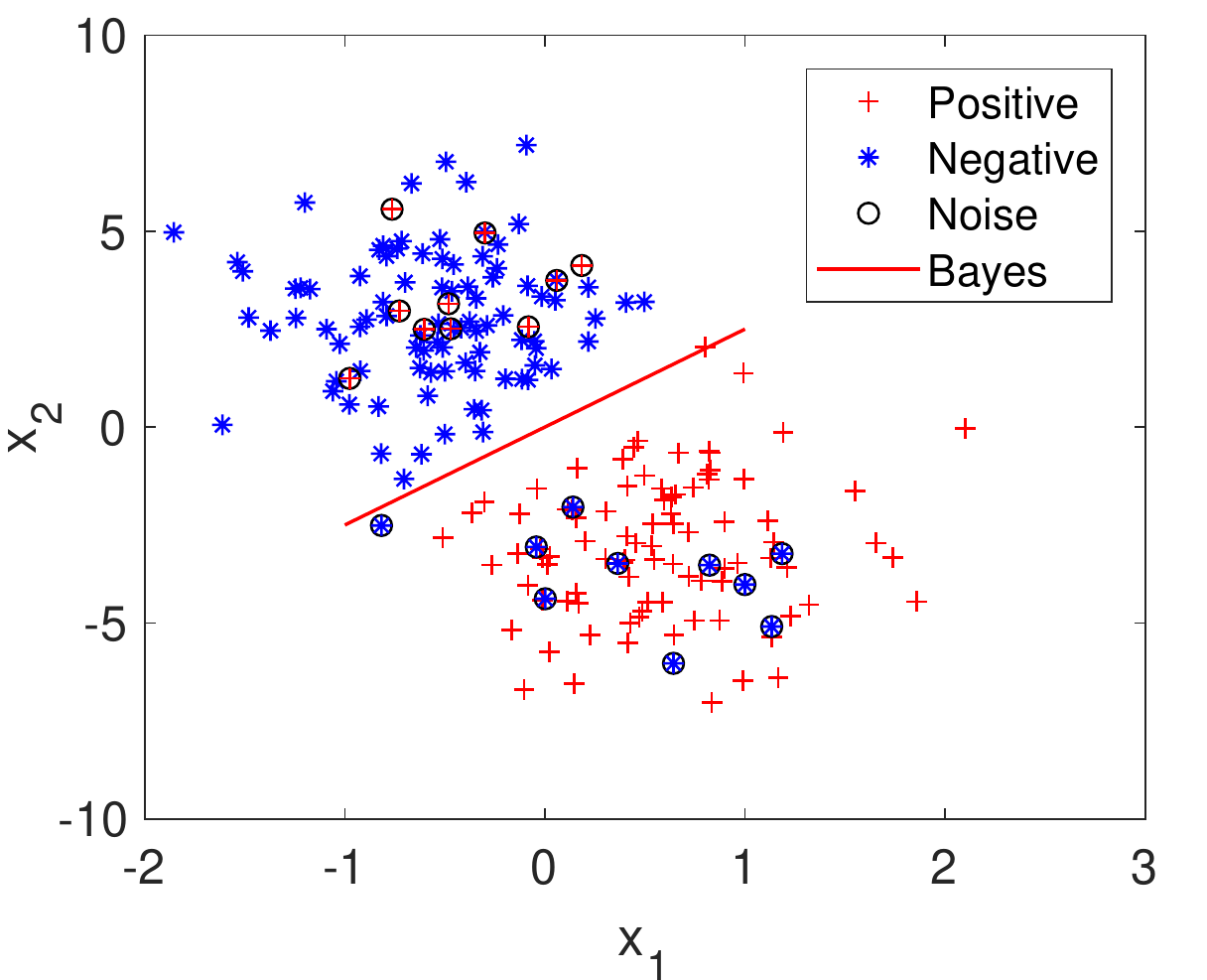}	
\caption{}
\end{subfigure}
 \caption{(a) A two dimensional training set with $200$ samples. (b)  Data in (a) but with $r$=10\% outliers. Blue stars: sampling samples in class $-1$. Red crosses: sampling samples in class $+1$. Red dashed lines: the Bayes classifier.}
\label{fig:rho-0}
\end{figure}

 For visualization, we first consider a two-dimensional example, where the features come from Gaussian distributions  \cite{HSS2014,XL2016}. One can observe that \ADMML\ performs extraordinarily in terms of delivering a considerably small number of support vectors.
\begin{example}[Synthetic data in $\R^2$ without outliers]\label{ex:syn-data-no-outliers} In this example, $m$ samples $\bfx_{i}, i\in\N_m$ with positive labels $y_i=+1$ are drawn from  $N(\bm{\mu}_{1},\Sigma_{1})$ and samples  $\bfx_{i}$ with negative labels $y_i=-1$ are drawn from $N(\bm{\mu}_{2},\Sigma_{2})$, where $\bm{\mu}_{1}=[0.5,-3]^{\top},\bm{\mu}_{2}=[-0.5,3]^{\top}$ and $\Sigma_{1}=\Sigma_{2}=\left[\begin{array}{cc}
0.2&0\\
0&3
\end{array}\right]$.  We generate $2m$ samples with two classes having equal numbers, and then evenly split all samples into a training set and a testing set.
\end{example}
Data generated in this way has centralized features of each class. For this example, the corresponding Bayes classifier is  $2.5x_1-x_2+0=0$. We display Bayes classifier and 200 training samples in \cref{fig:rho-0} (a), where samples are no extra noises  contaminated. We then add outliers on the data generated in \cref{ex:syn-data-no-outliers} as follows.

\begin{example}[Synthetic data in $\R^2$ with  outliers]\label{ex:syn-data-outliers} Firstly, $2m$  samples with  two  classes having equal numbers are generated as in \cref{ex:syn-data-no-outliers}. Then in each class, we randomly flip $r$ percentage of labels. For instance, in $m$ samples with positive labels $+1$, we change $mr$ labels to $-1$.  This means $r$ percentage of $2m$ samples are flipped their labels, namely $2rm$ outliers are generated. Here $r$ is the flapping ratio. Finally, the $2m$ samples are evenly split  into a training set and a testing set. In \cref{fig:rho-0} (b), the training set with $r$=10\% outliers are presented.
\end{example}
To solve these two  examples, {ten} solvers are applied to calculate the classifier $w_1x_1+w_2x_2+b=0$. Since data are generated randomly, to avoid randomness, we report average results  of \ACC, \NSV,~{\SWS,~\TNI}~and~\CPU~over 10 times.
%

{\bf (d) Synthetic data without outliers.}   {Ten} solvers are applied to solve \cref{ex:syn-data-no-outliers} with both the training and testing sample sizes being $m\in\{2000,4000,\cdots, 10000\}$. Average results are reported in \cref{Table-LU}, where "$--$" represents that the results are not obtained if one solver takes time longer than two hour (denote "$>2h$") or the required memory  is out of the capacity of our laptop (denote "$**$"), { and "3(8)" means the number of outer iterations (the average number of inner subproblem iterations).}  It can be clearly seen that all algorithms achieve desirable \ACC\ and \ADMML\ gets slightly better ones. When it comes to \NSV, the result is significant different. Obviously, \ssvm, \psvm\ and \logi\ take all samples as the support vectors, {
while \hsvm, \rsvm, \sosvm, \pega, \svrg~and~\katy~have a small number of  the support vectors.} It is evidently that \ADMML\ uses a considerably small number of the support vectors.
 { As for \SWS~and \TNI, \ssvm,~\psvm,~\rsvm,~\sosvm~and~\logi~take all samples as the working set, while most of them use a small \TNI~except for \psvm. By contrast, \ADMML~and others select a very small portion of samples as the working set, and \ADMML~uses a small \TNI~(no more than 50 for all cases).}
 Because of this,  \ADMML\ consumes the shortest \CPU\ time. 

{\bf (e) Synthetic data with outliers.} For \cref{ex:syn-data-outliers}, we fix $m=5000, n=2$ while alter the flapping ratio $r$ from  $\{0,0.05,0.1,0.15,0.2\}$ to see the robustness of each method to outliers. Average results are presented in \cref{Table-LU-2}. Apparently, the more outliers, the smaller \ACC\ for each solver. There is no big difference of \ACC\ generated by  {ten} solvers.  Again, \ADMML\ gets slightly better \ACC, being more robust to outliers than the others. Similar observations to that in \cref{Table-LU} can be seen for \NSV, { \SWS~and \TNI}. Moreover, the more outliers are added, {the more examples
become support vectors for \hsvm, \sosvm, \pega, \svrg\ and \katy}, { and bigger values of \TNI\ are generated by \hsvm\ and  \psvm}. By contrast, \ADMML\ makes use of fewer support vectors, { \SWS~and \TNI} when more outliers are added. Not surprisingly, \ADMML\ again runs  the fastest.  
\begin{table}[h]
\caption{Comparisons of 10 solvers  for solving \cref{ex:syn-data-no-outliers}, where $L_{0/1}$ stands for \ADMML.}\vspace{-3mm}
	\renewcommand{\arraystretch}{1.15}\addtolength{\tabcolsep}{-5pt}
	\label{Table-LU}
	\begin{center}
		\begin{tabular}{lp{0.65cm}crrrrrrrrr}
			\hline
&\multicolumn{6}{c}~~~~~~~~~~~~~~~~~~~~~{\ACC~(\%)} \\\cline{2-11}
			 $m$ &$L_{0/1}$~& {\hsvm} &{\ssvm} &{\psvm}&{\rsvm} &{\sosvm}&{\logi}&{\pega} &{\svrg}&{\katy}  \\\hline
				2000 	&	\textbf{97.05}&	\textbf{97.05}	&	97.00	&	\textbf{97.05}& \textbf{97.05}& \textbf{97.05} &97.03&
                   { 97.01}  &  {\textbf{97.05}} &	 {\textbf{97.05}}	\\
				4000  	&	\textbf{97.35}&	97.25	&	97.30	&	97.30& 97.33& 97.32&97.25&

{ 97.26} & {97.33}  &	{\textbf{97.35}}	\\
		        6000	&	\textbf{97.33}&	97.28	&	\textbf{97.33}	&	97.24	&		\textbf{97.33}& $--$&97.22&

{ 97.16} & {\textbf{97.33}}  &{97.30}	\\

				8000	&	\textbf{96.96}&	96.91	&	96.89&		96.91			&	\textbf{96.96}&$--$ &\textbf{96.96}&

{96.93}  &{96.94}  &{\textbf{96.96}}\\
				10000&	\textbf{97.20~}&	97.18	&	97.16	&	97.19			&	\textbf{97.20}&$--$ &97.18&

{97.16}  &{97.18}  &{97.18}\\\hline
&\multicolumn{6}{c}~~~~~~~~~~~~~~~~~~~~~{\NSV} \\\cline{2-11}
         		2000&	\textbf{7}	&	187	&	2000	&	2000	&	96&146 &	2000&

         {198} &{184}&{192}\\
          		4000&	\textbf{10}	&	301	&	4000	&	4000	&	141& 289	&4000&

          {325} &{332}&{295}		\\
                6000&	\textbf{18}	&	439	&	6000	&	6000	 &	201&$--$ &6000	&

                {453} &{444}&{452}	\\
         		8000&	\textbf{26}	&	571	&	8000	&	8000	 &	223&$--$ &8000	&

         {566} &{579}&{563}\\
         		10000&	\textbf{22}	&	658	&	10000	&	10000	 &	240&$--$ &10000	&{669} &{675}&{648}\\\hline

                  &\multicolumn{6}{c}~~~~~~~~~~~~~~~~~~~~~{ {\SWS}} \\ \cline{2-11}
	{2000}&	{22}&	{2}	&	{2000}	&	{2000}	&	{2000}&{2000} &	{2000}&
        { \bf{1}} &{ \bf{1}}&{ \bf{1}}\\
          		{4000}&	{31}	&	{2}	&	{4000}	&	{4000}	&	{4000}& {4000}	&{4000}&

          { \bf{1}} &{ \bf{1}}&{ \bf{1}}		\\
                {6000}&{35}	&	{2}	&	{6000}	&	{6000}	 &	{6000}&{6000} &{6000}	&

                { \bf{1}} &{ \bf{1}}&{ \bf{1}}	\\
         		{8000}&	{38}	&	{2}	&	{8000}	&	{8000}	 &	{8000}&{8000} &{8000}	&

         { \bf{1}} &{ \bf{1}}&{ \bf{1}}\\
         		{10000}&	{46}	&	{2}	&	{10000}	&	{10000}	 &	{10000}&{10000} &{10000}	&{ \bf{1}} &{ \bf{1}}&{ \bf{1}}\\\hline
                  &\multicolumn{6}{c}~~~~~~~~~~~~~~~~~~~~~{{\TNI}} \\\cline{2-11}
                {2000}&	{20}	&	{259}	&	{14}	&	{1216}	&	{3(8)}& {2(13)} & {\bf{8}}&{4000} & {4000} & {4000}		\\
         		{4000}&	{28}	&	{463}	&	{14}	&	{2325}	&	{3(16)}& {3(18)}& {\bf{9}}	& {8000} &{8000}  &{8000}	\\
                {6000}&	{34}	&	{639}	&	{15}	&	{3750}	&	{4(15)}& {$--$}&{\bf{9}}&{12000}  & {12000} &{12000}		\\
 		        {8000}&	{40}	&	{772}	&	{16}	&	{5247}	&	{4(21)}&{$--$} &{\bf{9}}& {16000} & {16000} &{16000}\\
         		{10000}&	{47}	&	{961}	&	{16}	&	{6326}	&	{5(23)}& {$--$}&{\bf{10}}& {20000} &{20000} &{20000}\\\hline
          &\multicolumn{6}{c}~~~~~~~~~~~~~~~~~~~~~{\CPU~(seconds)} \\\cline{2-11}
                2000&	\textbf{0.002}	&	0.014	&	0.221	&	9.642	&	3.969& 132.5 & 0.034&

                 {0.028} & {0.024}    & {0.025}\\

         		4000&	\textbf{0.006}	&	0.022	&	0.626	&	67.58	&	16.29& 2043  & 0.112&

         { 0.089}   & {0.087}   &	{0.088}\\
                6000&	\textbf{0.008}	&	0.036	&	1.200	&	209.9	&	31.44& $>2h$ & 0.204&

                { 0.133}   & {0.126}   &{0.131}\\
 		        8000&	\textbf{0.013}	&	0.069	&	2.342	&	493.2	&	65.25& $>2h$ & 0.536&

{ 0.194} & {0.185}  &{0.188}\\
         		10000&	\textbf{0.018}	&	0.094	&	3.951	&	775.3	&	124.7& $>2h$ & 0.938&

         { 0.281} & {0.266} & {0.268}	\\\hline
		\end{tabular}
	\end{center}
\end{table}\vspace{-0.3cm}
\begin{table}[h]
\caption{Comparisons of 10 solvers for solving \cref{ex:syn-data-outliers}}\vspace{-4mm}
	\renewcommand{\arraystretch}{1.15}\addtolength{\tabcolsep}{-4pt}
	\label{Table-LU-2}
	\begin{center}
		\begin{tabular}{lp{0.65cm}crrrrrrrr}
			\hline
&\multicolumn{5}{c}~~~~~~~~~~~~~~~~~~~~~{\ACC~(\%)} \\\cline{2-11}
			 $r$ &$L_{0/1}$&{\hsvm}  &{\ssvm}  &{\psvm}  &{\rsvm}&{\sosvm} &{\logi}&{{\pega}}&{{\svrg}}&{{\katy}}   \\\hline
				0.00 	&	\textbf{97.16}	&	97.08	&	97.10	&	\textbf{97.16}	&	\textbf{97.16}& 97.12 & 97.08&

{97.03}  & {\textbf{97.16}} & {\textbf{97.16}}\\

				0.05  	&	\textbf{92.65}	&	92.46	&	92.50	&	92.60	        &	\textbf{92.65}& 92.57 &	92.58&

{92.54}  & {92.30}  & {92.35}\\
	 	        0.10	&	\textbf{87.98}	&	87.78	&	87.78	&	87.90		    &	87.90& 87.90&87.70&

{87.68} & {87.46} & {87.45}\\
			    0.15	&	\textbf{83.06}	&	82.86	&	82.80	&	82.98			&	\textbf{83.06}& 83.04 &82.93&

{82.98}  & {82.88} &{82.88}\\
				0.20	&	\textbf{78.30}	&	78.16	&	78.12	&	78.28			&	78.28& 78.20&78.16&

{78.21} &{ 78.17} & {78.18}\\\hline
&\multicolumn{5}{c}~~~~~~~~~~~~~~~~~~~~~{\NSV} \\\cline{2-11}
				0.00 	&	\textbf{21}	&	364	&	5000	&	5000	&	184& 329	&5000&

{372} &{359}&{357}\\
				0.05  	&	\textbf{20}	&	947	&	5000	&	5000	&	175& 874	&5000&

{942} &{953}&{945}\\
	         	0.10	&	\textbf{17}	&	1385	&	5000	&	5000		&	170& 1015 &5000&

{1365} &{1373}&{1389}\\
				0.15	&	\textbf{16}	&	1795	&	5000	&	5000			&	161& 1657&5000&
{1790} &{1781}&{1792}\\
				0.20	&	\textbf{13}	&	2160	&	5000	&	5000			&	137& 1989 &5000&

{2177}&{2175}&{2187}\\	\hline
       &\multicolumn{5}{c}~~~~~~~~~~~~~~~~~~~~~{ {\SWS}} \\\cline{2-11}
				{0.00} 	&	{34}	&	{2}	&	{5000}	&	{5000}	&	{5000}& {5000}	&{5000}&

{ \bf{1}} &{ \bf{1}}&{ \bf{1}}\\
				{0.05}  	&	{31}	&	{2}	&	{5000}	&	{5000}	&	{5000}& {5000}	&{5000}&

{ \bf{1}} &{ \bf{1}}&{ \bf{1}}\\
	         	{0.10}	&	{30}	&	{2}	&	{5000}	&	{5000}		&	{5000}& {5000} &{5000}&

{ \bf{1}} &{ \bf{1}}&{ \bf{1}}\\
				{0.15}	&	{28}	&	{2}	&	{5000}	&	{5000}			&	{5000}& {5000}&{5000}&
{ \bf{1}} &{ \bf{1}}&{ \bf{1}}\\
				{0.20}	&	{27}	&	{2}	&	{5000}	&	{5000}			&	{5000}& {5000} &{5000}&

{ \bf{1}}&{ \bf{1}}&{ \bf{1}}\\	\hline
&\multicolumn{5}{c}~~~~~~~~~~~~~~~~~~~~~{ {\TNI}} \\\cline{2-11}
	    {0.00}&	{32}	&	{584}	&	{15}	&	{3042}	&	{3(25)}& {3(21)}& {\bf{9}}&

{10000}  & {10000} & {10000}	\\
				{0.05}&	{30}	&	{3726}	& 	{15} 	&	{3126}	&	{4(18)}& {3(21)}&	{\bf{9}}&

{10000} & {10000}  & {10000}\\
	            {0.10}&	{29}	&	{5128}	&	{15} 	&	{3268}		&	{4(17)}&{3(21)} &{\bf{9}}&

{10000} & {10000} & {10000}\\
	        	{0.15}&	{26}	&	{8423}	&	{15} 	&	{3373}			&	{5(13)}& {3(21)}&{\bf{9}}&

{10000} & {10000} &{10000}\\
				{0.20}&	{25}	&	{10776}	&	{15} 	&	{3443}			&	{5(13)}&{3(21)} &{\bf{9}}&
{10000} & {10000} & {10000}\\
       \hline
       &\multicolumn{5}{c}~~~~~~~~~~~~~~~~~~~~~{\CPU~(seconds)} \\\cline{2-11}
			    0.00&	\textbf{0.008}	&	0.027	&	0.801	&	93.11	&	22.53& 4047& 0.149&

{0.117}  & {0.108} & {0.112}	\\
				0.05&	\textbf{0.008}	&	0.075	&	0.823	&	101.3	&	20.99& 4069 &	0.131&

{0.119} & {0.114}  & {0.115}\\
	            0.10&	\textbf{0.006}	&	0.123	&	0.853	&	105.4		&	19.43&4084 &0.147&

{0.118} & {0.111} & {0.112}\\
	        	0.15&	\textbf{0.005}	&	0.172	&	0.885	&	108.3			&	18.96& 4092&0.152&

{0.118} & {0.110} &{0.111}\\
				0.20&	\textbf{0.005}	&	0.236	&	0.898	&	110.6			&	18.41&4094 &0.165&
{ 0.119} & {0.115} & {0.116}\\
       \hline
		\end{tabular}
	\end{center}
\end{table}
\subsection{Comparisons with Real Data}
We now apply these solvers  to deal with 14 real datasets. Their information are presented in \cref{Table-LU-3}, where the last six datasets have the testing data.
\begin{example}[Real data without outliers]\label{ex:real-data-no-outlier}
 We perform 10-fold cross validation for the first eight datasets. Each one is randomly split into ten parts, with one part being used for testing and the rest being used for training. We then record average results to evaluate performance. In our experiments,  all features are scaled to $[-1,1]$.\vspace{-0.3cm}

\end{example}
\begin{table}[h]
\caption{Descriptions of 14 real datasets. }\vspace{-1mm}
	\renewcommand{\arraystretch}{1}\addtolength{\tabcolsep}{-2pt}
	\label{Table-LU-3}
	\begin{center}
		\begin{tabular}{lrrrr}
			\hline
			 &  Training~data & Testing data&  Features   \\ \cline{2-4}
			Datasets & $m$& $m_t$ & $n$ \\\hline
            Colon-cancer (\texttt{col})& 62&0&2000\\
			Australian (\texttt{aus})    & 690         &0 &    14                 \\
			Two-norm (\texttt{two})    & 7400 & 0 &  20   \\
			Mushrooms (\texttt{mus})  & 8124 &  0 &           112    \\
			Adult (\texttt{adu})  &  17887    & 0 &13           \\
			Covtype.binaty (\texttt{cov})    & 581012         &0 &    54                 \\
			SUSY (\texttt{sus})  & 5000000&  0&           18    \\
			HIGGS (\texttt{hig})  &  11000000   & 0 &28           \\
			Lekemia (\texttt{lek})&38&34&7129\\
             Splice  (\texttt{spl})    &1000        &2175 &    60                 \\
           A6a  (\texttt{a6a})    & 11220        &21341 &    123                \\
			W6a (\texttt{w6a})    & 17188 & 32561 &  300    \\
			W8a (\texttt{w8a})  & 49749 &  14951 &          300    \\
			ijcnn1 (\texttt{ijc})  & 49990   & 91701 &22           \\
			\hline
		\end{tabular}
	\end{center}
\end{table}

\begin{table}[!htbp]
\caption{Comparisons of 10 solvers for solving \cref{ex:real-data-no-outlier}}\vspace{-4mm}
	\renewcommand{\arraystretch}{1.01}\addtolength{\tabcolsep}{-4.2pt}
	\label{Table-LU-4}
	\begin{center}
		\begin{tabular}{lrrrrrrrrrr}
			\hline
&\multicolumn{6}{c}~~~~~~~~~~~~~~~~~~~~~{\ACC~(\%)} \\\cline{2-11}
			 \text{Data} & $L_{0/1}$ &{\hsvm}  &{\ssvm}  &{\psvm}  &{\rsvm}&{\sosvm} &{\logi} &{{\pega}}&{{\svrg}}&{{\katy}}  \\\hline
    	\texttt{col}  	&	\textbf{90.23} 	&	64.52	        &	85.48	&	77.69	&	89.68  & 85.87  & 86.74 & {89.68}  &  {89.68}  & {89.68}		\\
        \texttt{aus} 	&	\textbf{86.23}	&	85.51	        &	85.80	&	85.80	&	86.02  & 85.98   &	86.18& {86.04} & {86.18} & {\textbf{86.23}}\\
	    \texttt{lek}  	&	\textbf{82.35}	&	58.82       	&	79.41	&	58.82	&	76.47  & \textbf{82.35}&\textbf{82.35} &{\textbf{82.35}} & {\textbf{82.35}} &{\textbf{82.35}}	\\
     	\texttt{spl}	&	85.52         &	\textbf{88.97} 	&	 85.75 	&	 85.52 	&	85.47  & 85.47    & 85.15 &

     {84.18}    & {85.44}  &{85.33}	\\
		\texttt{two}	&	\textbf{98.37}	&	98.02	        &	97.97	&	97.97	&	98.24  & $--$     & 97.78 &

{98.10}   & {\textbf{98.37}} &{98.24}\\
	    \texttt{mus}	&	\textbf{100.0}	&	\textbf{100.0}	&	\textbf{100.0}&	\textbf{100.0} &	\textbf{100.0}&$--$ &\textbf{100.0} &
 {\textbf{100.0}}   &   {\textbf{100.0}}  &  { \textbf{100.0}}\\
        \texttt{adu}	&	\textbf{83.90}	&	83.29	        &	83.01	&	83.07	&	83.79  & $--$     & 82.95 &

        { 83.29}       & {83.34}   &{\textbf{83.90}}\\		
     	\texttt{a6a}	&	\textbf{84.90} 	&	 84.18         	&	84.55 	&	84.69	&	 84.72 & $--$     & 84.76 &

     { 84.36}     &{84.72}  & {84.78}		\\
		\texttt{w6a}  	&	\textbf{97.93}	&	97.21	        &	97.58	&	97.21	&	97.86  & $--$     &	95.13     &

{97.24}     &{97.61}   & {97.57}	\\
     	\texttt{w8a} 	&	\textbf{98.54}	&	98.27	        & 	 $--$   &	$--$	&	$--$   & $--$     &	$--$ &

     {97.43}	& { 97.57} & {97.59}\\
		\texttt{ijc}  	&	\textbf{94.33}	&	92.73	        &	$--$	&	$--$	&	$--$   & $--$     & $--$ &

{93.49} & {93.35} & {93.56}\\
        \texttt{cov} 	&	\textbf{71.79}	&	$--$            &	$--$	&	$--$	&	$--$   & $--$     & $--$ &

        {68.93} & {69.83}&{69.77}\\
	    \texttt{sus} 	&	\textbf{67.58}	&	$--$	        &	$--$	&	$--$	&	$--$   & $--$     &	$--$ &

{64.28} & {65.62} & {65.86}\\
	    \texttt{hig}	&	\textbf{65.21}	&	$--$	        &	$--$    &	$--$	&	$--$   & $--$     & $--$&

{58.12}& {59.13}& {59.46}\\\hline

&\multicolumn{6}{c}~~~~~~~~~~~~~~~~~~~~~{\NSV } \\\cline{2-11}
		 \texttt{col}  	&	\textbf{34}  	&	46	            &	54	    &	54	    &	38     & 40     & 54    &

{46} &{45} &{46}\\
		 \texttt{aus}  	&	\textbf{24}	    &	203	            &	621 	&	621 	&	89	   & 177     &	621 &

{198} &{195} &{202}	\\
          \texttt{lek}   &	\textbf{26}  	&	31	            &	38    	&	38	    &	29	   & 31	      & 38 &

          {33}&{34} &{31}\\
    	\texttt{spl}	&   \textbf{70} 	&	 607 	        &	 1000 	&	 1000 	&	 87    & 332      & 1000  &
    {632} &{615}&{612}	\\
		\texttt{two}	&\textbf{30}    	&	758          	&	6600	&	6600	&	108    & $--$     &6600&

{783} &{775}&{788}\\
		\texttt{mus}	& \textbf{135}	&	550	&	7311	&	7311	&	506	&	$--$&   7311   &

{578}&  {575}  &{568}\\
        \texttt{adu}	&   \textbf{113}	&	6379	        &	16098	&	16098	&	1247   & $--$     & 16098&

        {6407}&  {6386}&{6394}\\
		\texttt{a6a}  	    &	\textbf{370}	&	4346	        &	11220	&	11220	&	1247   & $--$     & 11220&

{4562}&{4575}&{4582}	\\
       \texttt{w6a}  	 &	\textbf{429}	&	1128          	&	17188	&	17188	&	946    & $--$    &17188	&

       {1146}&{1152}&{1138}\\
       \texttt{w8a} 	 &	\textbf{867}	&	2857	        &	$--$	&	$--$	&	$--$   & $--$     & $--$&

       {2582}&{2579}&{2561}\\
	   \texttt{ijc} 	  &	\textbf{215}	&	8508	&	$--$	&	$--$	&	$--$		&	$--$	&$--$ &

{8535}&{8612}&{8608}\\
        \texttt{cov} 	  &	\textbf{137}	&	$--$         	&	$--$	&	$--$	&	$--$   & $--$     & $--$  &

        {$>$3e5}&{$>$3e5}&{$>$3e5}\\
		\texttt{sus} 	 &	\textbf{730}	&	$--$	        &	$--$	&	$--$	&	$--$   & $--$     & $--$&
{$>$2e6}&{$>$2e6}&{$>$2e6}\\
		\texttt{hig}	&  \textbf{1338}	&	$--$	        &	$--$	&	$--$	&	$--$   & $--$     & $--$ &
{$>$5e6}&{$>$5e6}&{$>$5e6}\\\hline

       &\multicolumn{6}{c}~~~~~~~~~~~~~~~~~~~~~{ {{\SWS}}} \\\cline{2-11}
		 \texttt{col}  	&	{37}  	&	{2}	            &	{54}	    &	{54}	    &	{54}     & {54}     & {54}    &

{ \bf{1}} &{ \bf{1}} &{ \bf{1}}\\
		 \texttt{aus}  	&	{66}	    &	{2}	            &	{621} 	&	{621} 	&	{621}	   & {621}     &	{621} &

{ \bf{1}} &{ \bf{1}} &{ \bf{1}}	\\
          \texttt{lek}   &	{29}  	&	{2}	            &	{38}    	&	{38}	    &	{38}	   & {38}	      & {38} &

          { \bf{1}}&{ \bf{1}} &{ \bf{1}}\\
    	\texttt{spl}	&   {94} 	&	 {2} 	        &	{1000}	&	{1000}	&	{1000}   &{1000}     &{1000} &
    { \bf{1}} &{ \bf{1}}&{ \bf{1}}	\\
		\texttt{two}	&{136}    	&	{2}         	&	{6600}	&	{6600}	&	{6600}  &{$--$}&{6600}&

{ \bf{1}} &{ \bf{1}}&{ \bf{1}}\\
		\texttt{mus}	&{ 772}	&	{2}	&	{7311}	&	{7311}	&	{7311}	&	{$--$}&   {7311}   &

{ \bf{1}}&  { \bf{1}}  &{ \bf{1}}\\
        \texttt{adu}	& { 1105}	&	{2}	        &	{16098}	&	{16098}	&	{16098}   & {$--$}     & {16098}&

        { \bf{1}}&  { \bf{1}}&{ \bf{1}}\\
		\texttt{a6a}  	    &{569}	&	{2}	        &	{11220}	&	{11220}	&	{11220}   & {$--$}     & {11220}&

{ \bf{1}}&{ \bf{1}}&{ \bf{1}}	\\
       \texttt{w6a}  	 &	{656}	&	{2}          	&	{17188}	&	{17188}	&	{17188}    & {$--$}    &{17188}	&

       { \bf{1}}&{ \bf{1}}&{ \bf{1}}\\
       \texttt{w8a} 	 &	{1284}	&	{2}	        &	{$--$}	&	{$--$}	&	{$--$}   & {$--$}     & {$--$}&

       { \bf{1}}&{ \bf{1}}&{ \bf{1}}\\
	   \texttt{ijc} 	  &{829}	&	{2}	&	{$--$}	&	{$--$}	&	{$--$}		&	{$--$}	&{$--$} &

{ \bf{1}}&{ \bf{1}}&{ \bf{1}}\\
        \texttt{cov} 	  &	{1520}	&	{$--$}         	&	{$--$}	&	{$--$}	&	{$--$}   & {$--$}     & {$--$}  &

        { \bf{1}}&{ \bf{1}}&{ \bf{1}}\\
		\texttt{sus} 	 &	{2814}	&	{$--$}	        &	{$--$}	&	{$--$}	&	{$--$}   & {$--$}     & {$--$}&
{ \bf{1}}&{ \bf{1}}&{ \bf{1}}\\
		\texttt{hig}	& { 3225}	&	{$--$}	        &	{$--$}	&	{$--$}	&	{$--$}   & {$--$}     & {$--$} &
{ \bf{1}}&{ \bf{1}}&{ \bf{1}}\\\hline
              &\multicolumn{6}{c}~~~~~~~~~~~~~~~~~~~~~{ {{\tt TNI}}} \\\cline{2-11}
			\texttt{col}  	&	{30} 	&	{41} 	&	{ \bf{2}}  	&	{31} 	& {2(2)}  & {2(4)} &{4}&

{108} & { 108}  &{ 108}\\
			\texttt{aus}  	&	{25}	&	{423} 	&	{17} 	&	{869} 	& {2(7)} &{3(26)} &	{\bf{6}}&	{1242} & {1242} & {1242}\\
            \texttt{lek}  	&	{18} 	&	{89}	&	{\bf{2}} 	&	{42}	&	 {2(2)}	 &{3(17)}&	{25}&

            { 76} &{76} &{76}\\
         	\texttt{spl}	&	{63} 	&	{ 595} 	&	{28}  	&	 {1276} 	&{2(9)} & {4(28)}	  & {\bf{9}} &

         {2000}   & {2000}  & {2000}\\
	    	\texttt{two}	&	{50}	&	{660}	&	 {75}	&	{3417}	&	{4(11)}  &	{$--$}&{\bf{12}} &
{13200}   & {13200}  & {13200}\\
			\texttt{mus}	&	{21}	&	{1623}	&	{106}  	&	{3685}	&	{4(12)}  &	{$--$}  & {\bf{18}} &

{14622}   & {14622} &{14622}\\
            \texttt{adu}	&	{26} 	&	{4766}	&	{157} 	&	{7720}		&	{5(21)}  &	{$--$} & {\bf{15}}&

            {32196}  & {32196}   &{32196}\\			
		    \texttt{a6a}  	&	{183}	&	{3032}	&	{289} 	&	{6873}			& {5(27)}	&{$--$} &	{\bf{16}}&

{22440}     &{22440}  & {22440}\\
            \texttt{w6a}  	&	{121}	&	{1450}	&	{404} 	&	{14417}		&	{7(32)} &	{$--$} &	{\bf{28}}  &

            {34376}     &{34376}   & {34376}	\\
      	    \texttt{w8a} 	&	{\bf{195}}	&	{8124}	&	{$--$}	&	{$--$}	&	{$--$}		&	{$--$} & {$--$} &

      {99498}	&  {99498} & {99498}	\\
			\texttt{ijc}  &{\bf{146}}	&	{6681}	&	{$--$}	&	{$--$}	&	{$--$}	&		{$--$} &{$--$} &

{99980} &{99980} &{99980}	\\
           	\texttt{cov} &	{\bf{103}}	&	{$--$}	&	{$--$}	&	{$--$}	&	{$--$}& 	{$--$} &{$--$} &

           {1.05$e^6$} & {1.05$e^6$}  &{1.05$e^6$}\\
			\texttt{sus} 		&	{\bf{117}}	&	{$--$}	&	{$--$}	&	{$--$}	&	{$--$}	&		{$--$} &{$--$}	&

{9.0$e^6$} & {9.0$e^6$} & {9.0$e^6$}	\\
			\texttt{hig}	&	{\bf{124}}	&	{$--$}	&	{$--$}	&	{$--$}		&	{$--$}& 	{$--$}&{$--$}& {1.98$e^7$} & {1.98$e^7$} & {1.98$e^7$}\\
       \hline
       &\multicolumn{6}{c}~~~~~~~~~~~~~~~~~~~~~{{\tt CPU} (seconds)} \\\cline{2-11}
			\texttt{col}  	&	0.021 	&	0.009 	&	\textbf{0.001} 	&	0.010 	&0.003 &1.488 & 0.182 &

{0.015} & { 0.012}  &{ 0.014}\\
			\texttt{aus}  	&	0.005	&	0.014 	&	0.033 	&	0.874 	&0.650  &87.23&	0.021&	{\textbf{0.004}} & {\textbf{0.004}} & {\textbf{0.004}}\\
            \texttt{lek}  	&	0.072 	&	0.057	&	\textbf{0.004}	&	0.010	&	0.008	 &54.36&	36.10&

            { 0.029} &{ 0.024} &{0.026}\\
         	\texttt{spl}	&	0.043 	&	 0.117 	&	 0.083 	&	 7.976 	&	0.631& 384.2 & 0.151 &

         {0.036}   & {\textbf{0.032}}  & {0.033}\\
	    	\texttt{two}	&	\textbf{0.054}	&	0.265	&	2.506	&	516.7	&	139.2 &	$>2h$&1.591 &
{0.171}   & {0.164}  & {0.166}\\
			\texttt{mus}	&	\textbf{0.074}	&	0.997	&	3.419	&	769.5	&	153.4 &	$>2h$  & 6.942 &

{0.422}   & {0.412} &{0.416}\\
            \texttt{adu}	&	\textbf{0.576}	&	3.775	&	24.58	&	1633.4		&	1013.2 &	$>2h$ & 5.032&

            { 0.775}  & {0.732}   &{ 0.744}\\			
		    \texttt{a6a}  	&	\textbf{0.172}	&	4.405	&	40.64	&	1472.5			&	1037.3&$>2h$ &	6.046&

{1.083}     &{1.025}  & {1.031}\\
            \texttt{w6a}  	&	\textbf{0.226}	&	1.532	&	170.9	&	5947.2			&	2747.4&	$>2h$ &	41.21  &

            {1.314}     &{1.186}   & {1.232}	\\
      	    \texttt{w8a} 	&	\textbf{2.576}	&	64.33	&	$**$	&	$**$	&	$>2h$		&	$>2h$ & $**$ &

      { 4.863}	&  {4.227} & {4.316}	\\
			\texttt{ijc}  &	\textbf{0.573}	&	36.95	&	$**$	&	$**$	&	$>2h$	&		$>2h$ &$**$ &

{1.526} &{1.247} &{ 1.316}	\\
           	\texttt{cov} &	\textbf{3.870}	&	$**$	&	$**$	&	$**$	&	$**$& 	$**$ &$**$ &

           {14.37} & {13.88}  &{13.91}\\
			\texttt{sus} 		&	\textbf{10.38}	&	$**$	&	$**$	&	$**$	&	$**$	&		$**$ &$**$	&

{137.6} & {132.4} & {133.7}	\\
			\texttt{hig}	&	\textbf{14.26}	&	$**$	&	$**$	&	$**$		&	$**$& 	$**$&$**$& {281.3} & {269.5} & {270.1}\\
       \hline
		\end{tabular}
	\end{center}
\end{table}
 \begin{example}[Real data with outliers]\label{ex:real-data-outlier} To see the influence of the real data with outliers, we select six datasets from small  sizes to moderate  sizes in \cref{Table-LU-3}. They are \texttt{col}, \texttt{aus},  \texttt{two}, \texttt{mus},  \texttt{adu} and \texttt{a6a}.   Same processes as in \cref{ex:real-data-no-outlier} are then applied into the first five datasets.  Finally, $r$ percentage of training and testing samples are randomly treated as outliers (i.e., their labels are flipped).
\end{example}

{\bf (f) Real data without outliers.} The average results are recorded in \cref{Table-LU-4}, {where ``$>3e5$" represents the number greater than $300000$.}
It can be clearly seen that \ADMML\ outperforms the others in terms of the highest \ACC, smallest \NSV\ and shortest \CPU\ for   most datasets, { and uses a small~\SWS~and \TNI.} For instance, \ADMML\ predicts more than 90\% samples correctly for \texttt{col} whilst {\hsvm} and {\psvm} only get less than 80\% correct predictions.  Compared with those generated by the other {nine} solvers, \NSV\ from \ADMML\ is relatively small. { As for \SWS, \ADMML\ takes  a small samples as the working set, which testifies that our constructed working set strategy is very effective to reduce the cost of per iteration. As for \TNI, \ADMML\ uses  a few \TNI~compared with \psvm,~\pega,~\svrg~and~\katy.} { For the computational speed, {\pega}, {\svrg} and \katy\ present the advantage of
\CPU\ for dealing with small scale datasets.} The \ADMML\ runs super fast for datasets in  big sizes,   0.573 seconds v.s. 36.95 seconds by {\hsvm} for data \texttt{ijc}. In addition, it only needs 14.26 seconds  for the dataset \texttt{hig} with more than ten million samples.  Overall, it seems that the bigger $m$ is, the more evident the advantage of \ADMML\ becomes.

\begin{figure}[h]
\centering
\begin{subfigure}{0.24\textwidth}
	\centering
	\includegraphics[height=2.7cm,width=4cm]{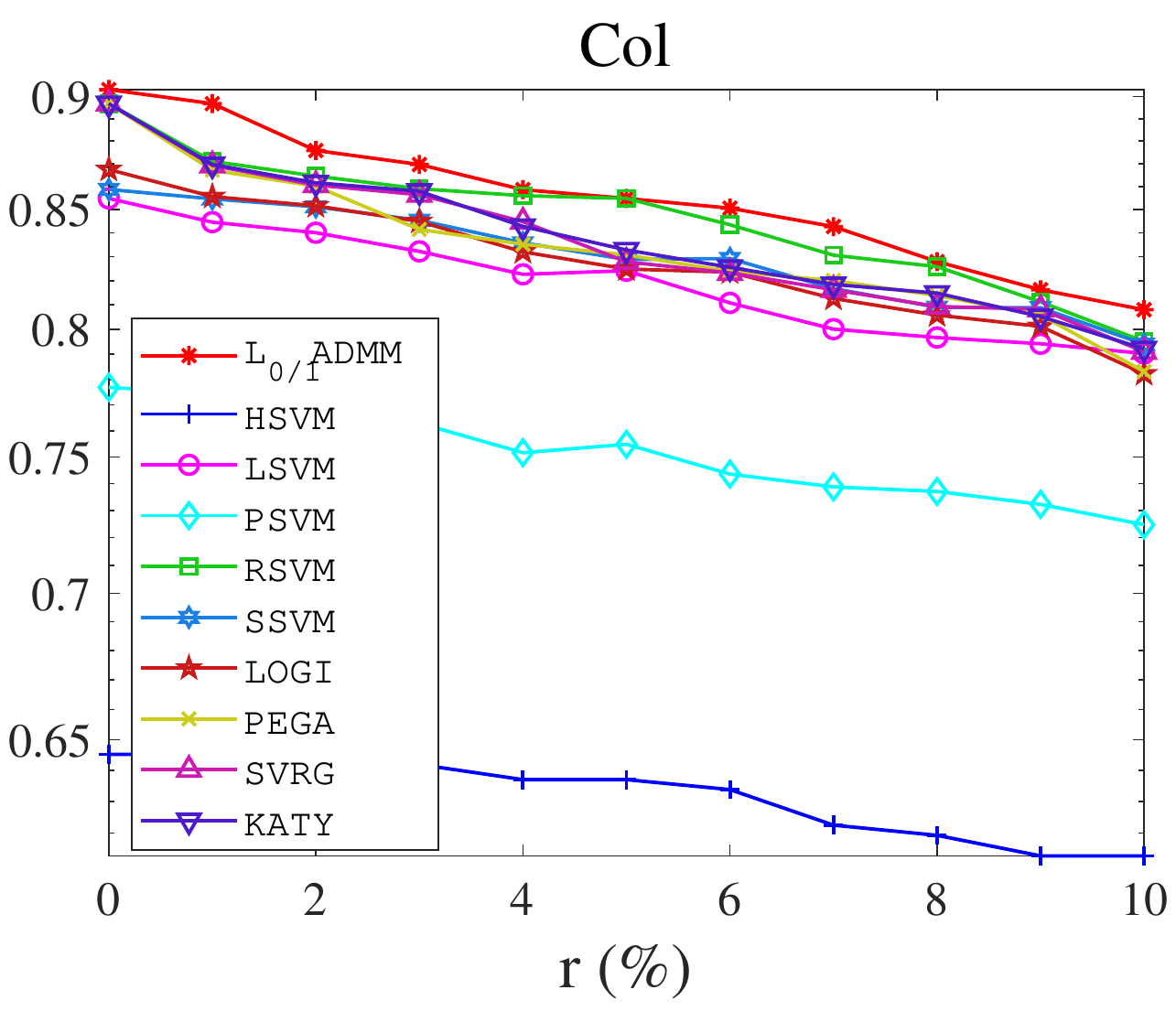}
\end{subfigure}%
\begin{subfigure}{0.24\textwidth}
	\centering
	\includegraphics[height=2.7cm,width=4cm]{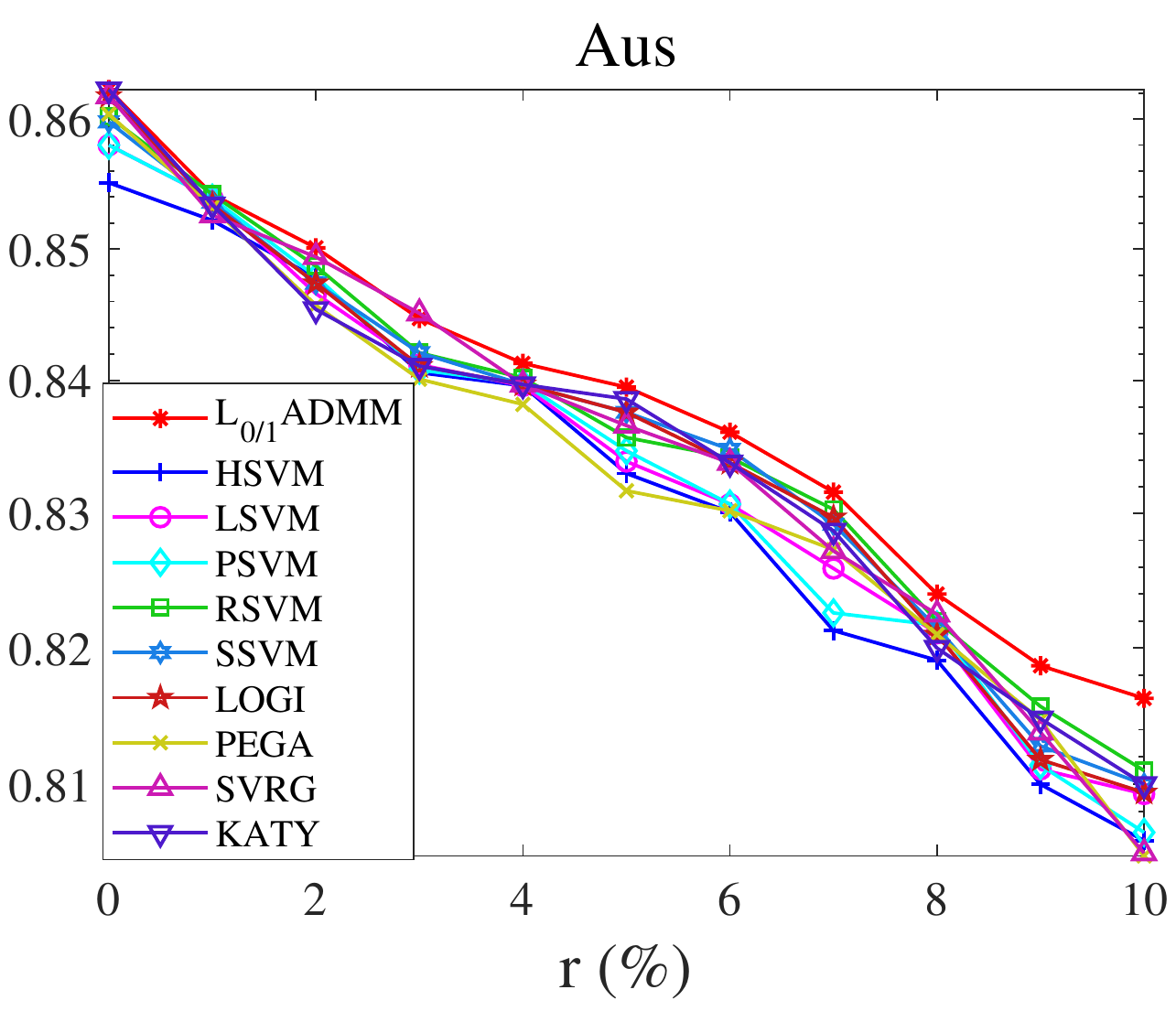}	
\end{subfigure}
\begin{subfigure}{0.24\textwidth}
	\centering
	\includegraphics[height=2.7cm,width=4cm]{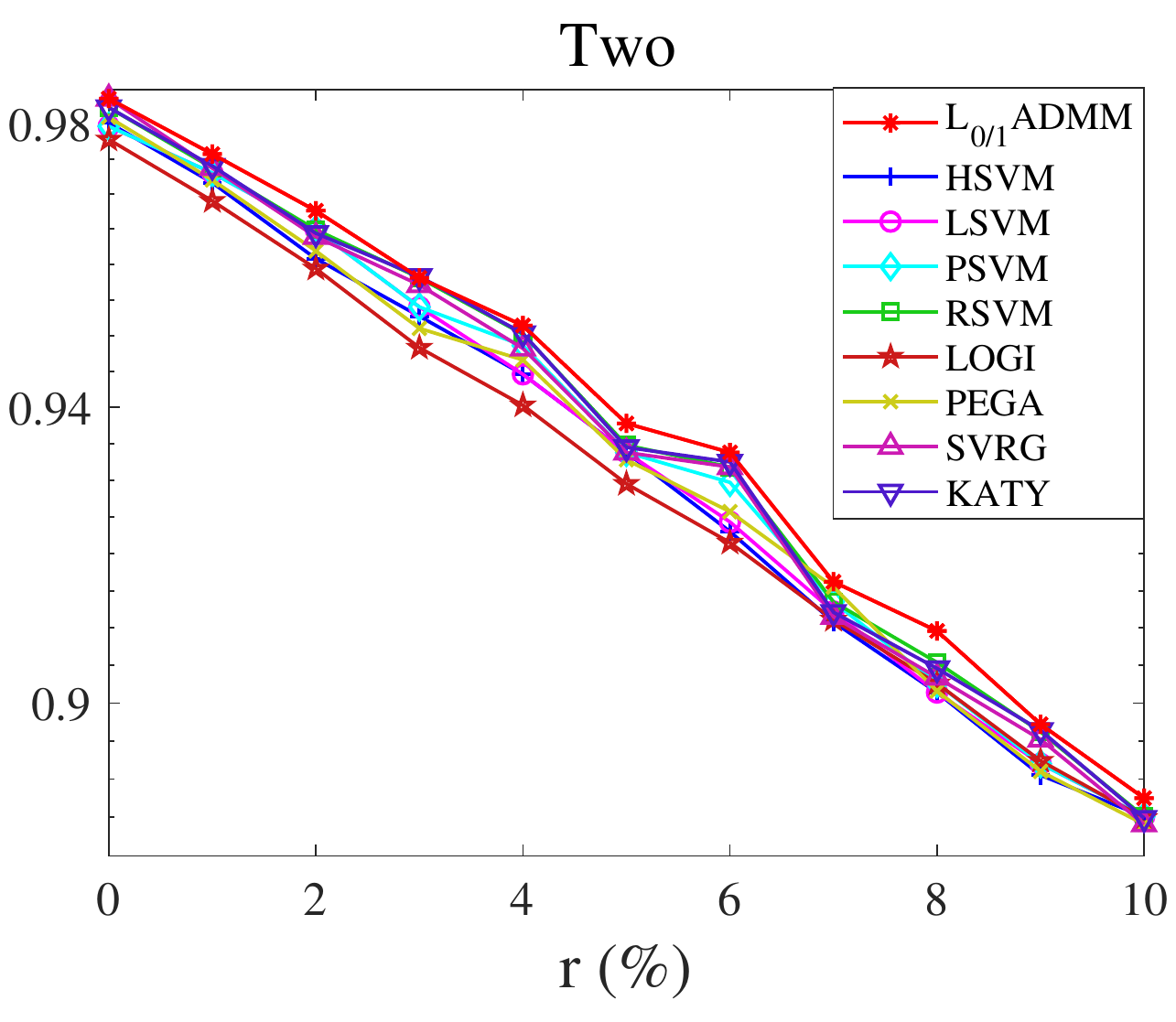}
\end{subfigure}%
\begin{subfigure}{0.24\textwidth}
	\centering
	\includegraphics[height=2.7cm,width=4cm]{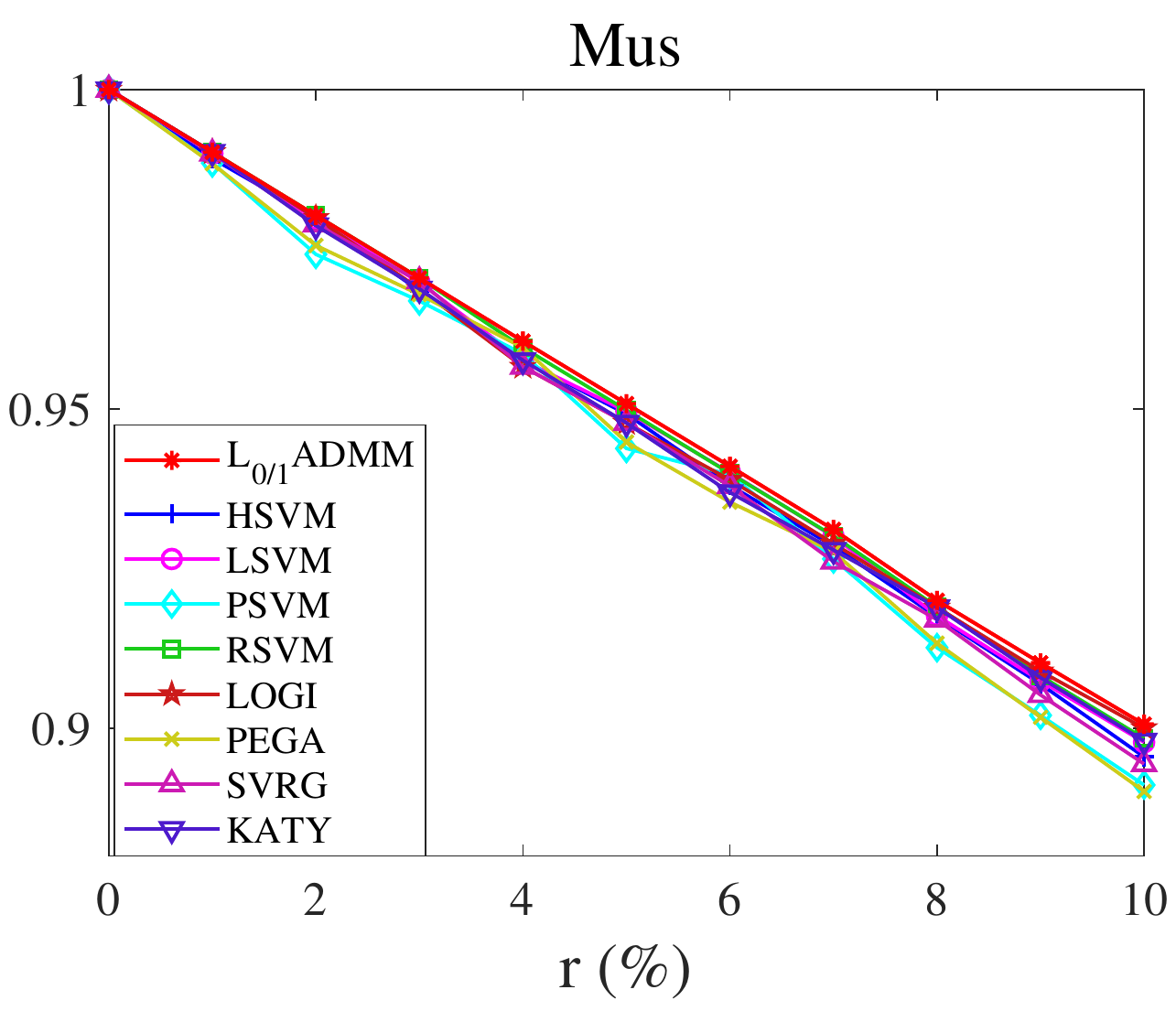}	
\end{subfigure}
\begin{subfigure}{0.24\textwidth}
	\centering
	\includegraphics[height=2.7cm,width=4cm]{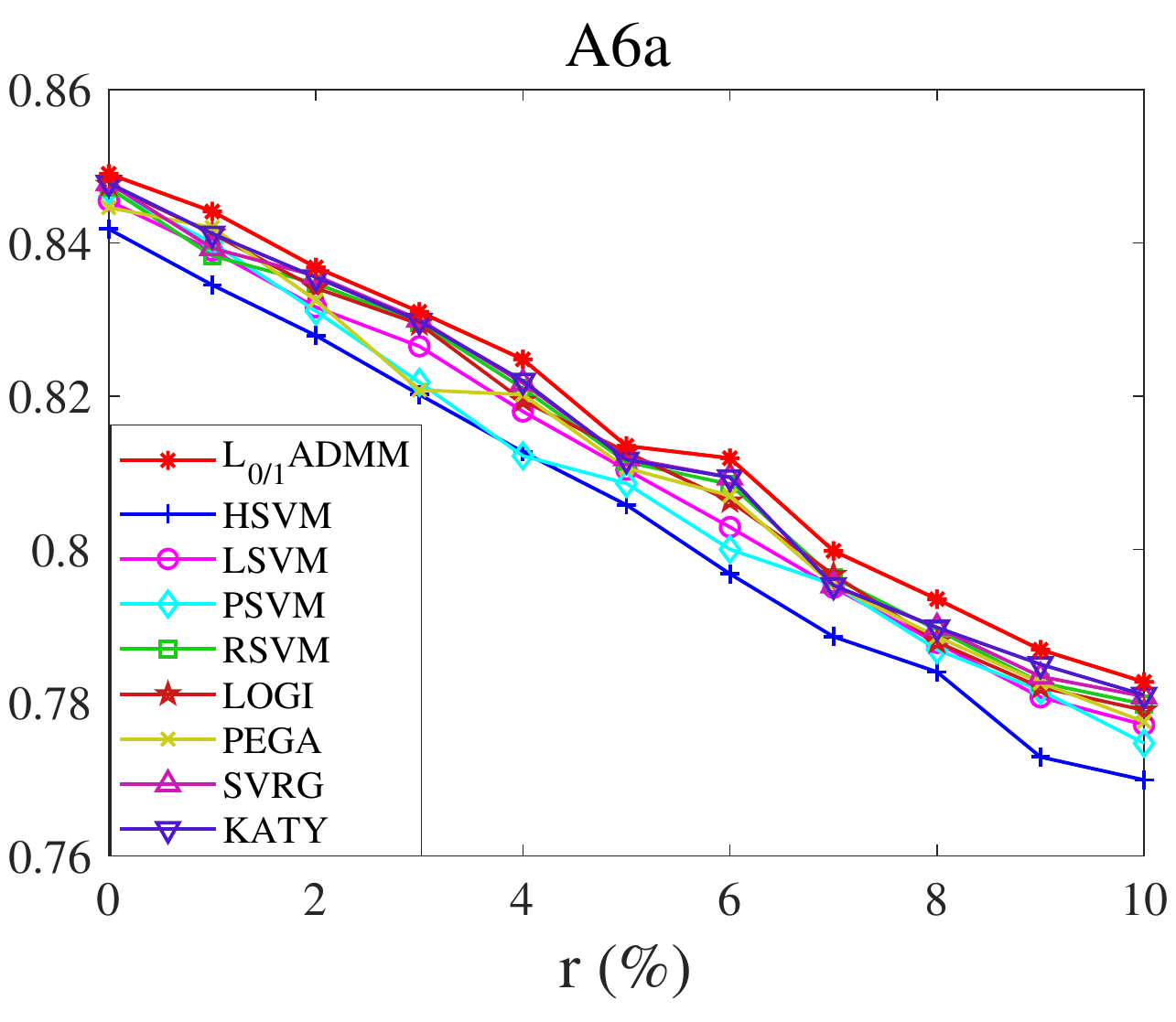}
\end{subfigure}%
\begin{subfigure}{0.24\textwidth}
	\centering
	\includegraphics[height=2.7cm,width=4cm]{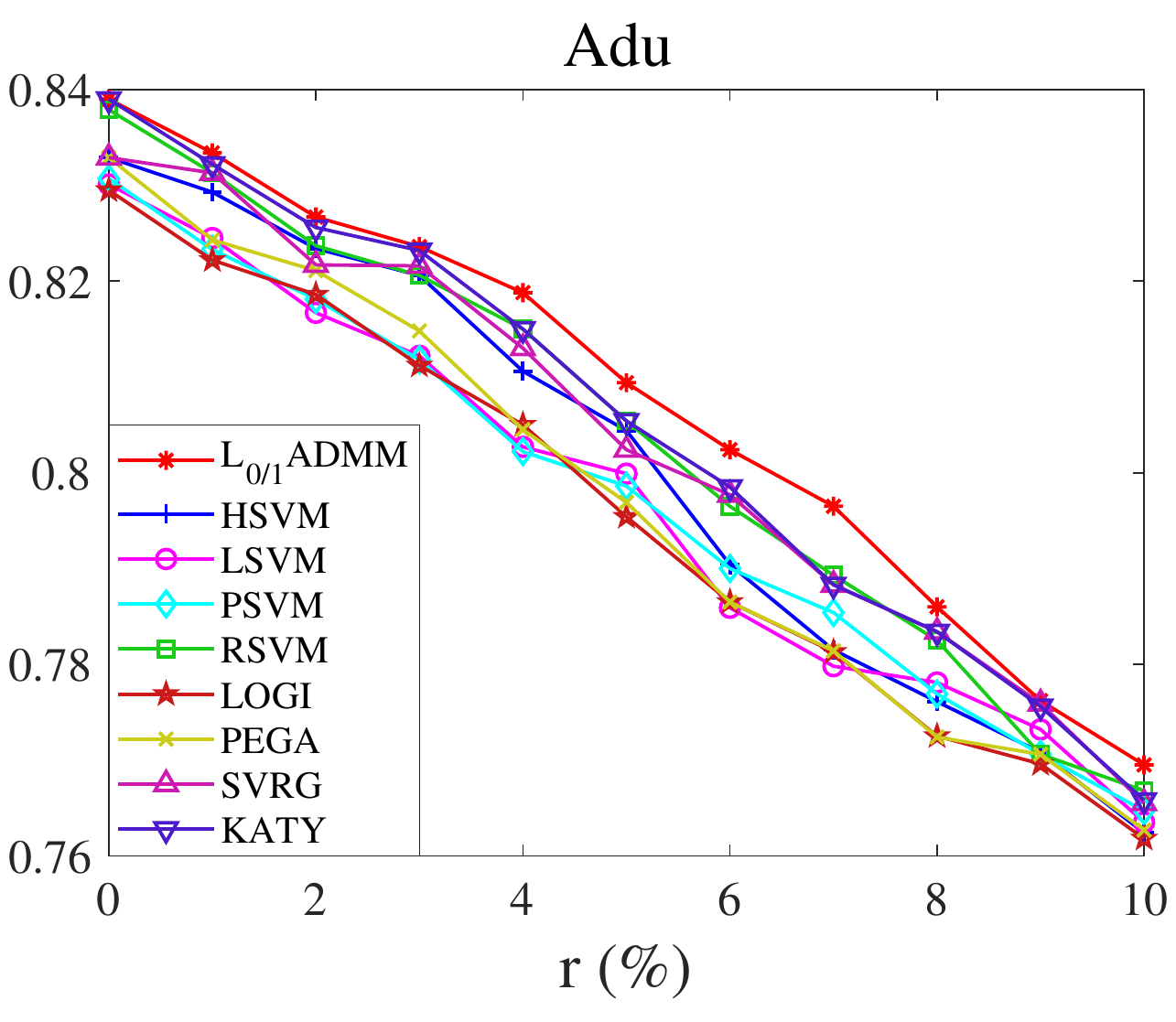}	
\end{subfigure}
 \caption{ \ACC\ vs. $r$  of all solvers for solving six datasets.}
\label{fig:acc}
\end{figure}

\begin{figure}[h]
\centering
\begin{subfigure}{0.24\textwidth}
	\centering
	\includegraphics[height=2.7cm,width=4cm]{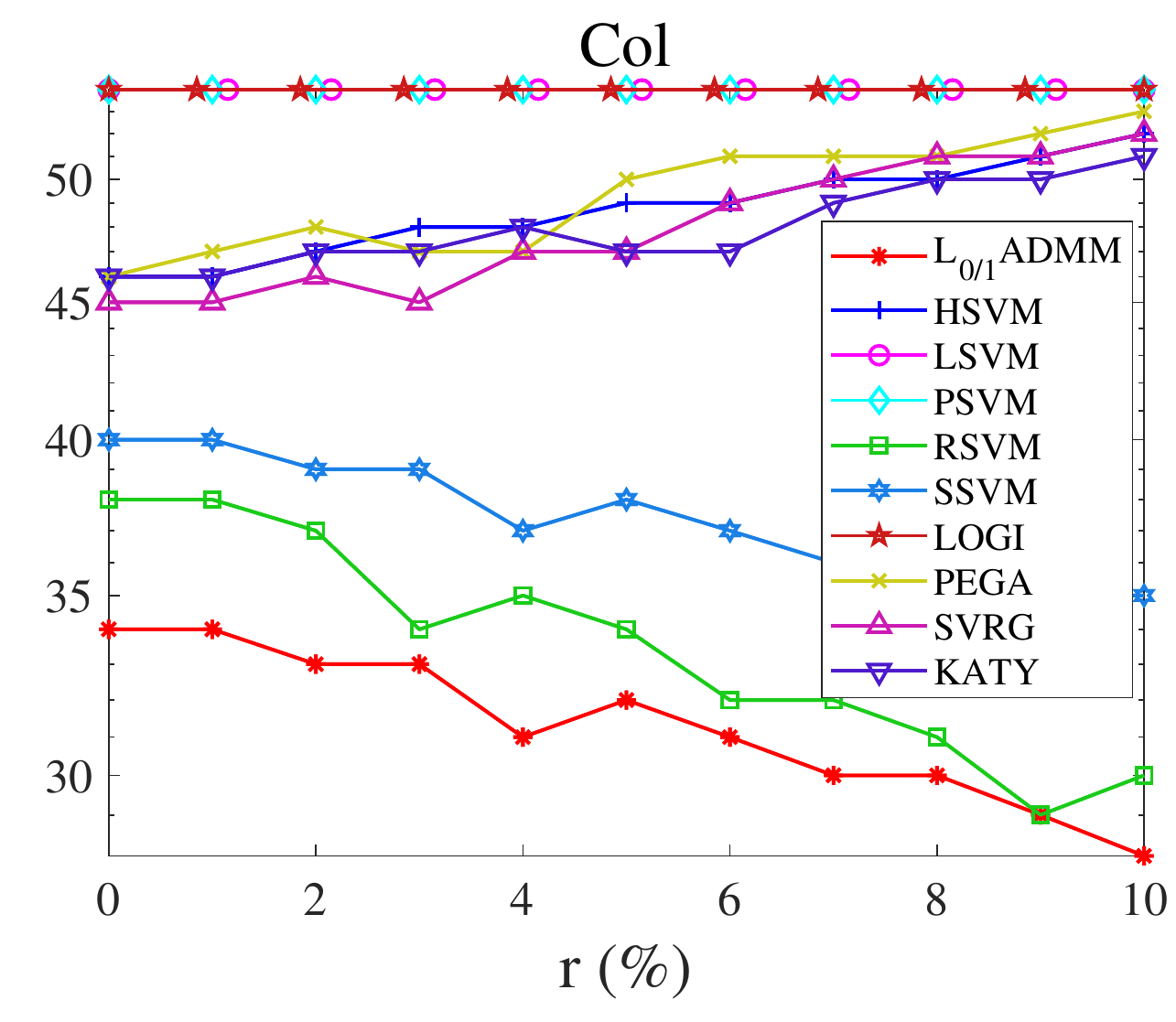}
\end{subfigure}%
\begin{subfigure}{0.24\textwidth}
	\centering
	\includegraphics[height=2.7cm,width=4cm]{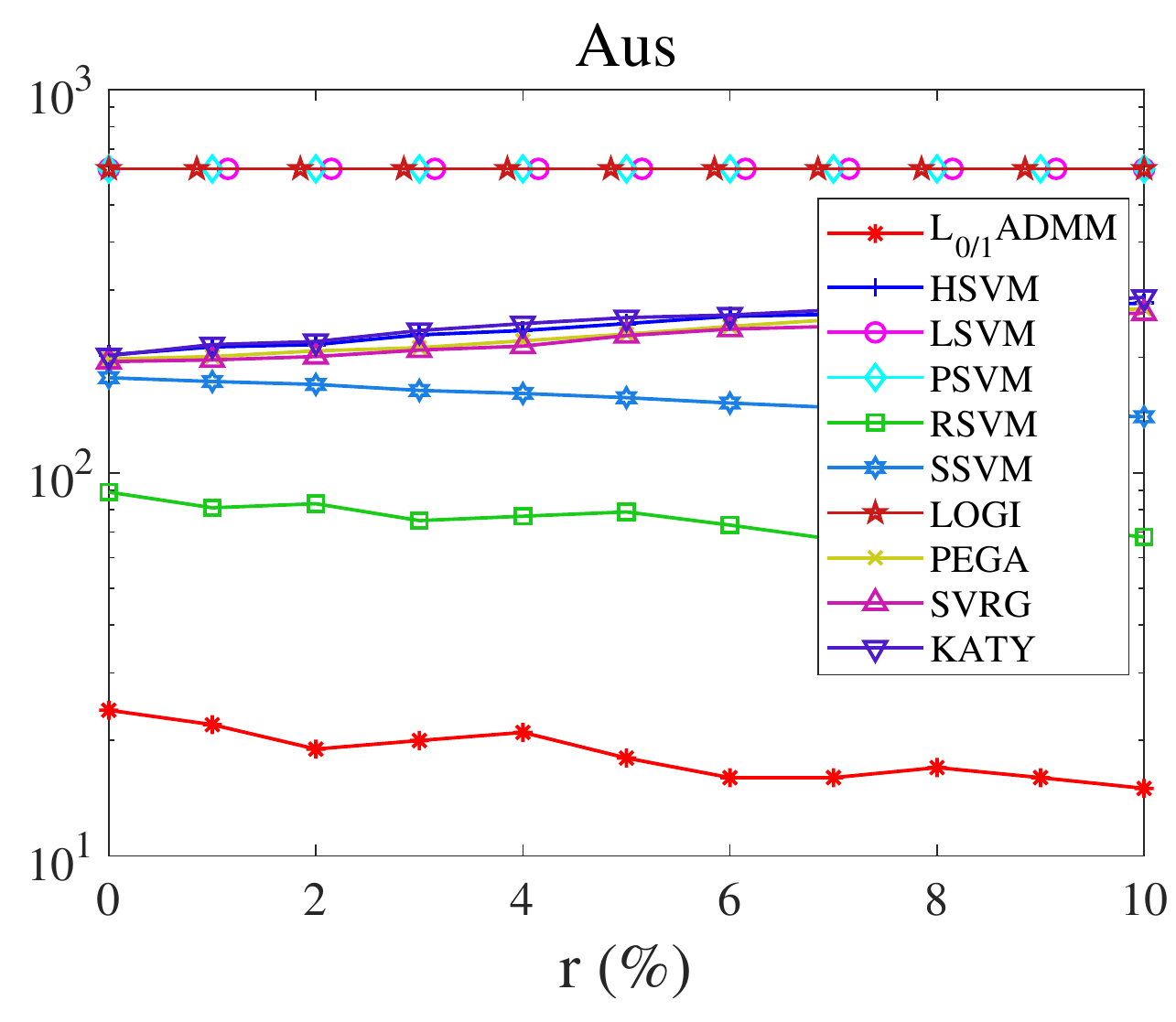}	
\end{subfigure}
\begin{subfigure}{0.24\textwidth}
	\centering
	\includegraphics[height=2.7cm,width=4cm]{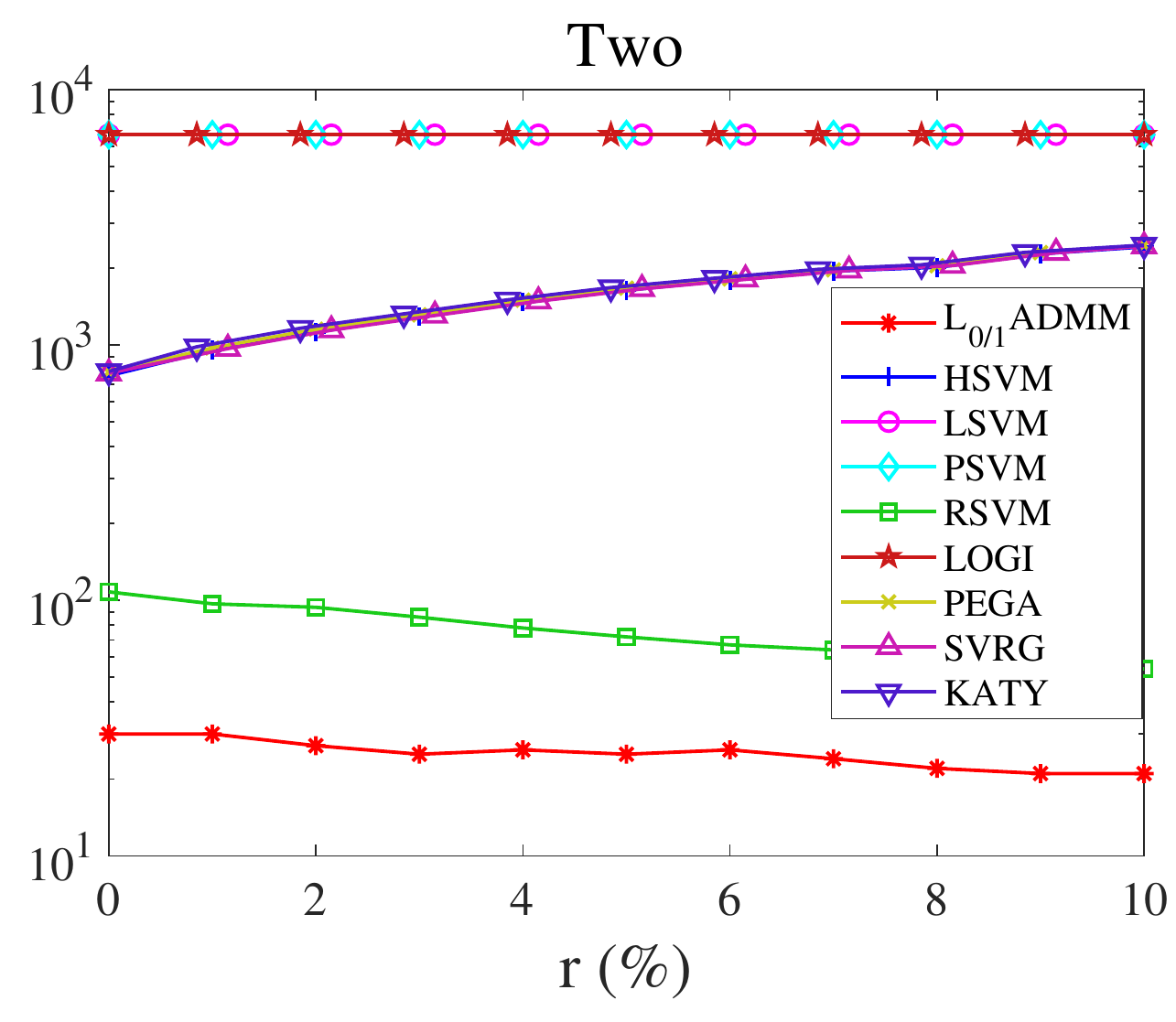}
\end{subfigure}%
\begin{subfigure}{0.24\textwidth}
	\centering
	\includegraphics[height=2.7cm,width=4cm]{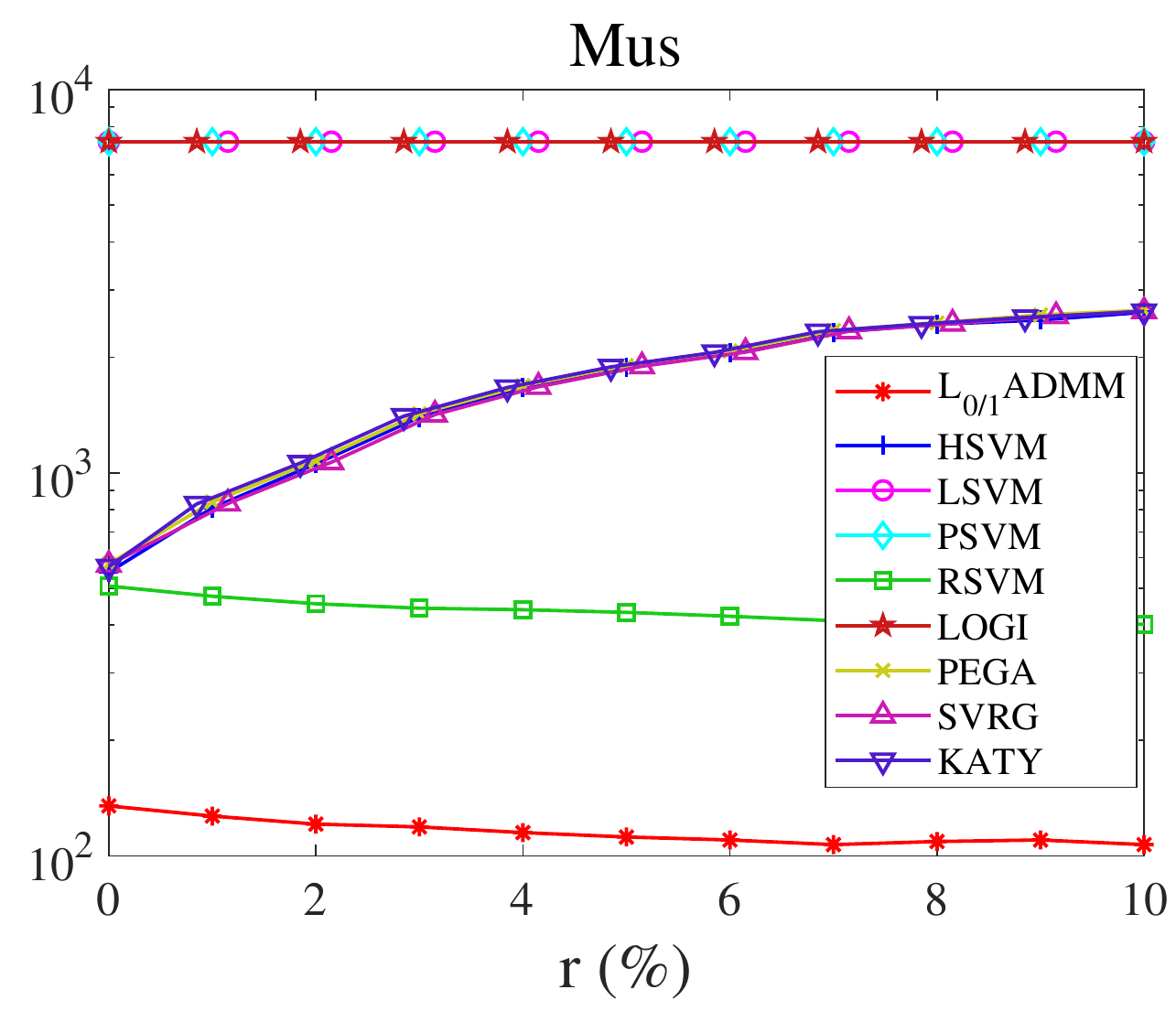}	
\end{subfigure}
\begin{subfigure}{0.24\textwidth}
	\centering
	\includegraphics[height=2.7cm,width=4cm]{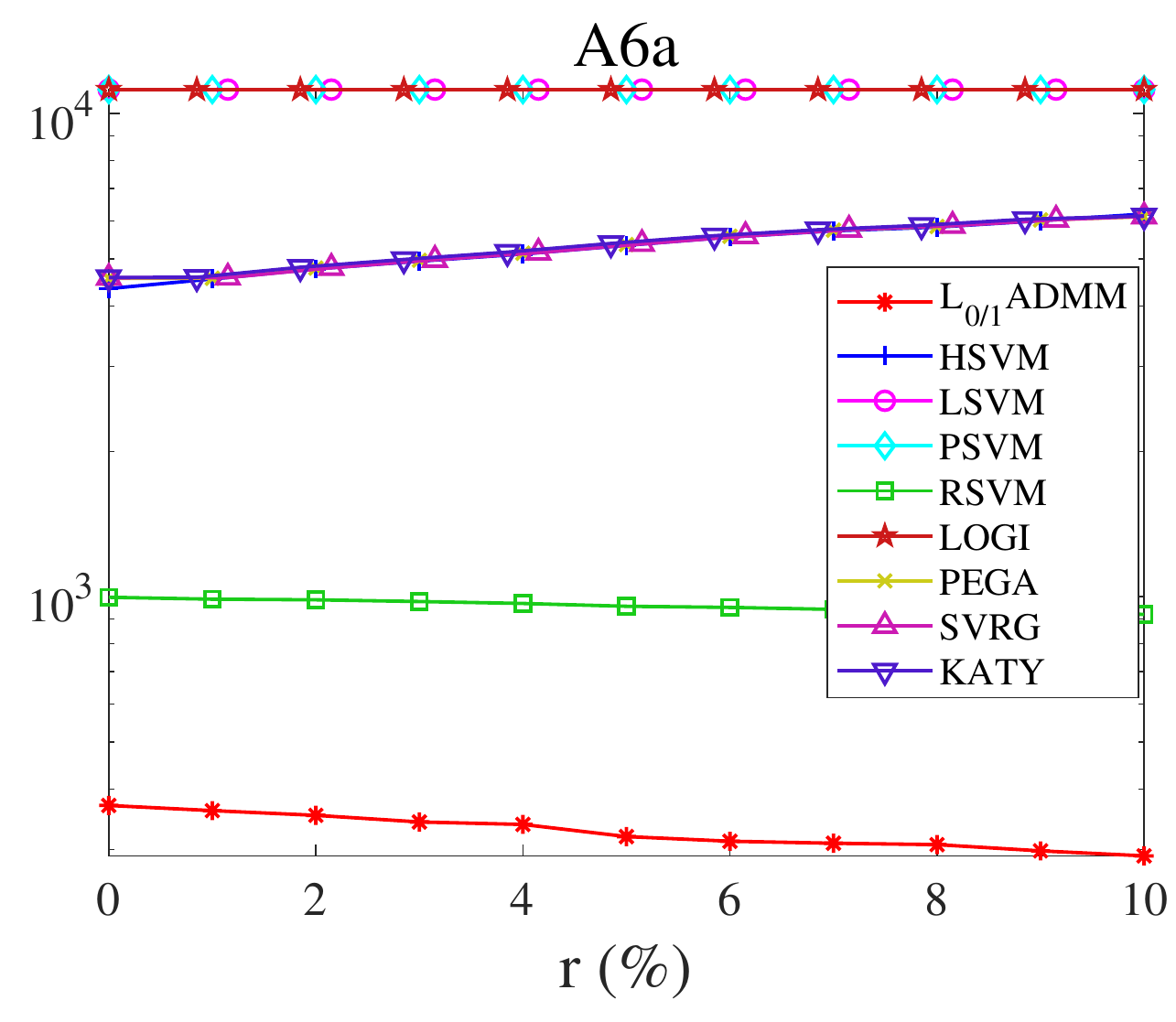}
\end{subfigure}%
\begin{subfigure}{0.24\textwidth}
	\centering
	\includegraphics[height=2.7cm,width=4cm]{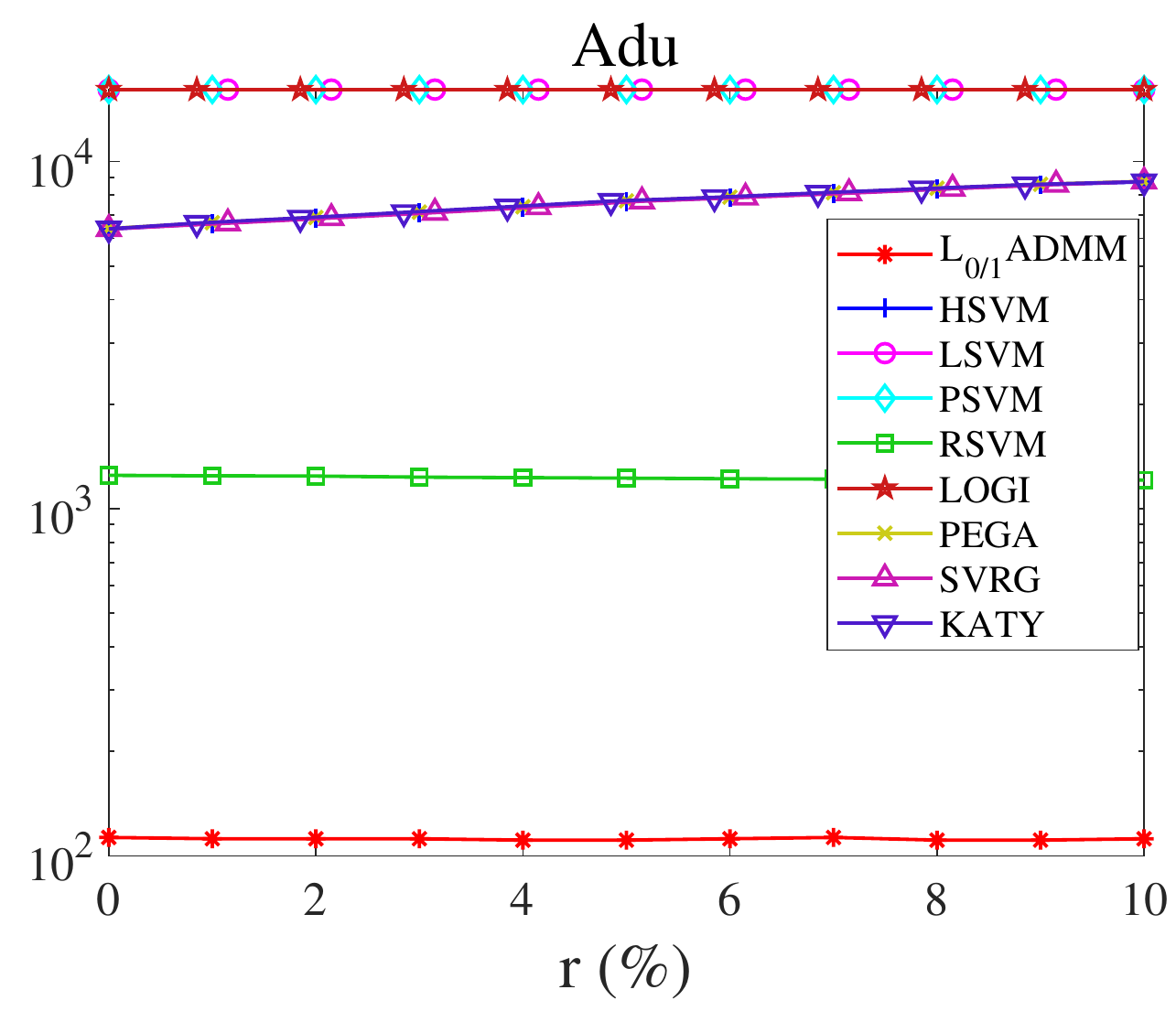}	
\end{subfigure}
 \caption{ \NSV\ vs. $r$  of all solvers for solving six datasets.}
\label{fig:nsv}
\end{figure}

\begin{figure}[h]
\centering
\begin{subfigure}{0.24\textwidth}
	\centering
	\includegraphics[height=2.7cm,width=4cm]{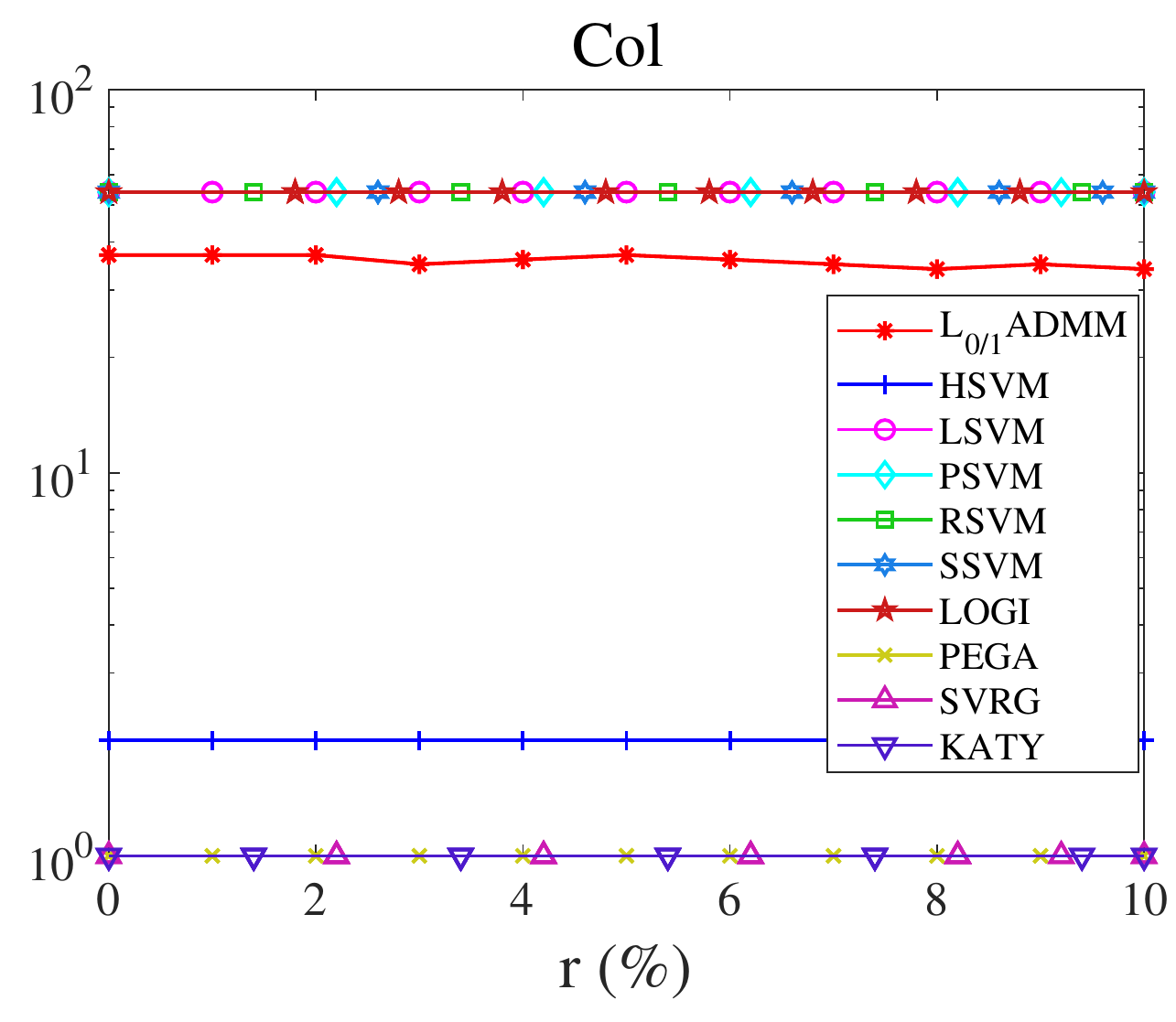}
\end{subfigure}%
\begin{subfigure}{0.24\textwidth}
	\centering
	\includegraphics[height=2.7cm,width=4cm]{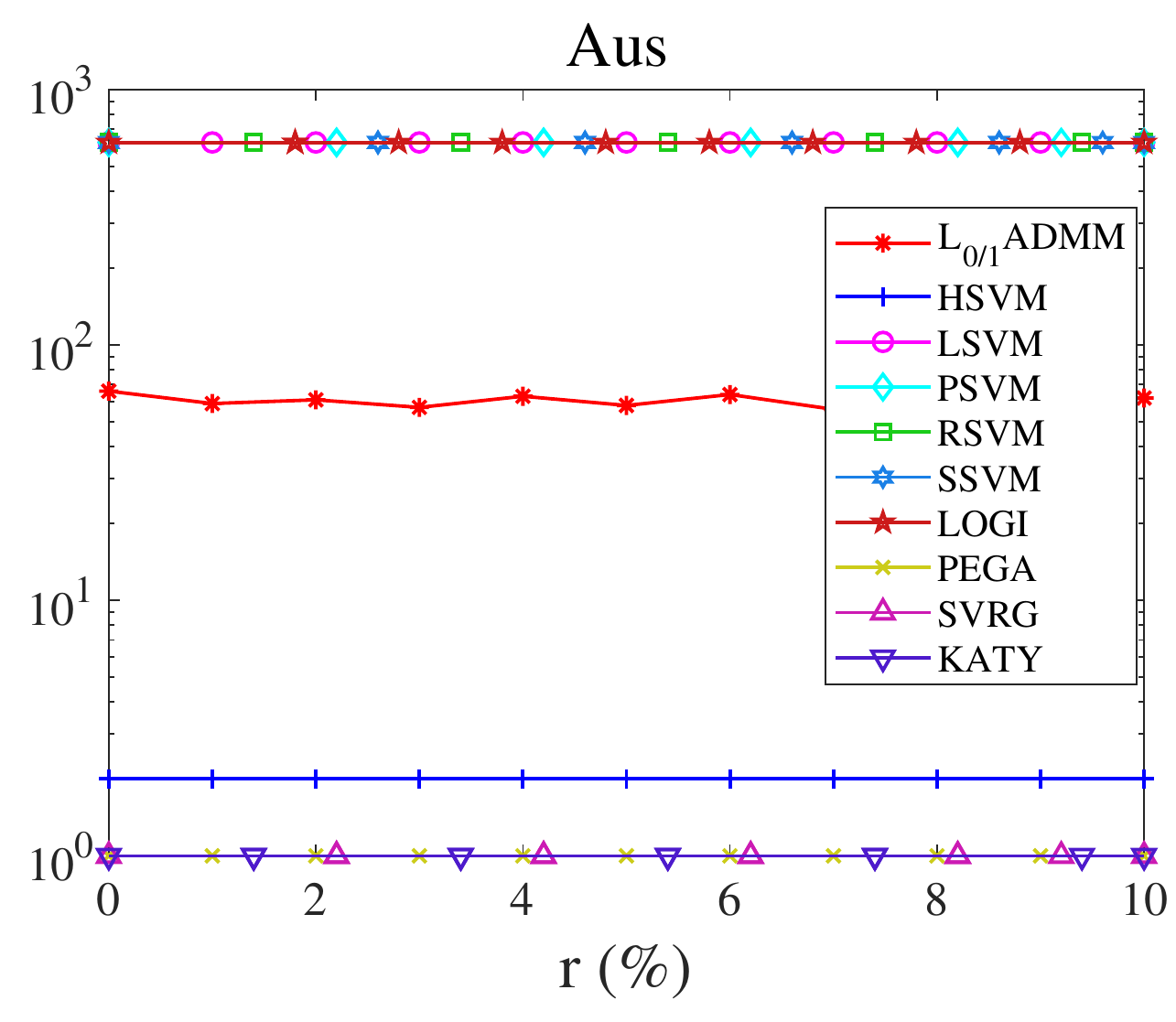}	
\end{subfigure}
\begin{subfigure}{0.24\textwidth}
	\centering
	\includegraphics[height=2.7cm,width=4cm]{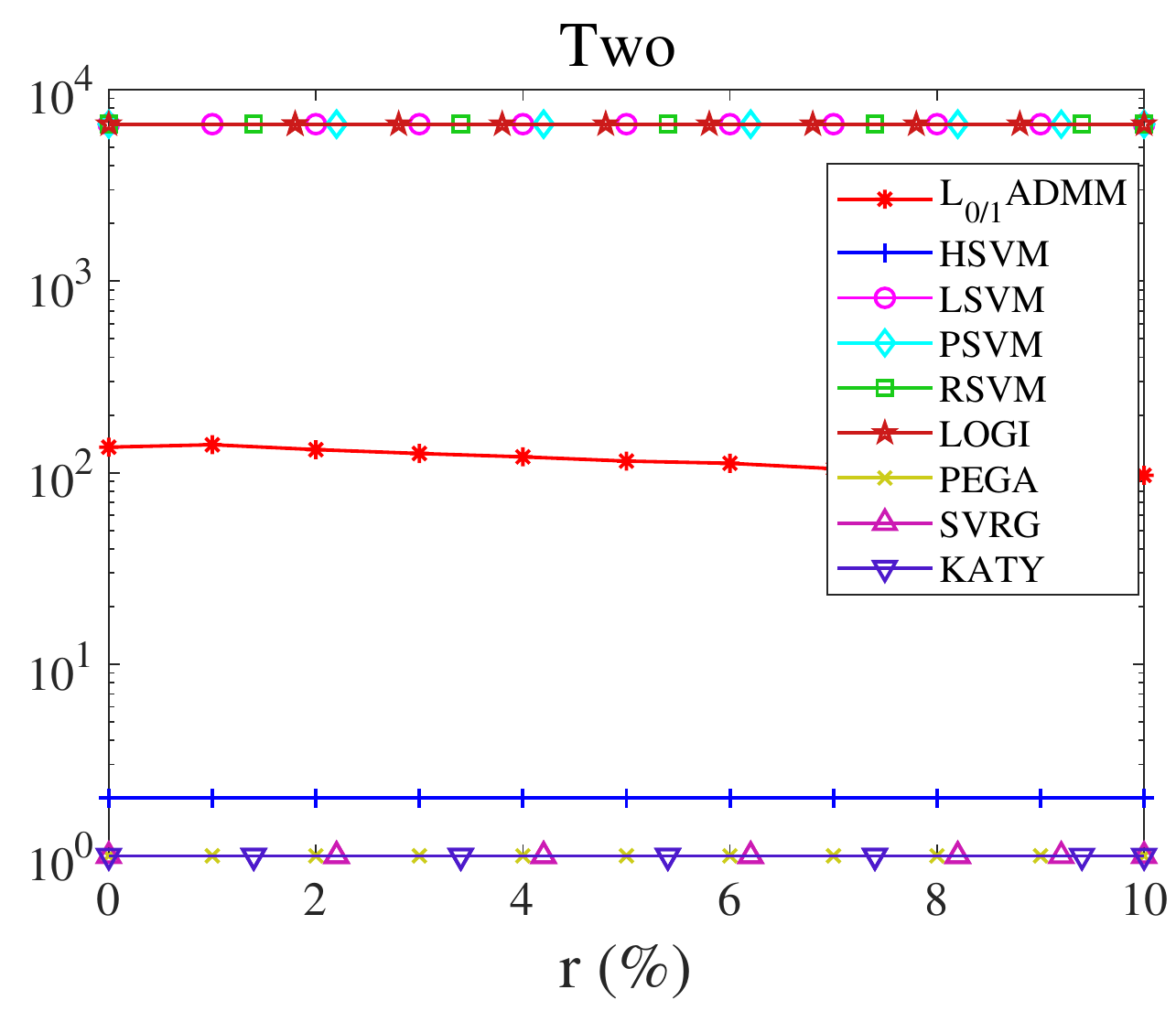}
\end{subfigure}%
\begin{subfigure}{0.24\textwidth}
	\centering
	\includegraphics[height=2.7cm,width=4cm]{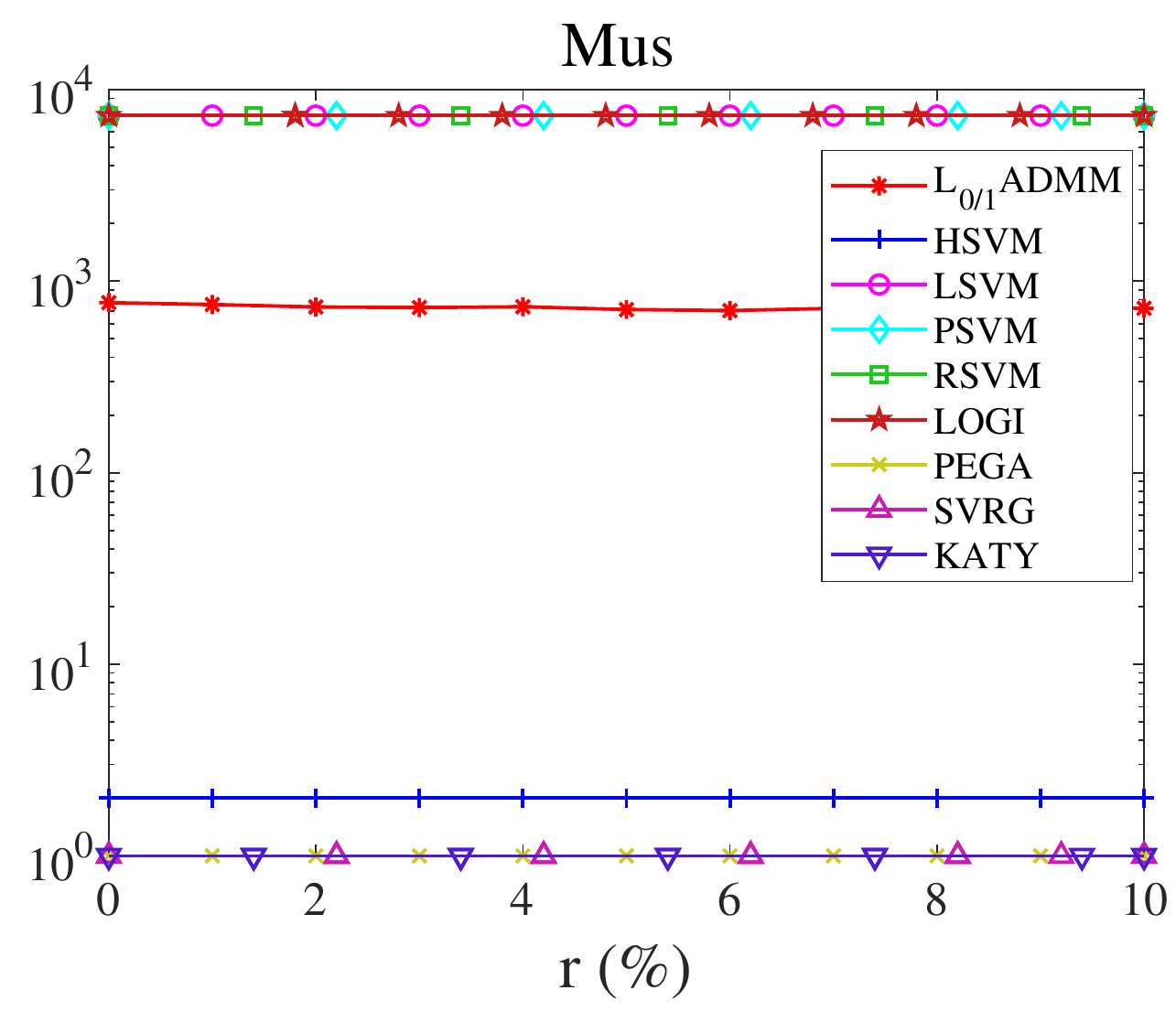}	
\end{subfigure}
\begin{subfigure}{0.24\textwidth}
	\centering
	\includegraphics[height=2.7cm,width=4cm]{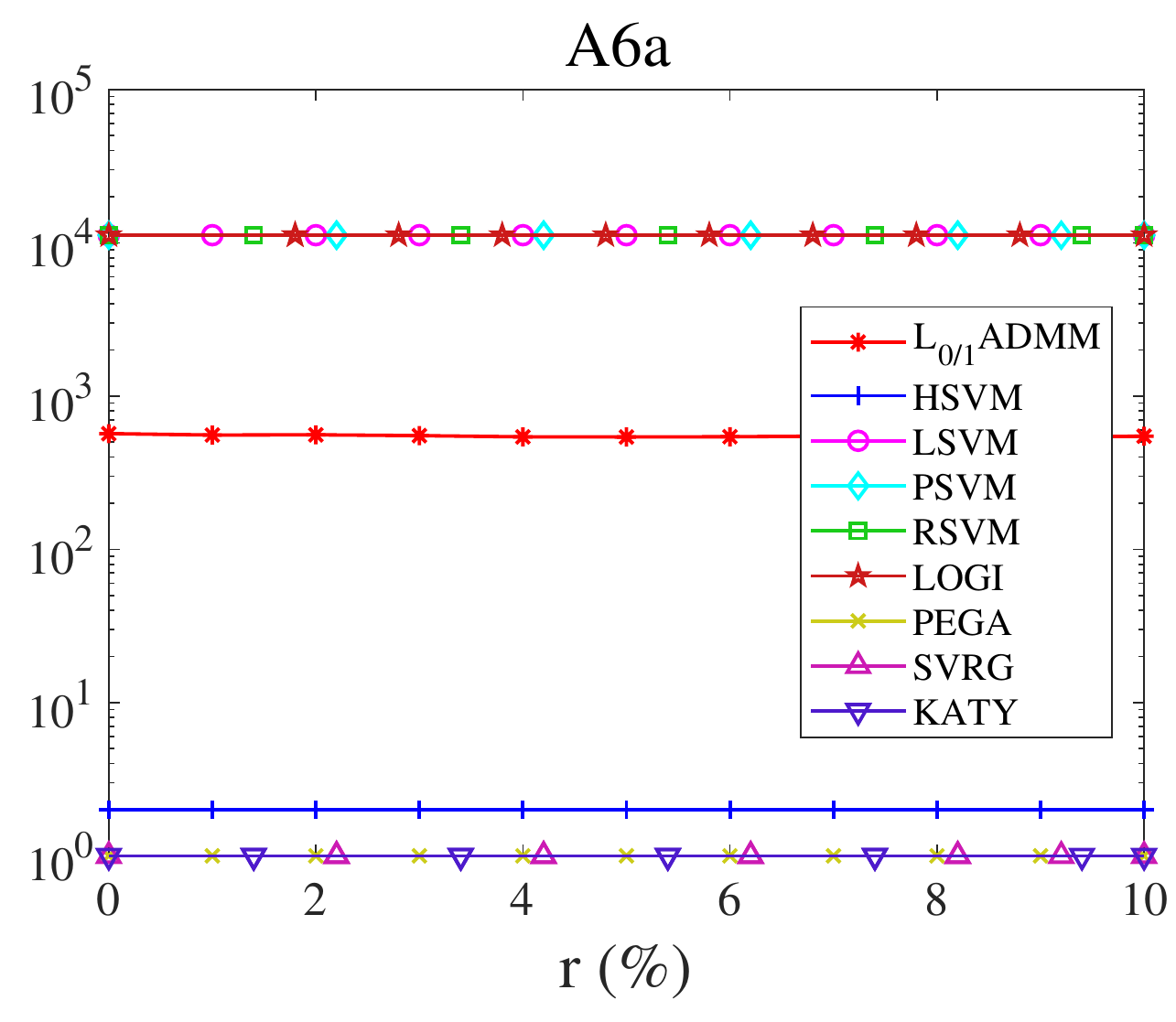}
\end{subfigure}%
\begin{subfigure}{0.24\textwidth}
	\centering
	\includegraphics[height=2.7cm,width=4cm]{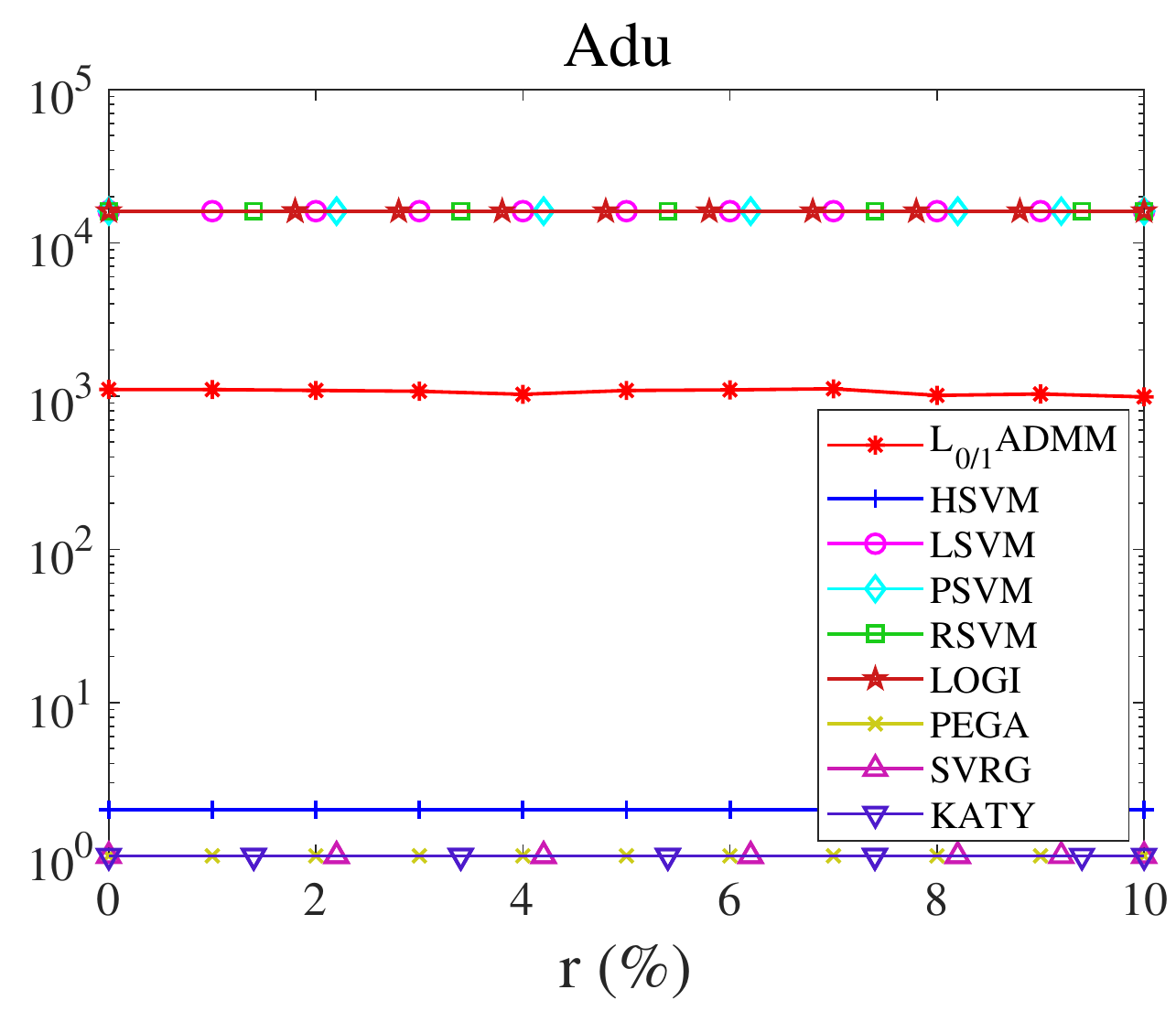}	
\end{subfigure}
 \caption{ {\SWS\ vs. $r$  of all solvers for solving six datasets.}}
\label{fig:sws}
\end{figure}

\begin{figure}[h]
\centering
\begin{subfigure}{0.24\textwidth}
	\centering
	\includegraphics[height=2.7cm,width=4cm]{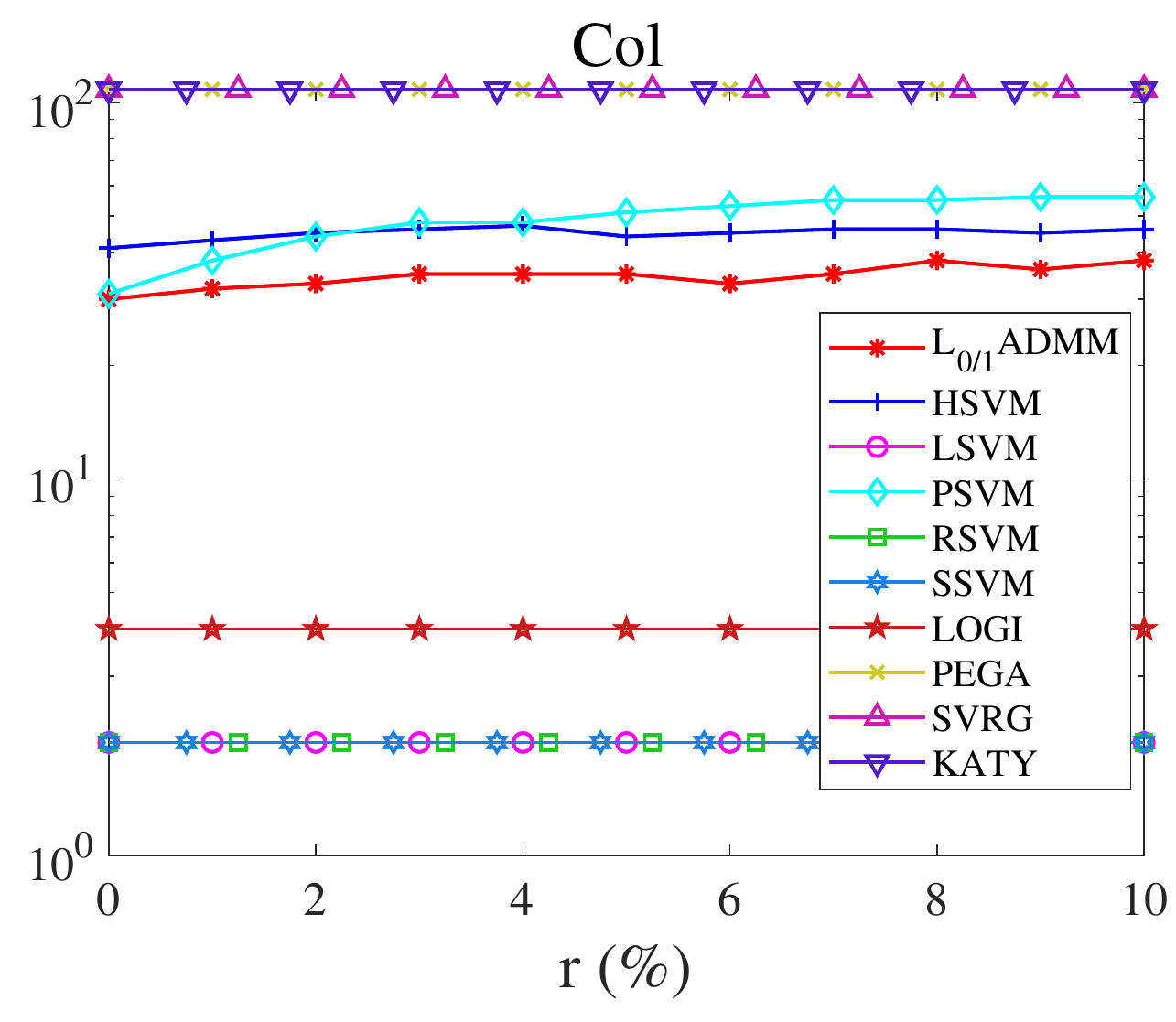}
\end{subfigure}%
\begin{subfigure}{0.24\textwidth}
	\centering
	\includegraphics[height=2.7cm,width=4cm]{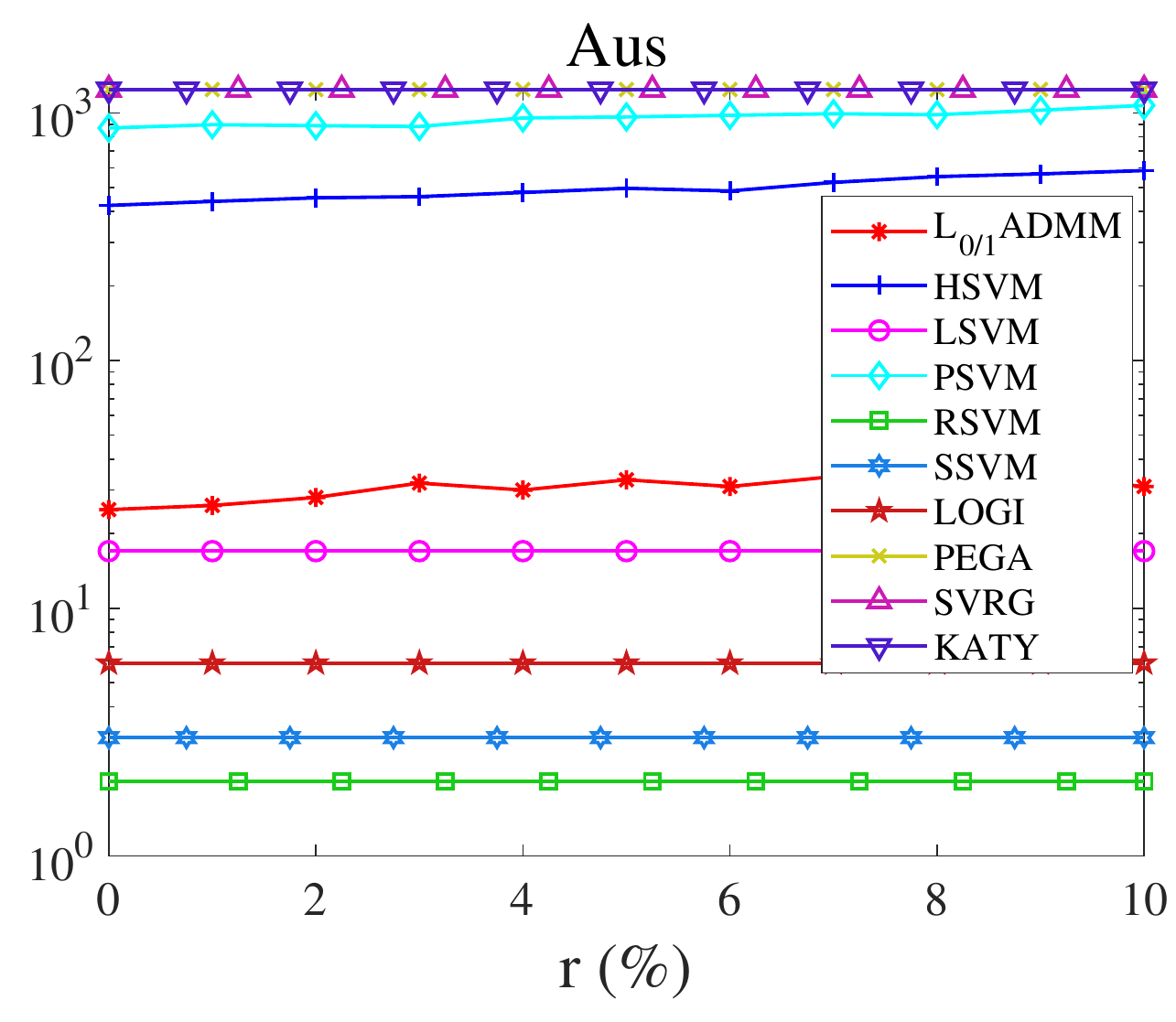}	
\end{subfigure}
\begin{subfigure}{0.24\textwidth}
	\centering
	\includegraphics[height=2.7cm,width=4cm]{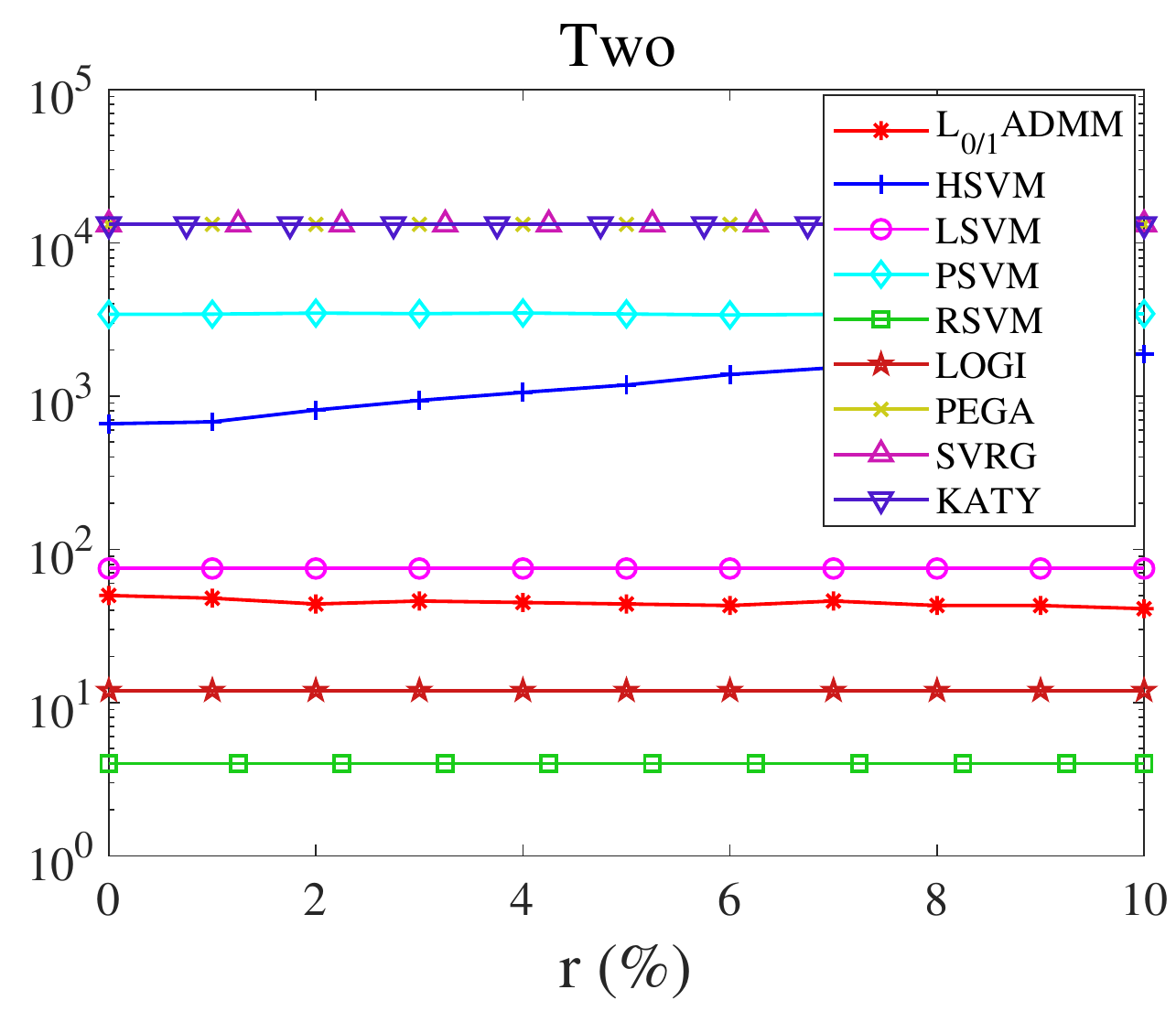}
\end{subfigure}%
\begin{subfigure}{0.24\textwidth}
	\centering
	\includegraphics[height=2.7cm,width=4cm]{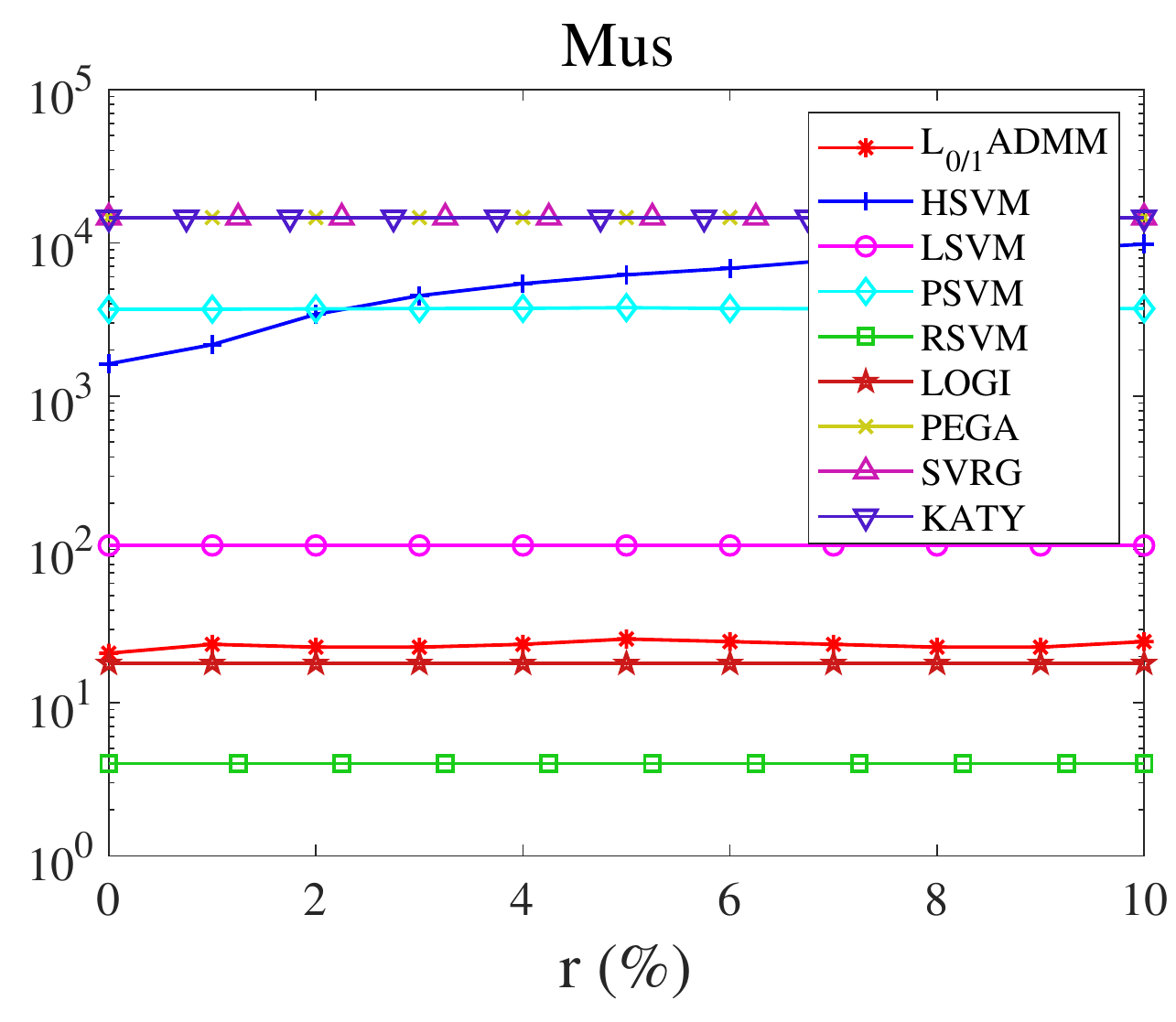}	
\end{subfigure}
\begin{subfigure}{0.24\textwidth}
	\centering
\includegraphics[height=2.7cm,width=4cm]{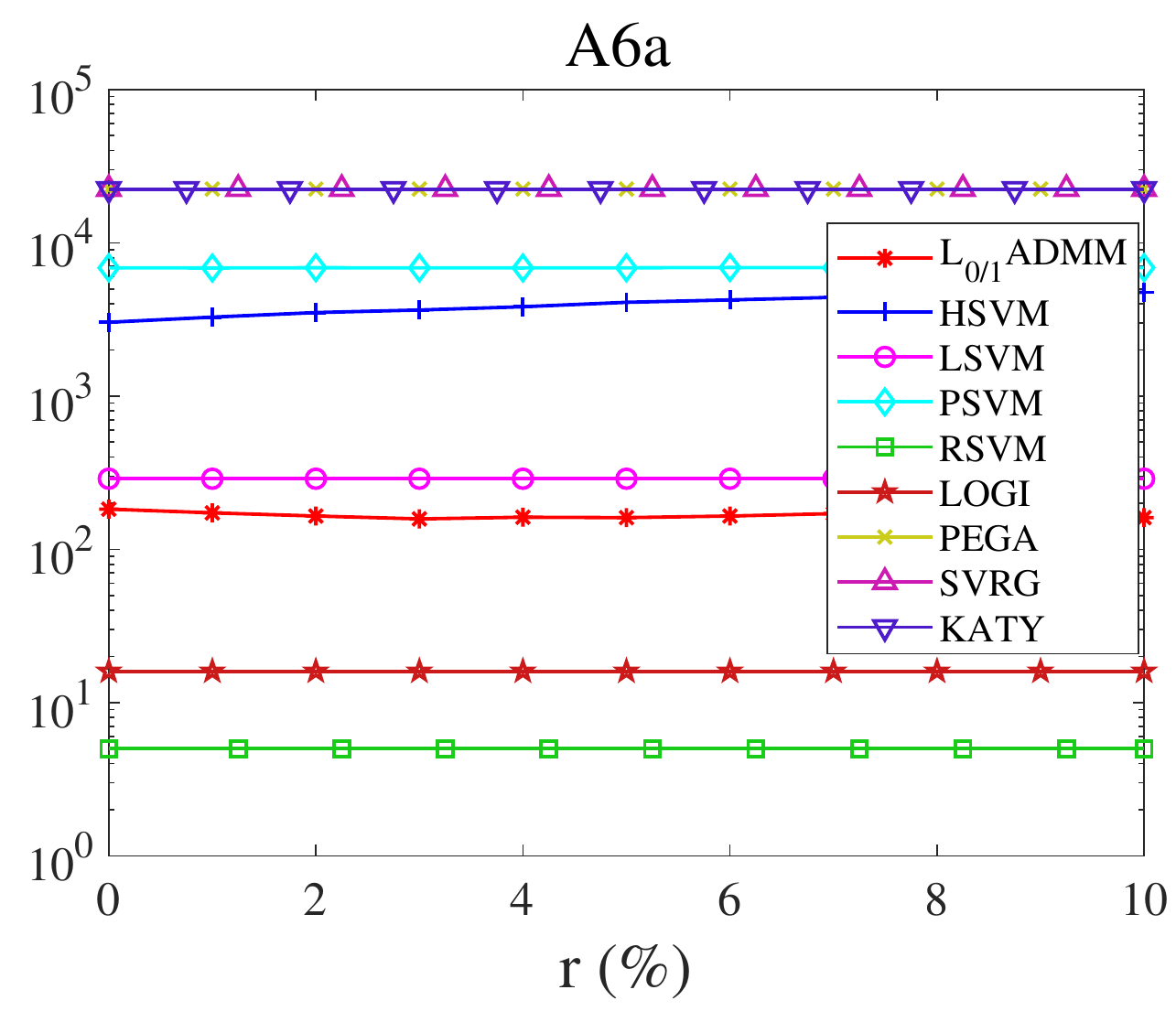}
\end{subfigure}%
\begin{subfigure}{0.24\textwidth}
	\centering	
\includegraphics[height=2.7cm,width=4cm]{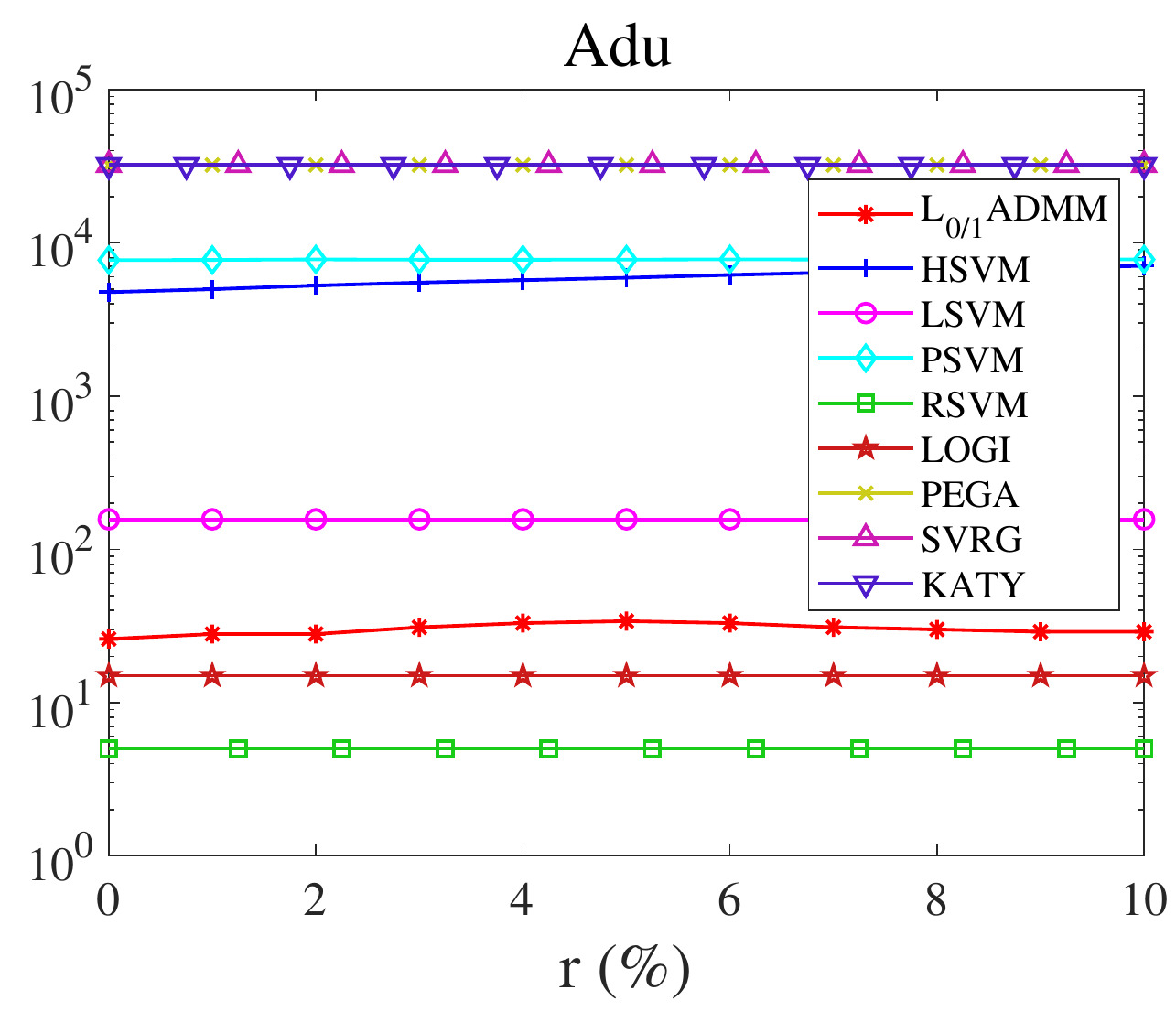}
\end{subfigure}
 \caption{ {\TNI\ vs. $r$  of all solvers for solving six datasets.}}
\label{fig:tni}
\end{figure}

\begin{figure}[h]
\centering
\begin{subfigure}{0.24\textwidth}
	\centering
	\includegraphics[height=2.7cm,width=4cm]{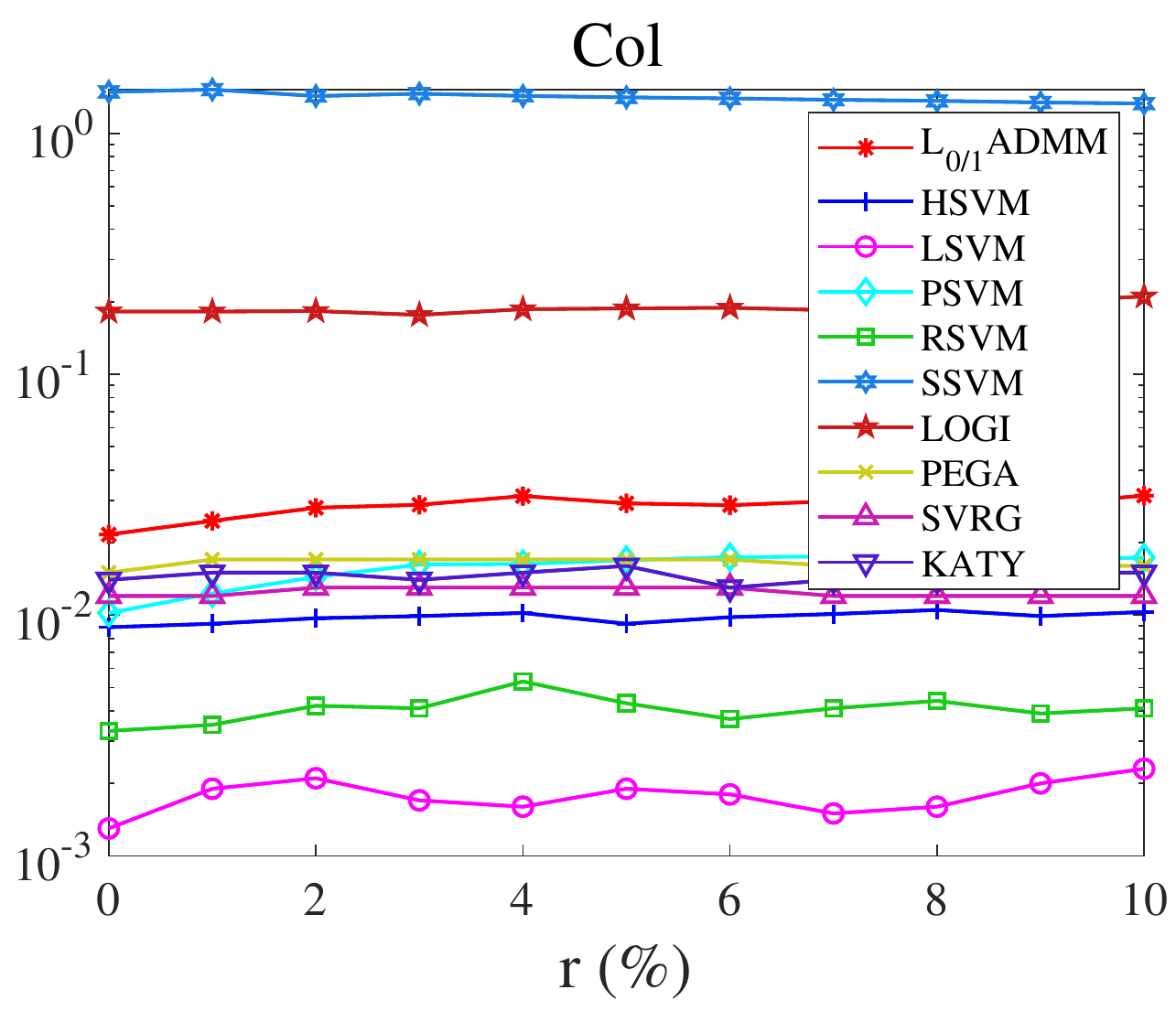}
\end{subfigure}%
\begin{subfigure}{0.24\textwidth}
	\centering
	\includegraphics[height=2.7cm,width=4cm]{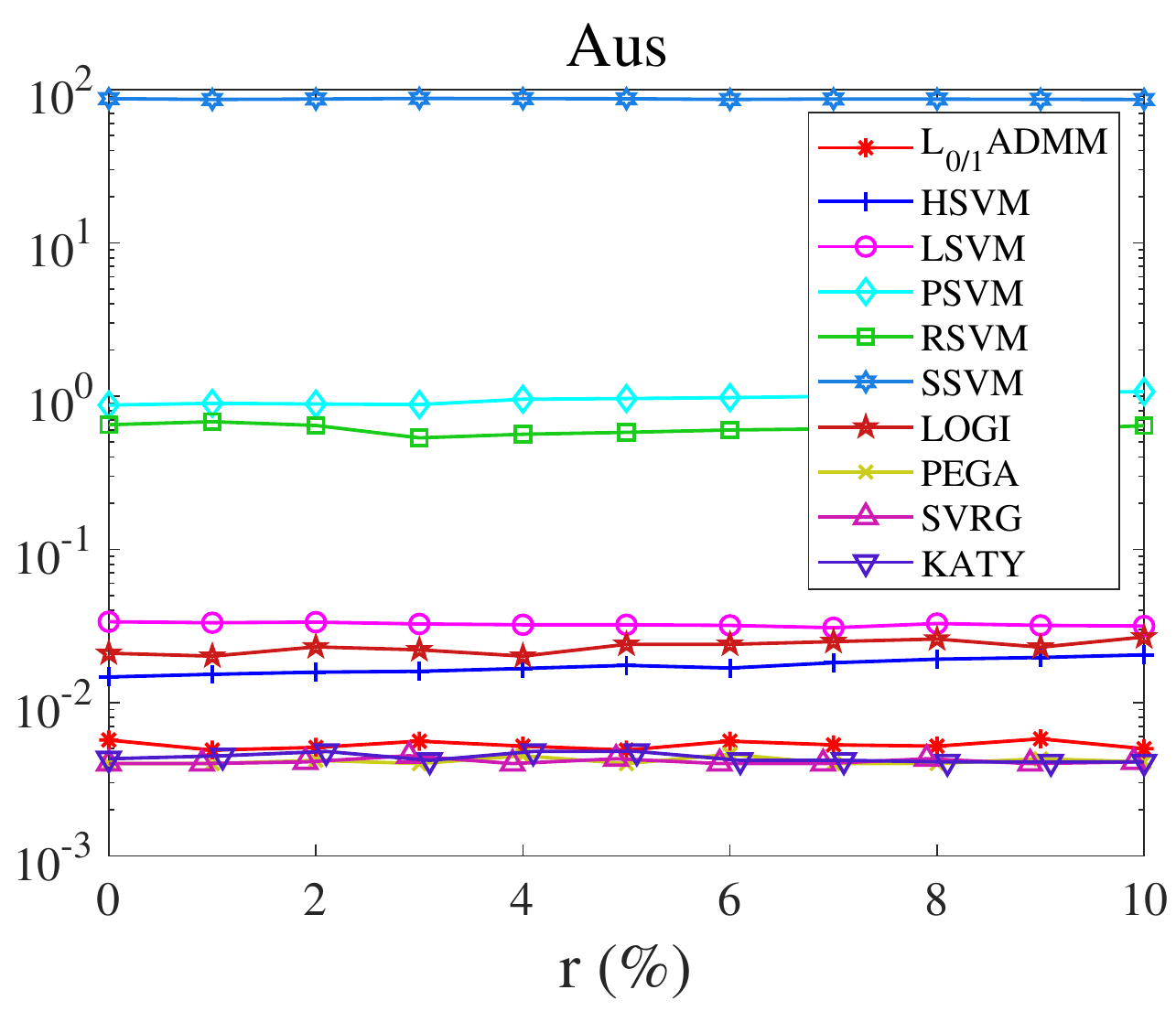}	
\end{subfigure}
\begin{subfigure}{0.24\textwidth}
	\centering
	\includegraphics[height=2.7cm,width=4cm]{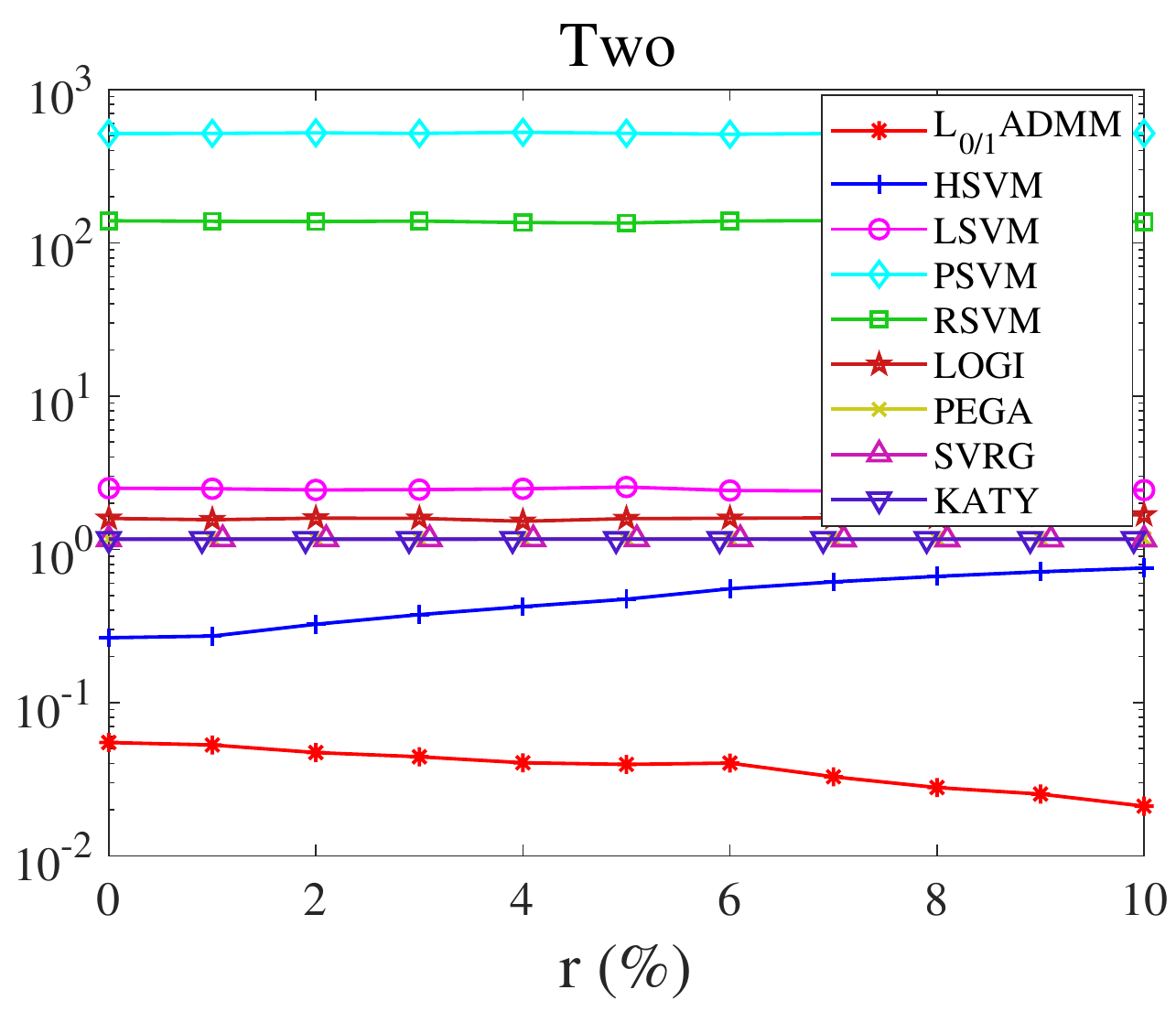}
\end{subfigure}%
\begin{subfigure}{0.24\textwidth}
	\centering
	\includegraphics[height=2.7cm,width=4cm]{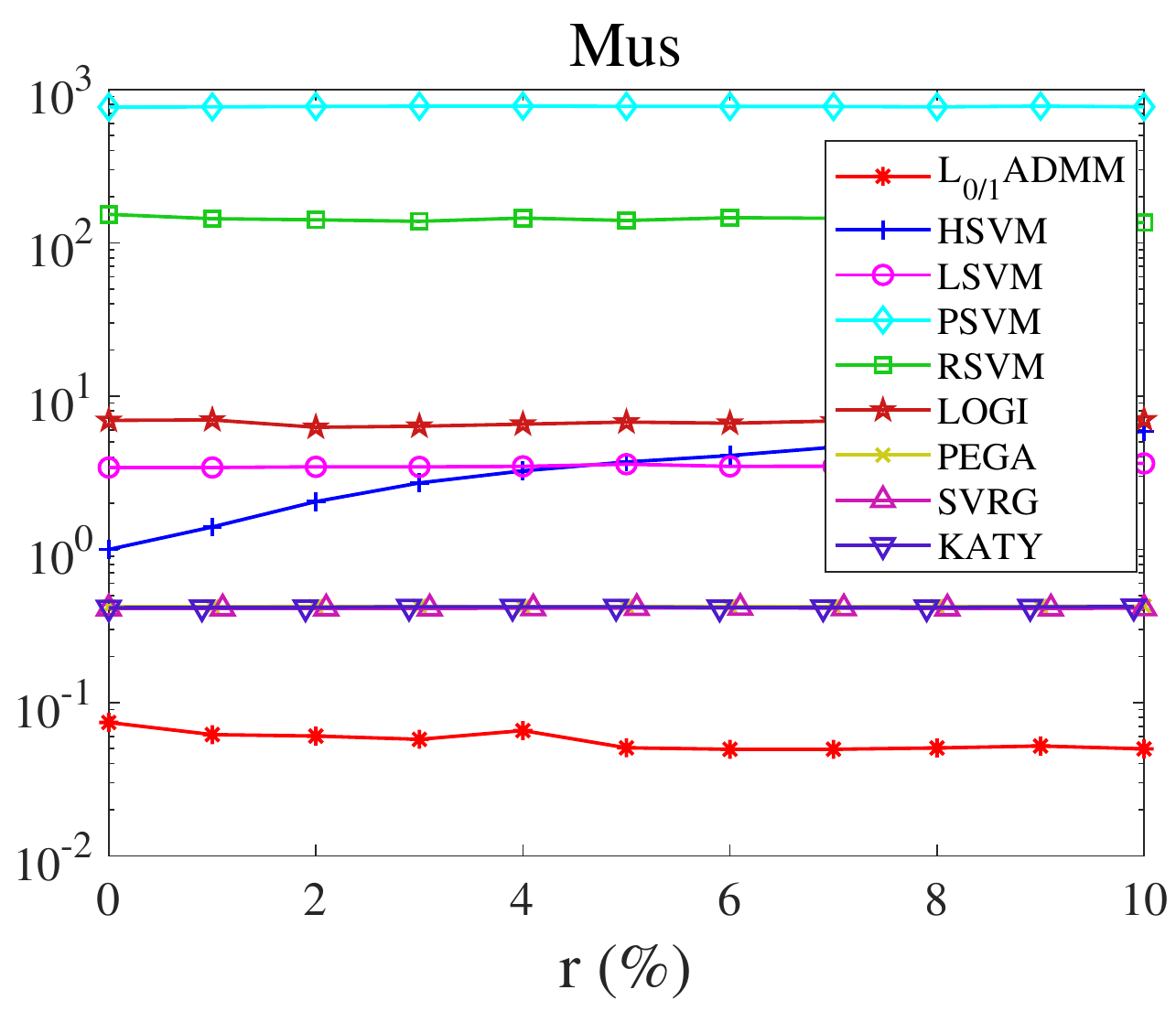}	
\end{subfigure}
\begin{subfigure}{0.24\textwidth}
	\centering
	\includegraphics[height=2.7cm,width=4cm]{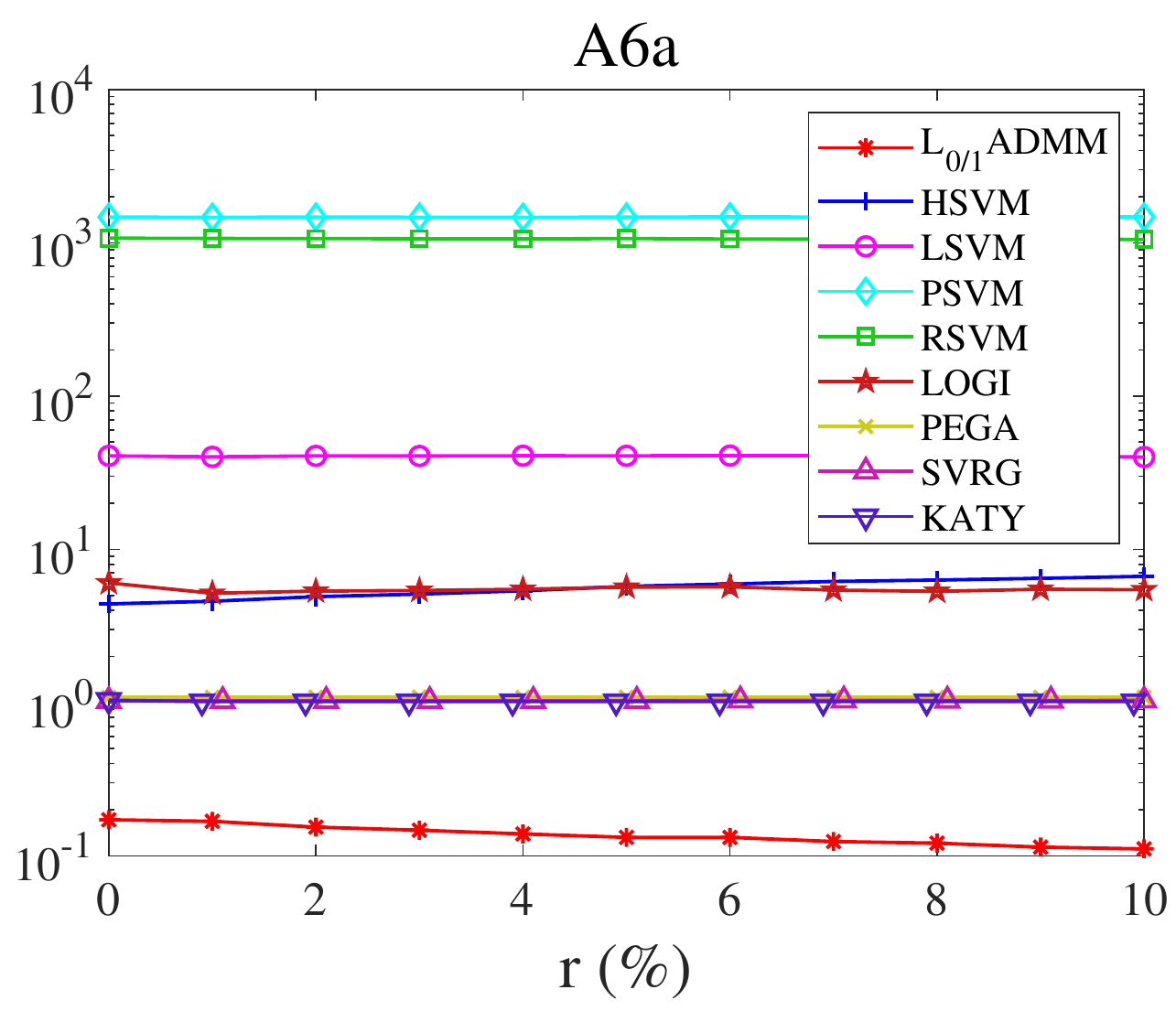}
\end{subfigure}%
\begin{subfigure}{0.24\textwidth}
	\centering
	\includegraphics[height=2.7cm,width=4cm]{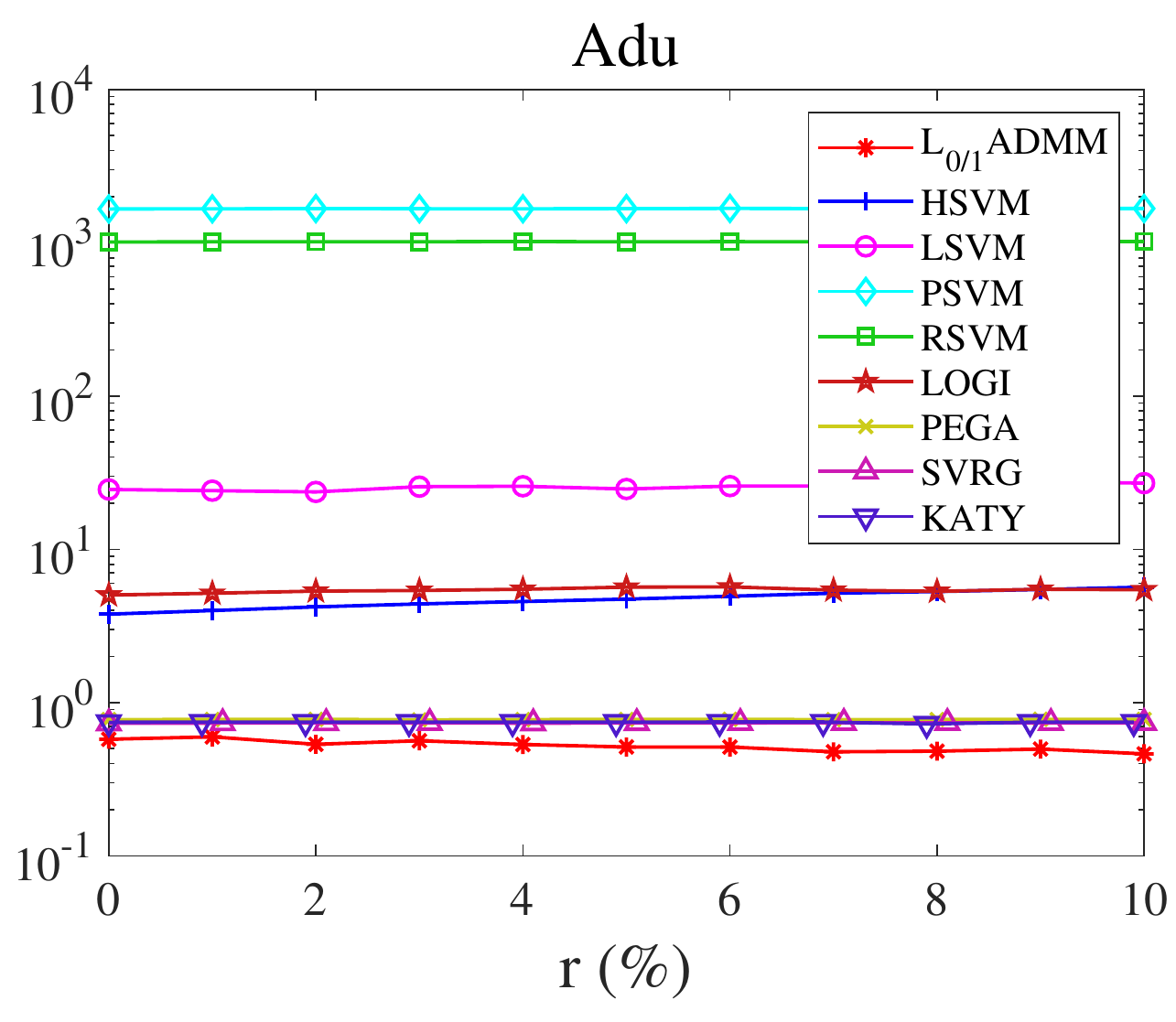}	
\end{subfigure}
 \caption{ \CPU\ vs. $r$  of all solvers for solving six datasets.}
\label{fig:time}
\end{figure}

{\bf (g) Real data with outliers.}  Finally, we would like to see the robustness of each solver to the outliers for real datasets in \cref{ex:real-data-outlier}. Again we alter the flapping ratio $r$ from $\{0.01, 0.02,\cdots,0.1\}$.
 It is shown in \cref{Table-LU-4} that  \sosvm\ takes too long time for datasets: \texttt{two}, \texttt{mus},  \texttt{adu} and \texttt{a6a}. Therefore, its results related to these datasets are omitted.  All lines of \ACC\ shown in \cref{fig:acc}  decline  with  $r$ ascending, and \ADMML\ achieves the highest \ACC.  As for \NSV\ in \cref{fig:nsv}, {\ssvm}, {\psvm} and {\logi} always treat all samples as support vectors.  {\hsvm, \sosvm, \pega, \svrg~and~\katy\ increase
\NSV\ with the rising of $r$.}  Lines from \ADMML\ and {\rsvm}  either decline or stabilize at a level with the rising of $r$, which means they are quite robust to $r$, namely robust to the outliers. What is more, \ADMML\ always renders the fewest \NSV.  { As for \SWS~in \cref{fig:sws},  with the rising of $r$, \ADMML\ stabilizes at a level for all datasets. As for \TNI~in
\cref{fig:tni}, all algorithms no big difference with the ascending of $r$ except for \hsvm.} For the computational speed, as demonstrated in \cref{fig:time}, \ADMML\ outperforms the others for all datasets except for  \texttt{col} and \texttt{aus} which have a very small size.

  \section{conclusion}
In this paper, we have explored an ideal soft-margin SVM model: $L_{0/1}$-SVM, which well captures the nature of the binary classification and guarantees a fewer number of support vectors than the other soft-margin SVM models. Despite the discreteness of the $L_{0/1}$-SVM, the establishment of the optimality theory, associated with the P-stationary point, makes it tractable numerically.  Based on the idea of  $L_{0/1}$  support vectors inspired by the P-stationary point, a working set was cast and integrated into the proximal ADMM, which tremendously speeds up the whole computation and reduces the number of support vectors. Consequently, the proposed method performed exceptionally well with fewer support vectors and faster computational speed, especially for datasets on large scales.

We feel that the established methodology and techniques might be able to extend to process the nonlinear kernel SVMs \cite{BW2019, YY2007,YZ2017} and problems from perception learning \cite{LL2007} and deep learning \cite{IB2016}. We leave these as future research.


 \section*{Acknowledgements}

{The authors would like to thank the associate editor and three anonymous
referees for their constructive comments, which have significantly improved
the quality of the paper.} This work is supported by the National Natural Science Foundation of China (11971052, 11926348-9, 61866010, 11871183), and the Natural Science Foundation of Hainan Province (120RC449).


 \vspace{-15 mm}

\begin{IEEEbiography}[{\includegraphics[width=1in , clip,keepaspectratio]{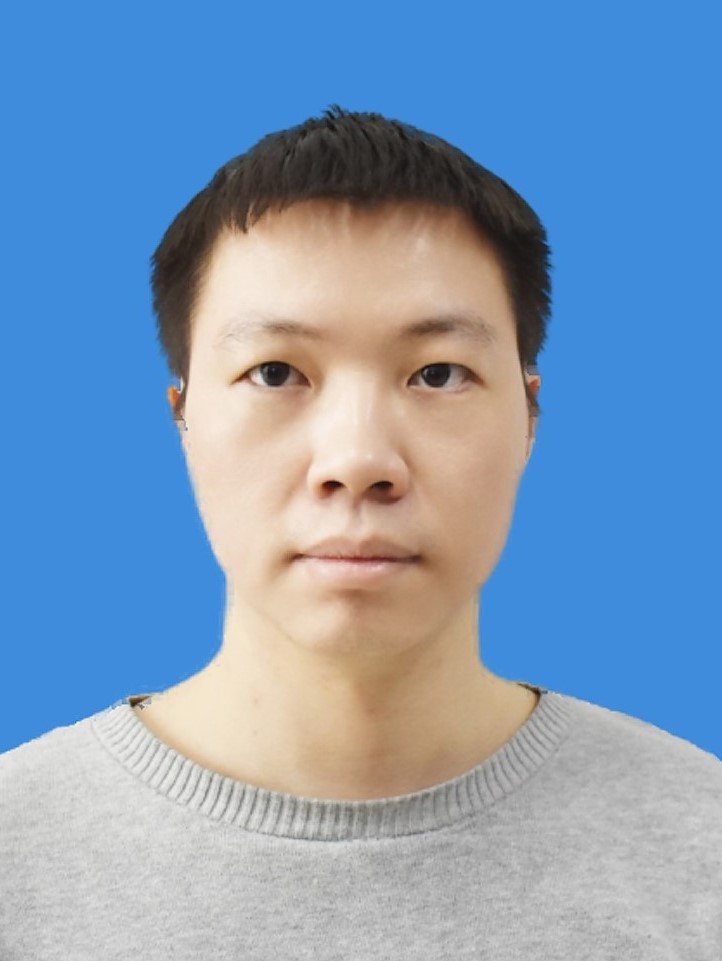}}]{Huajun Wang} received his M.Sc. degree in Department of Mathematics from Guilin University of Electronic Technology, China, in 2017.  He is currently a Ph.D. candidate of Department of Applied Mathematics at the Beijing Jiaotong University, China. His current research interests include large-scale classification optimization problems,  machine learning, 0-1 loss optimization and numerical computing.
\end{IEEEbiography} \vspace{-15mm}

\begin{IEEEbiography}[{\includegraphics[width=1in, clip,keepaspectratio]{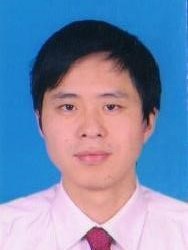}}]{Yuanhai Shao} received his B.Sc. degree in College of Mathematics from Jilin University, and received Ph.D. degree in College of Science from China Agricultural University, China, in 2006 and 2011, respectively. Currently, he is a professor at the Management School, Hainan University.
 His research interests include optimization methods, machine learning, and data mining. He has published over 100 refereed papers.
\end{IEEEbiography} \vspace{-15mm}

\begin{IEEEbiography}[{\includegraphics[width=1in,clip,keepaspectratio]{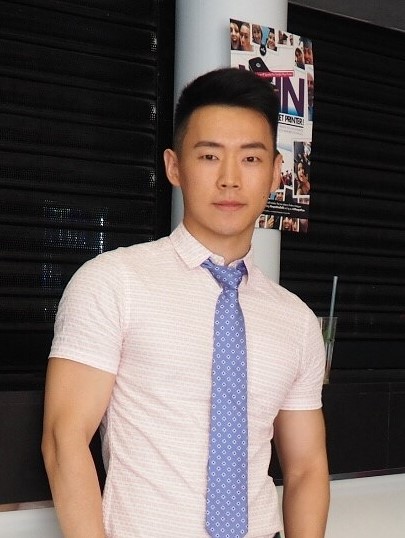}}]{Shenglong Zhou} received the B.Sc. degree in information and computing science in 2011 and the M.Sc. degree in operational research in 2014 from Beijing Jiaotong University, China,  and the Ph.D. degree in operational research in 2018 from the University of Southampton, the United Kingdom, where he was the Research Fellow from 2017 to 2019 and is currently a Teaching Fellow.  His research interests  include the theory and methods of optimization in the fields of sparse, low-rank matrix and bilevel optimization.
\end{IEEEbiography}\vspace{-15mm}

\begin{IEEEbiography}[{\includegraphics[width=1in,clip,keepaspectratio]{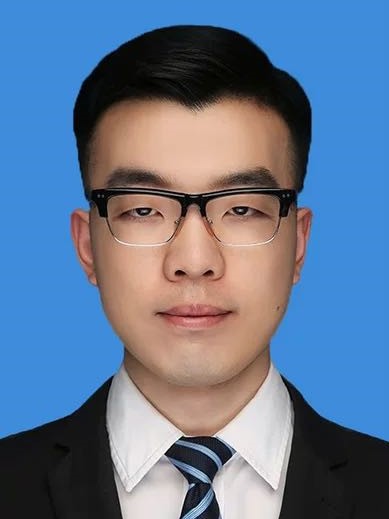}}]{Ce Zhang} received the B.Sc. degree in information and computing science and the M.Sc. degree in operational research from Beijing Jiaotong University, Beijing, China, in 2016 and 2019, respectively. He research interests include optimization methods, machine learning and applications in data and image processing.
\end{IEEEbiography}\vspace{-15mm}

\begin{IEEEbiography}[{\includegraphics[width=1in,clip,keepaspectratio]{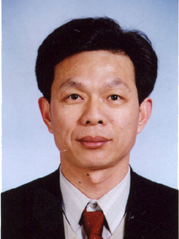}}]{Naihua Xiu} received the B.Sc. degree in mathematics from Hebei Normal University, Shijiazhuang, China, in 1982, and the Ph.D. degree in operational research and optimal control from the Institute of Applied Mathematics, Chinese Academy of Sciences, Beijing, China, in 1997. From 1997 to 1999, he was a Chinese Post-Doctoral Fellow with Beijing Jiaotong University, Beijing, where he was an Associate Professor in 1999 and has been a Professor in operational research since 2001. He was also a Research Fellow with the City University of Hong Kong, Hong Kong, from 2000 to 2002 and a Visiting Scholar with the University of Waterloo, Waterloo, ON, Canada, from 2006 to 2007.

His current research interests include machine learning, mathematical optimization, mathematics of operations research, and complementarity problems and variational inequalities. Dr. Xiu is the 9-10th Vice President of the Operations Research Society of China, and also serves as a member of Editorial Board for several journals such as Acta Mathematicae Applicatae Sinica, OR Transactions, Operations Research and Management, and Journal of the Operations Research Society of China.
\end{IEEEbiography}

\ifCLASSOPTIONcompsoc


\ifCLASSOPTIONcaptionsoff
  \newpage
\fi

\end{document}

